\documentclass{gtpart}     
\gtart  
\usepackage[english]{babel}
\usepackage{tikz}
\usetikzlibrary{arrows}
\usetikzlibrary{patterns}
\usepackage{caption}
\usepackage{xr, multicol, textcomp,url,hyperref}

\title{A formula for the $\Theta$-invariant from Heegaard diagrams}

\author{Christine Lescop}
\givenname{Christine}
\surname{Lescop}
\address{Institut Fourier, UJF Grenoble, CNRS\\100 rue des maths, BP 74, 38402 Saint-Martin d'H\`eres cedex}
\email{Christine.Lescop@ujf-grenoble.fr}
\urladdr{http://www-fourier.ujf-grenoble.fr/~lescop/}


\arxivreference{arXiv:1209.3219}

\volumenumber{}
\issuenumber{}
\publicationyear{}
\papernumber{}
\startpage{}
\endpage{}
\doi{}
\MR{}
\Zbl{}
\received{}
\revised{}
\accepted{}
\published{}
\publishedonline{}
\proposed{}
\seconded{}
\corresponding{}
\editor{}
\version{}

\def\theequation {\thesection.\arabic{equation}}
\def\hfl#1{\smash{\mathop{\hbox to 10mm{\rightarrowfill}}\limits^{\textstyle
#1}}}
\newtheorem{thm}{Theorem}[section]    
\newtheorem{lemma}[thm]{Lemma}          
 
\newtheorem{proposition}[thm]{Proposition}
\newtheorem{corollary}[thm]{Corollary}
\newtheorem{theorem}[thm]{Theorem}
\theoremstyle{definition}
\newtheorem{example}[thm]{Example}

\newtheorem{remark}[thm]{Remark}

\makeatletter

\newcommand{\tilth}{d_e}

\newcommand{\CI}{{\cal I}}
\newcommand{\CJ}{{\cal J}}

\newcommand{\CN}{{\cal N}}

\newcommand{\halfbra}{{|}}

\newcommand{\ZZ}{\mathbb{Z}}
\newcommand{\RR}{\mathbb{R}}
\newcommand{\QQ}{\mathbb{Q}}
\newcommand{\CC}{\mathbb{C}}

\newcommand{\Asc}{{\cal A}}
\newcommand{\Bdesc}{{\cal B}}
\newcommand{\CaC}{{\cal C}} 
\newcommand{\CD}{{\cal D}}
\newcommand{\fMorse}{f}
\newcommand{\funcf}{f}
\newcommand{\metrigo}{\mathfrak{g}} 
\newcommand{\twocycG}{G}
\newcommand{\Gdu}{G_{\uparrow \downarrow}}

\newcommand{\homotop}{h}
\newcommand{\handlebodA}{H_{\Asc}}
\newcommand{\handlebodB}{H_{\Bdesc}}
\newcommand{\handleboda}{H_a}
\newcommand{\handlebodb}{H_b} 
\newcommand{\lambdh}{\ell_2(\CD)}
\newcommand{\Link}{L}
\newcommand{\linhp}{{\cal L}_+}
\newcommand{\linhn}{{\cal L}_-}
\newcommand{\linh}{{\cal L}}
\newcommand{\lintheta}{{\cal L}_{\theta}}
\newcommand{\matchingo}{\mathfrak{m}}
\newcommand{\normbun}{\mathfrak{V}} 
\newcommand{\prop}{{\cal P}}
\newcommand{\preprop}{P}
\newcommand{\pointp}{p}
\newcommand{\rpt}{\mathbb{RP}^3}
\newcommand{\tang}{T}
\newcommand{\vecv}{\vec{v}}
\newcommand{\manU}{U}
\newcommand{\manV}{V}
\newcommand{\manW}{W} 
\newcommand{\pointw}{w}
\newcommand{\ComX}{X}
\newcommand{\ComY}{Y}

\newcommand{\bp}{\noindent {\sc Proof: }}
\newcommand{\eop}{\nopagebreak \hspace*{\fill}{$\diamond$} \medskip}

\begin{document}

\begin{abstract} The $\Theta$-invariant is the simplest $3$-manifold invariant defined with configuration space integrals. It is actually an invariant of rational homology spheres equipped with a combing over the complement of a point. It can be computed as the algebraic intersection of three propagators associated to a given combing $\ComX$ in the $2$-point configuration space of a $\QQ$--sphere $M$. These propagators represent the linking form of $M$ so that $\Theta(M,\ComX)$ can be thought of as the cube of the linking form of $M$ with respect to the combing $\ComX$.
The invariant $\Theta$ is the sum of $6 \lambda(M)$ and $\frac{p_1(\ComX)}{4}$,
where $\lambda$\index{Nt}{lambda@$\lambda$} denotes the Casson-Walker invariant, and $p_1$\index{Nt}{pwone@$p_1$} is an invariant of combings, which is an extension of a first relative Pontrjagin class.
In this article, we present explicit propagators associated with Heegaard diagrams of a manifold, and we use these ``Morse propagators'', constructed with Greg Kuperberg, to prove a combinatorial formula for the $\Theta$-invariant in terms of Heegaard diagrams.
\end{abstract}
\maketitle

\tableofcontents

\section{Introduction}
\label{secintro}

In this article, a {\em $\QQ$--sphere\/} or {\em rational homology sphere\/} is a smooth closed oriented $3$-manifold that has the same rational homology as $S^3$.

\subsection{General introduction}
\label{subgenint}

The work of Witten \cite{witten} pioneered the introduction of many $\QQ$--sphere invariants. The Le-Murakami-Ohtsuki universal finite type invariant \cite{lmo} and the Kontsevich configuration space invariant \cite{ko}, which was proved to be equivalent to the LMO invariant for integer homology spheres by
G. Kuperberg and D. Thurston \cite{kt}, are among them.
The construction of the Kontsevich configuration space invariant for a $\QQ$--sphere $M$ involves a point $\infty$ in $M$,
an identification of a neighborhood of $\infty$ with a neighborhood of $\infty$ in $S^3=\RR^3 \cup \{\infty\}$, and a parallelization
$\tau$ of $(\check{M}=M \setminus \{\infty\})$ that coincides with the standard parallelization of $\RR^3$ near $\infty$.
The Kontsevich configuration space invariant is in fact an invariant of $(M,\tau)$.
Its degree one part
$\Theta(M,\tau)$ is
the sum of $6 \lambda(M)$ and $\frac{p_1(\tau)}{4}$, where $\lambda$ is the Casson-Walker invariant and $p_1$ is a Pontrjagin number associated with $\tau$, according to a Kuperberg Thurston theorem \cite{kt} generalized to rational homology spheres in \cite{lessumgen}.
Here, the Casson-Walker invariant $\lambda$ is normalized as in \cite{akmc,gm,mar} for integer homology spheres, and like $\frac{1}{2}\lambda_W$ for rational homology spheres where $\lambda_W$ is the Walker normalisation in \cite{wal}.

The invariant $\Theta(M,\tau)$ reads 
$$\Theta(M,\tau)=\int_{\check{M}^2\setminus \mbox{diag}(\check{M})^2} \omega(M,\tau)^3$$
for some closed $2$-form $\omega(M,\tau)$, which is often called a {\em propagator\/}.
As it is developed in \cite[Section 6.5]{lessumgen}, $\Theta(M,\tau)$ can also be written as the algebraic intersection
of three $4$-dimensional chains in a compactification $C_2(M)$\index{Nt}{Ctwo@$C_2(M)$} of $\check{M}^2\setminus \mbox{diag}(\check{M})^2$,
for chains that are Poincar\'e dual to $\omega(M,\tau)$ in the $6$--dimensional configuration space $C_2(M)$.
In this article, a {\em propagator\/} will be such a $4$-chain. For more precise definitions, see Subsection~\ref{subprop}.
A \emph{combing} of a $3$-manifold $M$ as above is an asymptotically constant nowhere zero section of the tangent bundle to $\check{M}$.

In Theorem~\ref{thmdefinvcomb}, we will prove that the invariant $\Theta$ is an invariant of combed $\QQ$--spheres $(M,X)$ rather than an invariant of parallelised punctured $\QQ$--spheres, so that $(4 \Theta(M,X)-24\lambda(M))$ is an extension of the Pontrjagin number $p_1$ to combings. The invariant $p_1$ of parallelizations coincides with the Hirzebruch defect of the parallelization $\tau$ studied in \cite{hirzebruchEM, km}.
This invariant $p_1$ of combings is studied in \cite{lescomb}, and it is shown to be the analogue of the Gompf $\theta$-invariant \cite[Section 4]{Go} of 
$\QQ$--sphere combings, for asymptotically constant combings of punctured $\QQ$--spheres. The variations of $\Theta$, $\theta$ and $p_1$ under various combing changes are described in \cite{lescomb}.

In Section~\ref{secpropout}, we describe explicit propagators associated with Morse functions or with Heegaard splittings. These ``Morse propagators'' have been obtained in collaboration with Greg Kuperberg.
Then we use these propagators to produce a combinatorial description of $\Theta$ in terms of Heegaard diagrams in Theorem~\ref{thmmain}.

Our Morse propagators and our techniques could be applied to compute more configuration space invariants, and they might be useful to relate finite type invariants to Heegaard Floer homology.

This article benefited from the stimulating visit of Greg Kuperberg in Grenoble in 2010-2011. It also benefited from the referees' comments.

\subsection{Conventions and notations}
Unless otherwise mentioned, all manifolds are oriented. Boundaries are oriented by the outward normal first convention.
Products are oriented by the order of the factors. More generally, unless otherwise mentioned, the order of appearance of coordinates or parameters orients manifolds or chains, which are linear combinations of manifolds.
The fiber of the normal bundle $\normbun(\manV)$ to an oriented submanifold $\manV$ is oriented so that the normal bundle followed by the tangent bundle to the submanifold induce the orientation of the ambient manifold, fiberwise.
The transverse intersection of two submanifolds $\manV$ and $\manW$ is oriented so that the normal bundle to $\manV\cap \manW$ is $(\normbun(\manV) \oplus \normbun(\manW))$, fiberwise.
When the dimensions of two such submanifolds add up to the dimension of the ambient manifold $\manU$, each intersection point $x$ is equipped with a sign $\pm 1$ that is $1$ if and only if $(\normbun_x(\manV) \oplus \normbun_x(\manW))$ (or equivalently $(T_x(\manV) \oplus T_x(\manW))$) induces the orientation of $\manU$. When $\manV$ is compact, the sum of the signs of the intersection points is the \emph{algebraic intersection number} $\langle \manV,\manW \rangle_{\manU}$.
 For a manifold $\manV$, $(-\manV)$ denotes the manifold $\manV$ equipped with the opposite orientation.

\section{The \texorpdfstring{$\Theta$}{Theta}-invariant}
\setcounter{equation}{0}
This section presents a complete definition of the invariant $\Theta$.

\subsection{On configuration spaces}
\label{subconf}

In this article, {\em blowing up\/} a submanifold $\manV$ means replacing it by its unit normal bundle. Locally, $\RR^c \times \manV$ is replaced with $[0,\infty[ \times S^{c-1} \times \manV$, where the fiber $\RR^c$ of the normal bundle is naturally identified with $\{0\} \cup \left(]0,\infty[ \times S^{c-1}\right)$. Topologically, this amounts to removing an open tubular neighborhood of the submanifold (thought of as infinitely small), but the process is canonical, so that the created boundary is the unit normal bundle to the submanifold and there is a canonical projection from the manifold obtained by blow-up to the initial manifold.

In a closed $3$-manifold $M$, we fix a point $\infty$ and define the blown-up manifold
$C_1(M)$ as the compact $3$-manifold obtained from $M$ by blowing up $\{\infty\}$. This space $C_1(M)$ is a compactification of $\check{M}=(M \setminus \{\infty\})$.

The {\em configuration space\/} $C_2(M)$\index{Nt}{Ctwo@$C_2(M)$} is the compact $6$--manifold with boundary and corners obtained from $M^2$ by blowing up $(\infty,\infty)$, and the closures of $\{\infty\} \times \check{M}$, $\check{M} \times \{\infty\}$ and the diagonal of $\check{M}^2$, successively.

Then the boundary $\partial C_2(M)$ of $C_2(M)$ contains the unit normal bundle to the diagonal of $\check{M}^2$. This bundle is canonically isomorphic to the unit tangent bundle $U\check{M}$\index{Nt}{UM@$U\check{M}$} of $\check{M}$ via the map
$$[(x,y)] \in \frac{\frac{T_m\check{M}^2}{\mbox{\tiny diag}} \setminus \{0\}}{\RR^{+\ast}} \mapsto [y-x] \in \frac{T_m\check{M} \setminus \{0\}}{\RR^{+\ast}}.$$

When $M$ is a rational homology sphere, the configuration space $C_2(M)$ has the same rational homology as $S^2$ (see the proof of Theorem~\ref{thmdefinvcomb} below) and $H_2(C_2(M);\QQ)$ has a canonical generator $[S]$\index{Nt}{S@$[S]$}
that is the homology class of a product $(x \times \partial B(x))$ where $B(x)$ is a ball embedded in $\check{M}$ that contains $x$ in its interior.
For a $2$-component link $(J,K)$ of $M$, the homology class $[J \times K]$ of
$J \times K$ in $H_2(C_2(M);\QQ)$ reads $lk(J,K)[S]$, where $lk(J,K)$ is the \emph{linking number} of $J$ and $K$, which is the algebraic intersection number of $J$ and a $2$-dimensional chain bounded by $K$ (see \cite[Proposition 1.6]{lesmek}).

\subsection{On propagators}
\label{subprop}

When $M$ is a rational homology sphere, a {\em propagator\/} of $C_2(M)$ is a $4$--cycle $\prop$ of $(C_2(M),\partial C_2(M))$ that is Poincar\'e dual to the preferred generator of
$H^2(C_2(M);\QQ)$ that maps $[S]$ to $1$.
For such a propagator $\prop$, for any $2$-cycle $\twocycG$ of $C_2(M)$,
$$[\twocycG]=\langle \prop,\twocycG\rangle_{C_2(M)} [S]$$
in $H_2(C_2(M);\QQ)$ where $\langle \prop,\twocycG\rangle_{C_2(M)}$ denotes
the algebraic intersection of $\prop$ and $\twocycG$ in $C_2(M)$.

Let $B$ and $\frac12 B$ be two balls in $\RR^3$ of respective radii $R$ and $\frac{R}2$, centered at the origin in $\RR^3$.
Identify a neighborhood of $\infty$ in $M$ with $S^3 \setminus (\frac12 B)$ in $(S^3=\RR^3 \cup \{ \infty\})$
so that $\check{M}$ reads $\check{M}=B_M \cup_{]R/2,R] \times S^2} (\RR^3 \setminus (\frac12 B))$ for a rational homology ball $B_M$\index{Nt}{BM@$B_M$} whose complement in $\check{M}$ is identified with $\RR^3 \setminus B$.
There is a canonical regular map
\index{Nt}{pzinfty@$p_{\infty}$}$$p_{\infty} \colon (\partial C_2(M) \setminus UB_M) \rightarrow S^2$$
that maps the limit in $\partial C_2(M)$ of a convergent sequence of ordered pairs of distinct points of $\left(\check{M}\setminus B_M\right)^2$ to the limit of the direction from the first point to the second one. See \cite[Lemma 1.1]{lesconst}.
Let $$\tau_s\colon \RR^3 \times \RR^3 \rightarrow T\RR^3$$ denote the standard parallelization of $\RR^3$.
In this article, a {\em combing} $\ComX$ of a $\QQ$--sphere $M$ is a section of $U\check{M}$ that is constant outside $B_M$, i.e. that reads $\tau_s((\check{M} \setminus B_M) \times \{\vecv(\ComX)\})$ for some fixed $\vecv(\ComX) \in S^2$ outside $B_M$.
The {\em propagator boundary\/} $\partial \prop_{\ComX}$ associated with such a combing $\ComX$ is the following $3$--cycle of $\partial C_2(M)$
$$\partial \prop_{\ComX}=p_{\infty}^{-1}(\vecv(\ComX)) \cup \ComX(B_M)$$
where the part
$\ComX(B_M)$ of $\partial C_2(M)$ is the graph of the restriction of the combing $\ComX$ to $B_M$
and a {\em propagator associated with the combing $\ComX$\/} is a $4$--chain $\prop_{\ComX}$ of $C_2(M)$ whose boundary reads
$\partial \prop_{\ComX}$. Such a $\prop_{\ComX}$ is indeed a propagator (because for a tiny sphere $\partial B(x)$ around a point $x$, $\langle x \times \partial B(x),\prop_{\ComX} \rangle_{C_2(M)}$ is the algebraic intersection in $U\check{M}$ of a fiber and the section $\ComX(\check{M})$, which is one).

\subsection{On the \texorpdfstring{$\Theta$}{Theta}-invariant of a combed \texorpdfstring{$\QQ$}{Q}--sphere}
\label{subTheta}

\begin{theorem}\label{thmdefinvcomb}
Let $\ComX$ be a combing of a rational homology sphere $M$, and let $(-\ComX)$ be the opposite combing.
Let $\prop_{\ComX}$ and $\prop_{-\ComX}$ be two associated transverse propagators. Then $\prop_{\ComX} \cap \prop_{-\ComX}$ is a two-dimensional cycle whose homology class is independent of the chosen propagators. It reads
$\Theta(M,\ComX)[S]$, where $\Theta(M,\ComX)$ is therefore a rational valued topological invariant of $M$ and of the homotopy class of $\ComX$. 

\end{theorem}
\bp Let us first show that $C_2(M)$ has the same rational homology as $S^2$.
The space $C_2(M)$ is homotopy equivalent to $(\check{M}^2 \setminus \mbox{diag})$. Since $\check{M}$ is a rational homology $\RR^3$, the rational homology of $(\check{M}^2 \setminus \mbox{diag})$ is isomorphic to the rational homology of $((\RR^3)^2 \setminus \mbox{diag})$.
Since $((\RR^3)^2 \setminus \mbox{diag})$ is homeomorphic to $\RR^3 \times ]0,\infty[ \times S^2$ via the map $$(x,y) \mapsto (x,\parallel y-x \parallel, \frac{1}{\parallel y-x \parallel}(y-x)),$$ $((\RR^3)^2 \setminus \mbox{diag})$ is homotopy equivalent to $S^2$. 

In particular, since $H_3(C_2(M);\QQ)=0$, there exist propagators $\prop_{\ComX}$ and $\prop_{-\ComX}$ with the given boundaries $\partial \prop_{\ComX}$ and $\partial \prop_{-\ComX}$. By general position arguments \cite[Chapter 3]{hirsch}, $\prop_{\ComX}$ and $\prop_{-\ComX}$ can be assumed to be transverse. (Explicit transverse propagators will be constructed in Subsections~\ref{subdeformprop} and \ref{subpertudif}.)
Without loss, assume that $\prop_{\pm \ComX} \cap \partial C_2(M)=\partial \prop_{\pm \ComX}$. 
Since $\partial \prop_{\ComX}$ and $\partial \prop_{-\ComX}$ do not intersect, $\prop_{\ComX} \cap \prop_{-\ComX}$ is a $2$--cycle. Since $H_4(C_2(M);\QQ)=0$, the homology class of $\prop_{\ComX} \cap \prop_{-\ComX}$ in $H_2(C_2(M);\QQ)$ does not depend on the choices of $\prop_{\ComX}$ and $\prop_{-\ComX}$ with their given boundaries. Then it is easy to see that $\Theta(M,\ComX) \in \QQ$ is a locally constant function of the combing $\ComX$.
\eop

When $M$ is an integer homology sphere, a combing $\ComX$ is the first vector of a unique parallelization $\tau(\ComX)$ that coincides with $\tau_s$ outside $B_M$, up to homotopy. When $M$ is a rational homology sphere, and when
$\ComX$ is the first vector of a such a parallelization $\tau(\ComX)$, this parallelization is again unique.
In this case, the invariant $\Theta(M,\ComX)$ can be identified with the invariant $\Theta(M,\tau(\ComX))$ discussed in Subsection~\ref{subgenint} using \cite[Lemma 6.16]{lessumgen}.

Let $W$ be a connected compact $4$--dimensional manifold with corners with signature $0$ whose boundary is
$$\partial W = B_M \cup_{1\times \partial B_M} (-[0,1]\times S^2) \cup_{0 \times S^2} (-B^3)$$
and that is identified with an open subspace of one of the products $[0,1[ \times B^3$ or $]0,1] \times B_M$ near $\partial W$.
Then the Pontrjagin number $p_1(\tau(\ComX))$\index{Nt}{pwone@$p_1$} is the obstruction to extending the trivialization of $TW \otimes \CC$ induced by $\tau(\ComX)$ and $\tau_s$ on $\partial W$ to $W$. This obstruction lives in $H^4(W,\partial W;\pi_3(SU(4))=\ZZ)=\ZZ$. See \cite[Section 1.5]{lesconst} for more details.
In \cite{kt},
G. Kuperberg and D. Thurston proved that
$$\Theta(M,\ComX)=6 \lambda(M) + \frac{p_1(\tau(\ComX))}{4}$$
when $M$ is an integer homology sphere. This result was extended to $\QQ$--spheres by the author in \cite[Theorem 2.6 and Section 6.5]{lessumgen}.
Setting $p_1(\ComX)=(4 \Theta(M,\ComX)-24\lambda(M))$ extends the Pontrjagin number from parallelizations to combings so that the formula above is still valid for combings.

The following theorem is proved in \cite{lescomb}.

\begin{theorem}
\label{thmvarinvcomb}
Let $\ComX$ and $\ComY$ be two combings of $M$ such that the cycle $\partial \prop_{\ComY}$ is transverse to $\partial \prop_{\ComX}$
and to $\partial \prop_{-\ComX}$ in $\partial C_2(M)$.
Then the oriented intersection $\partial \prop_{\ComX} \cap \partial \prop_{\ComY}$ (resp. $\partial \prop_{\ComX} \cap \partial \prop_{-\ComY}$) is the graph of the restriction of $\ComX$ to an oriented link $L_{\ComX=\ComY}$ (resp. $L_{\ComX=-\ComY}$) in $U \check{M}$ and
$$\Theta(M,\ComY)-\Theta(M,\ComX)= \frac{p_1(\ComY)-p_1(\ComX)}{4}= lk(L_{\ComX=\ComY},L_{\ComX=-\ComY}).$$
\end{theorem}

\section{The formula for the \texorpdfstring{$\Theta$}{Theta}-invariant from Heegaard diagrams}
\setcounter{equation}{0}
\label{secformula}

\subsection{On Heegaard diagrams}
\label{subonhd}
Every closed $3$--manifold $M$ can be written as the union of two handlebodies $\handlebodA$ and $\handlebodB$ glued along their common boundary, which is a genus $g$ surface as
$$M = \handlebodA \cup_{\partial \handlebodA} \handlebodB$$
where $\partial \handlebodA= -\partial \handlebodB$.
Such a decomposition is called a {\em Heegaard decomposition\/} or  a \emph{Heegaard splitting} of $M$.
A {\em system of meridian disks\/} for $\handlebodA$ is a system of $g$ disjoint disks $D(\alpha_i)$\index{Nt}{Dalpha@$D(\alpha_i)$} properly embedded in $\handlebodA$ such that the union of the boundaries $\alpha_i$ of the $D(\alpha_i)$ does not separate $\partial \handlebodA$.
Let $(D(\alpha_i))_{i\in \{1,\dots,g\}}$ be such a system for $\handlebodA$ and let $(D(\beta_j))_{j\in \{1,\dots,g\}}$ be such a system for $\handlebodB$. Then the surface equipped with the collections of the curves $\alpha_i$ and the curves $\beta_j=\partial D(\beta_j)$ determines $M$.
When the collections $(\alpha_i)_{i\in \{1,\dots,g\}}$ and
$(\beta_j)_{j\in \{1,\dots,g\}}$ are transverse, the data collection
$$\CD=(\partial \handlebodA, (\alpha_i)_{i\in \{1,\dots,g\}}, (\beta_j)_{j\in \{1,\dots,g\}})$$ is called a {\em genus g Heegaard diagram.\/}
Figure~\ref{figheegrpt} shows two Heegaard diagrams of $\rpt$ (or $SO(3)$).

\begin{figure}\begin{center}
\begin{tikzpicture} \useasboundingbox (.2,-1) rectangle (12,1); 
\draw [red,->]  (1.88,.45) node[left]{\scriptsize $\alpha_1$}  (2,.9) .. controls (1.9,.9) and (1.75,.7) .. (1.75,.5);  
\draw  [red] (1.75,.5) .. controls (1.75,.3) and (1.9,.1) .. (2,.1);
\draw [red,dashed] (2,.9) .. controls (2.1,.9) and (2.25,.7) .. (2.25,.5) .. controls (2.25,.3) and (2.1,.1) .. (2,.1);
\draw plot[smooth] coordinates{(1.4,.1) (1.6,0) (2,-.1) (2.4,0) (2.6,.1)};
\draw plot[smooth] coordinates{(1.6,0) (2,.1) (2.4,0)};
\fill (.9,-.4) circle (0.05);
\draw (.85,-.25) node{\scriptsize $\pointw$};
\draw (.6,0) .. controls (.6,.45) and (1.3,.9)  .. (2,.9) .. controls (2.7,.9) and (3.4,.45) .. (3.4,0) .. controls  (3.4,-.45) and (2.7,-.9) .. (2,-.9) .. controls (1.3,-.9) and (.6,-.45)  ..  (.6,0);
\draw [blue,->]  (.6,0) .. controls (.6,.1) and (1.2,.7)  .. (2,.7) .. controls (2.6,.7) and (3.15,.35) .. (3.15,0);
\draw  [blue,->]    (2.98,0) node{\scriptsize $\beta_1$} (3.15,0) .. controls  (3.15,-.35) and (2.6,-.7) .. (2,-.7) .. controls (1.4,-.7) and (1,-.35)  ..  (1,-.1)  .. controls (1,.1) and (1.5,.3)  .. (2,.3) .. controls (2.3,.3) and (2.8,.15) .. (2.8,0) ;
\draw  [blue] (2.8,0) .. controls  (2.8,-.1) and (2,-.2) .. (2,-.1);
\draw [blue,dashed] (2,-.1)  .. controls (2,-.2) and (.6,-.1) ..(.6,0);
\draw  (3.4,-.7) node{\scriptsize $\CD_1$}  (1.75,.78) node{\scriptsize $d$}   (1.9,.4) node{\scriptsize $c$};
\begin{scope}[xshift=5cm] 
\draw [red,->]  (1.88,.45) node[left]{\scriptsize $\alpha_1$}  (2,.9) .. controls (1.9,.9) and (1.75,.7) .. (1.75,.5);  
\draw  [red] (1.75,.5) .. controls (1.75,.3) and (1.9,.1) .. (2,.1);
\draw [red,dashed] (2,.9) .. controls (2.1,.9) and (2.25,.7) .. (2.25,.5) .. controls (2.25,.3) and (2.1,.1) .. (2,.1);
\draw plot[smooth] coordinates{(1.4,.1) (1.6,0) (2,-.1) (2.4,0) (2.6,.1)};
\draw plot[smooth] coordinates{(1.6,0) (2,.1) (2.4,0)};
\fill (.9,-.4) circle (0.05);
\draw (.85,-.25) node{\scriptsize $\pointw$};
\draw (.6,0) .. controls (.6,.45) and (1.3,.9)  .. (2,.9) .. controls (2.7,.9) and (3.4,.3) .. (3.6,.3) (3.6,-.3) .. controls  (3.4,-.3) and (2.7,-.9) .. (2,-.9) .. controls (1.3,-.9) and (.6,-.45)  ..  (.6,0);
\draw [blue,->]  (.6,0) .. controls (.6,.1) and (1.2,.7)  .. (2,.7) .. controls (2.6,.7) and (3.15,.15) .. (3.6,.15) ;
\draw [blue,->]  (3.6,.15) .. controls (4.05,.15) and (4.7,.7) .. (5.2,.7) .. controls (5.7,.7) and (6.4,.4) .. (6.4,0) ;
\draw  [blue,->]  (6.4,0) .. controls (6.4,-.4) and  (5.7,-.7) .. (5.2,-.7) ..  controls (4.7,-.7) and (4.05,-.15) .. (3.6,-.15);
\draw  [blue,->]    (3.5,0) node{\scriptsize $\beta_1$} (3.6,-0.15) .. controls  (3.15,-.15) and (2.6,-.7) .. (2,-.7) .. controls (1.4,-.7) and (1,-.35)  ..  (1,-.1)  .. controls (1,.1) and (1.5,.3)  .. (2,.3) .. controls (2.3,.3) and (2.8,.15) .. (2.8,0) ;
\draw  [blue] (2.8,0) .. controls  (2.8,-.1) and (2,-.2) .. (2,-.1);
\draw [blue,dashed] (2,-.1)  .. controls (2,-.2) and (.6,-.1) ..(.6,0);
\draw    (1.75,.78) node{\scriptsize $d$}   (1.9,.4) node{\scriptsize $c$};
\draw  (3.6,-.7) node{\scriptsize $\CD_2$};
\end{scope}
\begin{scope}[xshift=8.2cm] 
\draw (.4,.3) .. controls (.6,.3) and (1.3,.9)  .. (2,.9) .. controls (2.7,.9) and (3.4,.45) .. (3.4,0) .. controls  (3.4,-.45) and (2.7,-.9) .. (2,-.9) .. controls (1.3,-.9) and (.6,-.3)  ..  (.4,-.3);
\draw [red,->]  (1.88,.45) node[left]{\scriptsize $\alpha_2$}  (2,.9) .. controls (1.9,.9) and (1.75,.7) .. (1.75,.5);  
\draw  [red] (1.75,.5) .. controls (1.75,.3) and (1.9,.1) .. (2,.1);
\draw [red,dashed] (2,.9) .. controls (2.1,.9) and (2.25,.7) .. (2.25,.5) .. controls (2.25,.3) and (2.1,.1) .. (2,.1);
\draw plot[smooth] coordinates{(1.4,.1) (1.6,0) (2,-.1) (2.4,0) (2.6,.1)};
\draw plot[smooth] coordinates{(1.6,0) (2,.1) (2.4,0)};
\draw [blue,->] (2,-.25) node{\scriptsize $\beta_2$} (2,-.4) .. controls (1.6,-.4) and (1.2,-.2)  ..  (1.2,0)  .. controls (1.2,.2) and (1.6,.4)  .. (2,.4) .. controls (2.4,.4) and (2.8,.2) .. (2.8,0) .. controls (2.8,-.2) and (2.4,-.4)  .. (2,-.4);
\draw (1.75,.76) node{\scriptsize $f$}   (1.72,.22) node{\scriptsize $e$};
\end{scope}
\end{tikzpicture}
\caption{Two Heegaard diagrams of $\rpt$}
\label{figheegrpt}\end{center}
\end{figure}

We fix a genus g Heegaard diagram $\CD$. A {\em crossing \/} $c$ of $\CD$ is an intersection point of a curve $\alpha_{i(c)}$ and a curve $\beta_{j(c)}$. Its sign $\sigma(c)$\index{Nt}{ssigma@$\sigma(c)$} is $1$ if $\partial \handlebodA$ is oriented by the oriented tangent vector of $\alpha_{i(c)}$ followed by the oriented tangent vector of $\beta_{j(c)}$ at $c$. It is $(-1)$ otherwise.
The collection of crossings of $\CD$ is denoted by $\CaC$.

Fix a point $a_i$ inside each disk $D(\alpha_i)$ and a point $b_j$ inside each disk $D(\beta_j)$.
Then join $a_i$ to each crossing $c$ of $\alpha_i$ by a segment
$[a_i,c]_{D(\alpha_i)}$ oriented from $a_i$ to $c$ in $D(\alpha_i)$, so that these segments only meet at $a_i$ for different $c$. Similarly define segments $[c,b_{j(c)}]_{D(\beta_{j(c)})}$ from $c$ to $b_{j(c)}$ in $D(\beta_{j(c)})$.
Then for each $c$, define the {\em flow line \/}
$\gamma(c) = [a_{i(c)},c]_{D(\alpha_{i(c)})} \cup [c,b_{j(c)}]_{D(\beta_{j(c)})}$.

A Heegaard decomposition as above can be obtained from a Morse function $\funcf_M$ on $M$ with one minimum, one maximum, index one-critical points $a_i$ mapped to $1$ and index $2$ critical points $b_j$ mapped to $5$, by setting $\handlebodA=\funcf_M^{-1}(]-\infty,3])$ and $\handlebodB=\funcf_M^{-1}([3,+\infty[)$ \cite[Chapter 6]{hirsch}. For an appropriate (generic) metric, the descending manifolds of the $b_j$ intersect $\handlebodB$ as disks $D(\beta_{j})$
and the ascending manifolds of the $a_i$ intersect $\handlebodA$ as disks $D(\alpha_{i})$ so that the boundaries $\alpha_i$ of the $D(\alpha_{i})$ are transverse to  the boundaries $\beta_j$ of the $D(\beta_{j})$. The Morse function $\funcf_M$ and such a metric $\metrigo$ induce a Heegaard diagram of $M$ where the flow line $\gamma(c)$\index{Nt}{gamma(c)@$\gamma(c)$} above can be chosen as the closure of the actual flow line through $c$ for the gradient flow of $\funcf_M$. Conversely, for any Heegaard diagram, there exist a Morse function and a metric as above that produce this diagram.

An \emph{exterior point} of the diagram is a point of $\partial \handlebodA \setminus \left(\coprod_{i=1}^g \alpha_i \cup \coprod_{j=1}^g \beta_j\right)$ as in Figure~\ref{figheegrpt}.
Pick an exterior point $\pointw$\index{Nt}{w@$\pointw$} of the diagram, and let $\overline{\gamma(\pointw)}$ be the closure of the flow line through $w$ with respect to $\metrigo$. It goes from the minimum of $\funcf_M$ to its maximum. Identify a ball around $\overline{\gamma(w)}$ with a neighborhood of $\infty$ in $S^3$, so that the restriction of $\funcf_M$ to $B_M$ extends to $\check{M}$ as a Morse function $\funcf$ that is the standard height function outside $B_M$, that has no extremum, whose index one critical points $a_i$\index{Nt}{ai@$a_i$} are mapped to $1$, and whose index $2$ critical points $b_j$ are mapped to $5$. 

In Section~\ref{secpropout}, we describe an explicit propagator $\prop(\funcf,\metrigo)$ associated with a Morse function $\funcf$ of $\check{M}$ that satisfies these properties, and with a metric $\metrigo$ that is standard outside $B_M$.

A {\em matching\/} in a genus $g$ Heegaard diagram $(\partial \handlebodA,\{\alpha_i\}_{i=1,\dots,g},\{\beta_j\}_{j=1,\dots,g})$ is a set $\matchingo$\index{Nt}{m@$\matchingo$} of $g$ crossings such that every curve of the diagram contains one crossing of $\matchingo$. Thus a matching $\matchingo$ can be written as
$\matchingo=\{c_i; i \in \{1,2, \dots, g\}\}$
where the $c_i$ are crossings of $\alpha_i \cap \beta_{\rho^{-1}(i)}$ for a permutation $\rho$ of $\{1,2, \dots, g\}$.

The choice of a matching $\matchingo$ and of an exterior point $\pointw$ in a diagram $\CD$ of $M$ equips $\check{M}$ with a combing $\ComX(\pointw,\matchingo )=\ComX(\CD,\pointw,\matchingo )$, which is roughly obtained from the gradient vector of $\funcf$ by reversing this singular field along the flow lines through the points of $\matchingo$. The combing $\ComX(\pointw,\matchingo)$ of $\check{M}$ is precisely described in Subsection~\ref{subframedsel}. \footnote{The same data $(\CD,\pointw,\matchingo)$ can be used to define an Euler structure or a combing of the non-punctured $M$. Such a combing represents a Spin$^c$ structure. Matchings representing a given Spin$^c$-structure $\xi$ are the generators of a chain complex whose homology is a Heegaard-Floer homology of $(M,\xi)$.}
The propagator $\prop(\funcf,\metrigo)$ is modified near $\partial C_2(M)$ to become a propagator $\prop_{\ComX(\pointw,\matchingo)}$ associated with $\ComX(\pointw,\matchingo)$ in Subsection~\ref{subFcomb}.

Sections~\ref{seccompint} and \ref{secproof} are devoted to the computation of $\Theta(M,\ComX(\pointw,\matchingo ))$, performed by evaluating the homology class of the intersection of $\prop_{\ComX(\pointw,\matchingo)}$ and $\prop_{-\ComX(\pointw,\matchingo)}$, and by applying the definition of Theorem~\ref{thmdefinvcomb}.
The current section is devoted to presenting the combinatorial formula 
$$\Theta(M,\ComX(\CD,\pointw,\matchingo))=\lambdh + lk(\Link(\CD,\matchingo),\Link(\CD,\matchingo)_{\parallel}) -e(\CD,\pointw,\matchingo)$$
that we get from our computation.

The three ingredients of our formula are completely combinatorial. They can be read on the Heegaard diagram without referring to Morse functions. 
However, they also have a topological meaning, which explains the chosen notation and which makes them easier to apprehend. We first introduce the ingredients $lk(\Link(\CD,\matchingo),\Link(\CD,\matchingo)_{\parallel})$ and $\lambdh$ with their topological interpretations in Subsections~\ref{subparflow} and \ref{subfundcyc}, respectively, before giving their combinatorial expressions in Corollary~\ref{corcombexp} at the end of Subsection~\ref{subeval}. The combinatorial definition of $e(\CD,\pointw,\matchingo)$ is given in Subsection~\ref{subframedHd}.

Let $$[\CJ_{ji}]_{(j,i) \in \{1,\dots,g\}^2}=[\langle \alpha_i,\beta_j \rangle_{\partial \handlebodA}]^{-1}$$\index{Nt}{Jcal@$\CJ_{ji}$}
be the inverse matrix of the matrix of the algebraic intersection numbers $\langle \alpha_i,\beta_j \rangle_{\partial \handlebodA}$. 
$$\sum_{i=1}^g\CJ_{ji}\langle \alpha_i,\beta_k \rangle_{\partial \handlebodA}=\delta_{jk}=\left\{\begin{array}{ll} 1 \;& \mbox{if} \; j=k\\ 0 & \mbox{otherwise.}\end{array}\right.$$

Let $$\Link(\matchingo) = \Link(\CD,\matchingo)=\sum_{i=1}^g \gamma(c_i)- \sum_{c\in \CaC} \CJ_{j(c)i(c)} \sigma(c) \gamma(c).$$\index{Nt}{Lm@$\Link(\matchingo)= \Link(\CD,\matchingo)$}

Note that $\Link(\matchingo)$
is a cycle since
$$ \partial \Link(\matchingo)=\sum_{i=1}^g(b_i-a_i)- \sum_{(i,j) \in \{1,\dots,g\}^2} \CJ_{ji}\langle \alpha_i,\beta_j \rangle_{\partial \handlebodA} (b_j-a_i)=0.$$

The term $lk(\Link(\CD,\matchingo),\Link(\CD,\matchingo)_{\parallel})$ is the linking number of $\Link(\matchingo)$ with a canonical parallel $\Link(\matchingo)_{\parallel}$ of $\Link(\matchingo)$ that is defined in Subsection~\ref{subparflow} below.

\begin{example}
\label{exaheegone}
For the genus one Heegaard diagram $\CD_1$ of Figure~\ref{figheegrpt}, $\sigma(c)=1$, $\langle \alpha_1,\beta_1 \rangle_{\partial \handlebodA}=2$, $\CJ_{11}=\frac12$, we choose $\{c\}$ as a matching and $\Link(\{c\})=\frac12(\gamma(c)-\gamma(d))$.
 
For the genus two Heegaard diagram $\CD_2$ of Figure~\ref{figheegrpt}, $\langle \alpha_2,\beta_1 \rangle_{\partial \handlebodA}=1$, $\CJ_{11}=\frac12$, $\CJ_{22}=1$, $\CJ_{12}=0$, $\CJ_{21}=-\frac12$,  we choose the matching $\{c,e\}$ and $\Link(\{c,e\})=\frac12(\gamma(c)-\gamma(d))$.
\end{example}

\subsection{Parallels of flow lines}
\label{subparflow}

For a crossing $c \in \alpha_{i(c)} \cap \beta_{j(c)}$, $\gamma(c)_{\parallel}$ will denote the following chain. Consider a small meridian curve $m(c)$ of $\gamma(c)$ on $\partial \handlebodA$, it intersects $\beta_{j(c)}$ at two points: $c_{\Asc}^+$ on the positive side of $D(\alpha_{i(c)})$ and $c_{\Asc}^-$ on the negative side of $D(\alpha_{i(c)})$. The meridian $m(c)$ also intersects $\alpha_{i(c)}$ at $c_{\Bdesc}^+$ on the positive side of $D(\beta_{j(c)})$ and $c_{\Bdesc}^-$ on the negative side of $D(\beta_{j(c)})$.
Let $[c_{\Asc}^+,c_{\Bdesc}^+]$, $[c_{\Asc}^+,c_{\Bdesc}^-]$, $[c_{\Asc}^-,c_{\Bdesc}^+]$ and $[c_{\Asc}^-,c_{\Bdesc}^-]$ denote the four quarters of $m(c)$ with the natural ends and orientations associated with the notation, as in Figure~\ref{figgammapar}.

\begin{figure}[h]
\begin{center}
\begin{tikzpicture}
\useasboundingbox (-3.7,-1.4) rectangle (3.7,1.4);
\begin{scope}[xshift=-3cm]
\draw [blue,->] (-.2,1.3) node{\scriptsize $\beta_j$} (0,-1.2) -- (0,1.2);
\draw [red,->] (1.1,.1) node[right]{\scriptsize $\alpha_i$} (-1.2,0) -- (1.2,0);
\draw (0,0) circle (.8);
\draw [->] (0,-.8) arc (-90:-135:.8);
\draw [->] (0,-.8) arc (-90:-45:.8);
\draw [->] (0,.8) arc (90:135:.8);
\draw [->] (0,.8) arc (90:45:.8);
\draw (.2,.2) node{\scriptsize $c$} (-1.8,1.3) node[left]{\scriptsize $\sigma(c)=1$}  (-1,-.25) node{\scriptsize $c_{\Bdesc}^-$} (1.05,-.25) node{\scriptsize $c_{\Bdesc}^+$} (.3,-1) node{\scriptsize $c_{\Asc}^-$} (-.45,-.8) node[left]{\tiny $[c_{\Asc}^-,c_{\Bdesc}^-]$} (-.45,.8) node[left]{\tiny $[c_{\Asc}^+,c_{\Bdesc}^-]$} (.45,-.8) node[right]{\tiny $[c_{\Asc}^-,c_{\Bdesc}^+]$} (.45,.8) node[right]{\tiny $[c_{\Asc}^+,c_{\Bdesc}^+]$} (.3,1) node{\scriptsize $c_{\Asc}^+$};
\fill (0,-.8) circle (0.05) (0,.8) circle (0.05) (.8,0) circle (0.05) (-.8,0) circle (0.05) (0,0) circle (0.05);
\end{scope}
\begin{scope}[xshift=3cm]
\draw [blue,->] (-.2,1.3) node{\scriptsize $\beta_j$} (0,1.2) -- (0,-1.2);
\draw [red,->] (1.1,.1) node[right]{\scriptsize $\alpha_i$} (-1.2,0) -- (1.2,0);
\draw (0,0) circle (.8);
\draw [->] (0,-.8) arc (-90:-135:.8);
\draw [->] (0,-.8) arc (-90:-45:.8);
\draw [->] (0,.8) arc (90:135:.8);
\draw [->] (0,.8) arc (90:45:.8);
\draw (.2,.2) node{\scriptsize $c$} (1.8,1.3) node[right]{\scriptsize $\sigma(c)=-1$} (-1,-.25) node{\scriptsize $c_{\Bdesc}^+$} (1.05,-.25) node{\scriptsize $c_{\Bdesc}^-$} (.3,-1) node{\scriptsize $c_{\Asc}^-$} (-.45,-.8) node[left]{\tiny  $[c_{\Asc}^-,c_{\Bdesc}^+]$} (-.45,.8) node[left]{\tiny $[c_{\Asc}^+,c_{\Bdesc}^+]$} (.45,-.8) node[right]{\tiny $[c_{\Asc}^-,c_{\Bdesc}^-]$} (.45,.8) node[right]{\tiny $[c_{\Asc}^+,c_{\Bdesc}^-]$} (.3,1) node{\scriptsize $c_{\Asc}^+$};
\fill (0,-.8) circle (0.05) (0,.8) circle (0.05) (.8,0) circle (0.05) (-.8,0) circle (0.05) (0,0) circle (0.05);
\end{scope}
\end{tikzpicture}
\caption{$m(c)$, $c_{\Asc}^+$, $c_{\Asc}^-$, $c_{\Bdesc}^+$ and $c_{\Bdesc}^-$}
\label{figgammapar}
\end{center}
\end{figure}

For each point $a_i$, choose a point $a_i^+$ and a point $a_i^-$ close to $a_i$ outside $D(\alpha_i)$ so that $a_i^+$ is on the positive side of $D(\alpha_i)$ (the side of the positive normal) and $a_i^-$ is on the negative side of $D(\alpha_i)$.
Similarly fix points $b_j^+$ and $b_j^-$ close to the $b_j$ and outside the $D(\beta_j)$.

Let $\gamma_{\Asc}^+(c)$ (resp. $\gamma_{\Asc}^-(c)$) be an arc parallel to $[a_{i(c)},c]_{D(\alpha_{i(c)})}$ from $a_{i(c)}^+$ to $c_{\Asc}^+$ (resp. from $a_{i(c)}^-$ to $c_{\Asc}^-$) that does not meet $D(\alpha_{i(c)})$.
Let $\gamma_{\Bdesc}^+(c)$ (resp. $\gamma_{\Bdesc}^-(c)$) be an arc parallel to $[c,b_{j(c)}]_{D(\beta_{j(c)})}$ from $c_{\Bdesc}^+$ to $b_{j(c)}^+$ (resp. from $c_{\Bdesc}^-$ to $b_{j(c)}^-$) that does not meet $D(\beta_{j(c)})$.

$$\begin{array}{ll}\gamma(c)_{\parallel}=&\frac12(\gamma_{\Asc}^+(c)+\gamma_{\Asc}^-(c))+\frac12(\gamma_{\Bdesc}^+(c)+\gamma_{\Bdesc}^-(c))\\
&+\frac14([c_{\Asc}^+,c_{\Bdesc}^+]+[c_{\Asc}^+,c_{\Bdesc}^-]+[c_{\Asc}^-,c_{\Bdesc}^+]+[c_{\Asc}^-,c_{\Bdesc}^-])
.\end{array}$$

Since the superscripts $+$ and the $-$ play the same roles in the above formula, $\gamma(c)_{\parallel}$ does not depend on the orientations of the $\alpha_i$ and the $\beta_j$. 
Set $a_{i\parallel}=\frac12(a_i^++a_i^-)$ and $b_{j\parallel}=\frac12(b_j^++b_j^-)$. Then
$\partial \gamma(c)_{\parallel}=b_{j(c)\parallel}-a_{i(c)\parallel}$.

\medskip

Set $\Link(\matchingo)_{\parallel}=\sum_{i=1}^g \gamma(c_i)_{\parallel}- \sum_{c\in \CaC} \CJ_{j(c)i(c)} \sigma(c) \gamma(c)_{\parallel}$ and note that $\Link(\matchingo)_{\parallel}$ is a cycle disjoint from $\Link(\matchingo)$. The cycle $\Link(\matchingo)$ depends neither on the orientations of the $\alpha_i$ and the $\beta_j$, nor on their order. Permuting the roles of the $\alpha_i$ and the roles of the $\beta_j$ reverses the orientations of $\Link(\matchingo)$ and $\Link(\matchingo)_{\parallel}$ and leaves $lk(\Link(\matchingo),\Link(\matchingo)_{\parallel})$ unchanged.

\subsection{A 2-cycle \texorpdfstring{$\twocycG(\CD)$}{G(D)}\texorpdfstring{ of $C_2(M)$}{} associated with a Heegaard diagram}
\label{subfundcyc}

The term $\lambdh$\index{Nt}{lctwo@$\lambdh$} will be defined from the homology class of the $2$--cycle $\twocycG(\CD)$ of $C_2(M)$ associated with the Heegaard diagram in the following proposition~\ref{propChMorse}, by the equality $[\twocycG(\CD)]=\lambdh[S]$ in $H_2(C_2(M);\QQ)$.
This term $\lambdh$ can be thought of as the main term of the formula, the other ones can be thought of as correction terms.

\begin{proposition}
\label{propChMorse}
Set $$\twocycG(\CD)=\sum_{(c,d) \in \CaC^2}\CJ_{j(c)i(d)}\CJ_{j(d)i(c)}\sigma(c)\sigma(d) (\gamma(c) \times \gamma(d)_{\parallel}) -\sum_{c \in \CaC} \CJ_{j(c)i(c)}\sigma(c)(\gamma(c) \times \gamma(c)_{\parallel}).$$
Then $\twocycG(\CD)$ is a $2$--cycle of $C_2(M)$. Its homology class $[\twocycG(\CD)]$ depends neither on the orientations of the $\alpha_i$ and the $\beta_j$, nor on their order. Permuting the roles of the $\alpha_i$ and the roles of the $\beta_j$ does not change it either.
\end{proposition}
\bp
Let us first prove that $\twocycG(\CD)$ is a $2$-cycle. Let $d \in \CaC$.
For any $j$,
$$\sum_{c\in \beta_j}\CJ_{j(d)i(c)}\sigma(c)=\sum_{i=1}^g\CJ_{j(d)i}\langle \alpha_i,\beta_j \rangle=\delta_{jj(d)}$$
and, for any $i$, $\sum_{c\in \alpha_i}\CJ_{j(c)i(d)}\sigma(c)=\sum_{j=1}^g\CJ_{ji(d)}\langle \alpha_i,\beta_j \rangle=\delta_{ii(d)}$.
Therefore, for any $d \in \CaC$,
$$\partial\left(\sum_{c \in \CaC}\CJ_{j(c)i(d)}\CJ_{j(d)i(c)}\sigma(c)\gamma(c)\right)=\CJ_{j(d)i(d)}(b_{j(d)}-a_{i(d)})=\CJ_{j(d)i(d)}\partial \gamma(d)$$ and
$$\begin{array}{lll}\partial \twocycG(\CD)&=&\sum_{d \in \CaC}\sigma(d)\CJ_{j(d)i(d)} (\partial \gamma(d)) \times \gamma(d)_{\parallel} -\sum_{c \in \CaC} \CJ_{j(c)i(c)}\sigma(c) (\partial \gamma(c)) \times \gamma(c)_{\parallel}\\
&&-\sum_{c \in \CaC}\CJ_{j(c)i(c)}\sigma(c) \gamma(c) \times \partial \gamma(c)_{\parallel}+\sum_{c \in \CaC} \CJ_{j(c)i(c)}\sigma(c)\gamma(c) \times \partial \gamma(c)_{\parallel}\\
&=&0.\end{array}$$

Since changing the orientation of $\alpha_{i(c)}$ leaves $\CJ_{j(d)i(c)}\sigma(c)$ invariant and changing the orientation of $\beta_{j(c)}$ leaves $\CJ_{j(c)i(d)}\sigma(c)$
invariant, the cycle $\twocycG(\CD)$ does not depend on the orientations of the $\alpha_i$ and the $\beta_j$. It clearly does not depend on the numbering.
It is also easy to see that permuting the roles of the $\alpha_i$ and the $\beta_j$ reverses the orientations of the $\gamma(c)$, changes $\CJ$ to the transposed matrix and does not change the cycle $\twocycG(\CD)$ either.
\eop

Note that $\lambdh$ is additive under connected sum of Heegaard diagrams, and therefore it is invariant under stabilisation of diagrams, but, as Example~\ref{exaheegfive} will show, it is not an invariant of Heegaard splittings.
In the next subsection, we state Proposition~\ref{propeval} that yields combinatorial formulae both for $\lambdh$ and for $lk(\Link(\CD,\matchingo),\Link(\CD,\matchingo)_{\parallel})$.

\subsection{Evaluating some \texorpdfstring{$2$}{2}--cycles \texorpdfstring{of $C_2(M)$}{ }}
\label{subeval}

When $d$ and $e$ are (possibly equal) crossings of $\alpha_i$, $[d,e]_{\alpha_i}=[d,e]_{\alpha}$ denotes the set of crossings from $d$ to $e$ (including them) along $\alpha_i$, or the closed arc from $d$ to $e$ in $\alpha_i$ depending on the context. Then $[d,e[_{\alpha}=[d,e]_{\alpha} \setminus \{e\}$.

Now, for such a part $I$ of $\alpha_i$, $$\langle I,\beta_j \rangle= \sum_{c\in I \cap \beta_j} \sigma(c).$$
We shall also use the notation $\halfbra$ for ends of arcs to say that an end is ``half-contained'' in an arc, and that it must be counted with coefficient $1/2$.
(``$[d,e\halfbra_{\alpha}=[d,e]_{\alpha} \setminus \{e\}/2$'' and ``$\halfbra d,e\halfbra_{\alpha}=[d,e\halfbra_{\alpha} \setminus \{d\}/2$'' so that $\halfbra d,d\halfbra_{\alpha}=\emptyset$.)

We use the same notation for arcs $[d,e\halfbra_{\beta_j}=[d,e\halfbra_{\beta}$ of $\beta_j$.
For example, if $d$ is a crossing of $\alpha_i\cap\beta_j$, then
$$\langle [d,d\halfbra_{\alpha},\beta_j \rangle=\frac{\sigma(d)}{2}$$
and
$$\langle [c,d\halfbra_{\alpha},[e,d\halfbra_{\beta} \rangle= \frac{\sigma(d)}{4} + \sum_{c\in [c,d[_{\alpha} \cap [e,d[_{\beta}} \sigma(c).$$

\begin{example}
\label{exaheegtwo}
In the diagram $\CD_1$ of Figure~\ref{figheegrpt},
$\langle [c,c\halfbra_{\alpha}, [c,c\halfbra_{\beta}\rangle=\frac14, 
\langle [c,c\halfbra_{\alpha}, [c,d\halfbra_{\beta}\rangle=\langle [c,d\halfbra_{\alpha}, [c,c\halfbra_{\beta}\rangle=\frac12, \langle [c,d\halfbra_{\alpha}, [c,d\halfbra_{\beta}\rangle=\frac54,\langle [c,c\halfbra_{\alpha}, \beta_1\rangle=\frac12$ and $\langle [c,d\halfbra_{\alpha}, \beta_1\rangle=\frac32.$
\end{example}

The following proposition is proved in Subsection~\ref{subpropeval}.

\begin{proposition}
\label{propeval}
For every curve $\alpha_i$ (resp. $\beta_j$), choose a basepoint $\pointp(\alpha_i)$ (resp. $\pointp(\beta_j)$). These choices being made, for two crossings $c$ and $d$ of $\CaC$, set
$$\begin{array}{ll}\ell(c,d)=&\langle [\pointp(\alpha(c)),c\halfbra_{\alpha}, [\pointp(\beta(d)),d\halfbra_{\beta}\rangle\\& -\sum_{(i,j) \in \{1,\dots,g\}^2}\CJ_{ji} \langle [\pointp(\alpha(c)),c\halfbra_{\alpha},\beta_j \rangle \langle \alpha_i ,[\pointp(\beta(d)),d\halfbra_{\beta} \rangle\end{array}$$\index{Nt}{lcd@$\ell(.,.)$}
where $\alpha(c)=\alpha_{i(c)}$ and $\beta(c)=\beta_{j(c)}$.
Then, for any $2$--cycle $\twocycG=\sum_{(c,d) \in \CaC^2}g_{c d} (\gamma(c) \times \gamma(d)_{\parallel})$ of $C_2(M)$,
$$[\twocycG]=\sum_{(c,d) \in \CaC^2}g_{c d}
\ell(c,d)[S]=\sum_{(c,d) \in \CaC^2}g_{c d}
\ell(d,c)[S].$$
\end{proposition}

We have the following immediate corollary of Proposition~\ref{propeval}.
\begin{corollary} \label{corcombexp} For any choice of $\ell$ as in Proposition~\ref{propeval}
$$\lambdh=\sum_{(c,d) \in \CaC^2}\CJ_{j(c)i(d)}\CJ_{j(d)i(c)}\sigma(c)\sigma(d) \ell(c,d) -\sum_{c \in \CaC} \CJ_{j(c)i(c)}\sigma(c)\ell(c,c)$$\index{Nt}{lctwo@$\lambdh$}
and
$$\begin{array}{ll}lk(\Link(\CD,\matchingo),\Link(\CD,\matchingo)_{\parallel})=&\sum_{(i,j) \in \{1,\dots,g\}^2}\ell(c_i,c_j) \\&+  \sum_{(c,d) \in \CaC^2}\CJ_{j(c)i(c)}\CJ_{j(d)i(d)}\sigma(c)\sigma(d) \ell(c,d)\\&- \sum_{(i,c) \in \{1,\dots,g\} \times \CaC }\CJ_{j(c)i(c)} \sigma(c)(\ell(c_i,c)+\ell(c,c_i)).\end{array}$$
\end{corollary}
\bp Recall $[\Link(\matchingo) \times \Link(\matchingo)_{\parallel}]=lk(\Link(\matchingo),\Link(\matchingo)_{\parallel})[S]$ in $H_2(C_2(M);\QQ)$. \eop

\begin{example}
\label{exaheegthree}
Again,
consider the diagram $\CD_1$ of Figure~\ref{figheegrpt}.
Choose $\pointp(\alpha_1)=\pointp(\beta_1)=c$. Using Example~\ref{exaheegtwo}, we get
$$\ell(c,c)=\frac14 -\frac18=\frac18, \ell(d,d)=\frac54 - \frac 98=\frac18, \ell(c,d)=\ell(d,c)=\frac12 - \frac38=\frac18.$$
For the diagram $\CD_2$ of Figure~\ref{figheegrpt}, choose $\pointp(\alpha_1)=\pointp(\beta_1)=c$ 
and $\pointp(\alpha_2)=\pointp(\beta_2)=e$.
Then we still have $\ell(c,c)=\ell(c,d)=\ell(d,c)=\ell(d,d)=\frac18$. Furthermore, $\ell(e,e)=0$ , and, as a nonsymmetric example, $\ell(c,e)=0$ {and} $\ell(e,c)=\frac18.$
Then $lk(\Link(\{c\}),\Link(\{c\})_{\parallel})=lk(\Link(\{c,e\}),\Link(\{c,e\})_{\parallel})=0$, and
$\ell_2(\CD_1)=\ell_2(\CD_2)=0.$
\end{example}

\subsection{Combinatorial definition of \texorpdfstring{$e(\pointw,\matchingo)$}{e(w,m)}}
\label{subframedHd}

Recall that we fixed a matching
$\matchingo=\{c_i; i \in \{1,2, \dots, g\}\}$
where the $c_i$ are crossings of $\alpha_i \cap \beta_{\rho^{-1}(i)}$ for a permutation $\rho$ of $\{1,2, \dots, g\}$.
Select an exterior point $\pointw$ of $\CD$.
These choices being fixed, represent the Heegaard diagram $\CD$ in a plane by removing a topological disk around $\pointw$ and by cutting the surface $\partial \handlebodA$ along the $\alpha_i$. The boundary of the removed topological disk will be pictured as a rectangle, and each $\alpha_i$ gives rise to two boundary components of the planar surface, which are copies of $\alpha_i$ denoted by $\alpha^{\prime}_i$ and $\alpha^{\prime\prime}_i$. They are drawn as circles. The crossing $c_i$ is located at the points with upward tangents of $\alpha^{\prime}_i$ and $\alpha^{\prime\prime}_i$, while the other crossings are located near 
the points with downward tangents as in Figure~\ref{figplansplit}. The curves $\beta_j$ intersect this picture as families of arcs, which begin and end at crossings with the $\alpha^{\prime}_i$ and the $\alpha^{\prime\prime}_i$ where they are horizontal. A diagram with these properties is called a \emph{rectangular diagram} of $(\CD,\matchingo,\pointw)$.

\begin{figure}[h]
\begin{center}
\begin{tikzpicture}
\useasboundingbox (0,-.2) rectangle (8.4,1.9);
\draw [red,->] (1,.5) arc (-90:270:.4);
\draw [red] (1,.2) node{$\alpha^{\prime}_1$} (2.8,.2) node{$\alpha^{\prime\prime}_1$};
\draw [red,->] (2.8,.5) arc (270:-90:.4);
\draw  (4.2,.9) node{$\dots$};
\draw [red,->] (5.6,.5) arc (-90:270:.4);
\draw [red] (5.6,.2) node{$\alpha^{\prime}_g$};
\draw [red,->] (7.4,.5) arc (270:-90:.4);
\draw [red] (7.4,.2) node{$\alpha^{\prime\prime}_g$};
\draw (0,-.1) -- (0,1.8) -- (8.4,1.8) -- (8.4,-.1) -- (0,-.1);
\draw [blue,thick,->] (2.55,.9) -- (2.25,.9);
\draw [blue,thick,->] (1.55,.9) -- (1.25,.9);
\draw [blue,thick,->] (7.15,.9) -- (6.85,.9);
\draw [blue,thick,->] (6.15,.9) -- (5.85,.9);
\draw [blue,thick,->] (3.05,.85) -- (3.35,.85);
\draw [blue,thick,->] (3,.75) -- (3.3,.75);
\draw [blue,thick,->] (3.25,.65) -- (2.95,.65);
\draw [blue,thick,->] (.85,.65) -- (.55,.65);
\draw [blue,thick,->] (.45,.85) -- (.75,.85);
\draw [blue,thick,->] (.5,.75) -- (.8,.75);
\draw [blue,thick,->] (5.1,.75) -- (5.4,.75);
\draw [blue,thick,->] (5.05,.85) -- (5.35,.85);
\draw [blue,thick,->] (5.35,.95) -- (5.05,.95);
\draw [blue,thick,->] (5.1,1.05) -- (5.4,1.05);
\draw [blue,thick,->] (7.6,.75) -- (7.9,.75);
\draw [blue,thick,->] (7.65,.85) -- (7.95,.85);
\draw [blue,thick,->] (7.95,.95) -- (7.65,.95);
\draw [blue,thick,->] (7.6,1.05) -- (7.9,1.05);
\draw [blue,thick,dotted] (2.25,.9) .. controls (1.9,.9) and (2.5,1.5) .. (2.8,1.5) .. controls (4,.9) and (4.9,1.05) .. (5.1,1.05) (7.9,1.05).. controls (8.4,1.05) and (8.2,1.7) .. (2.8,1.7) .. controls (2,1.7) and (1.9,.9) .. (1.55,.9);
\draw (1.6,.7) node{\small $c_1$} (2.2,.7) node{\small $c_1$} (6.2,.7) node{\small $c_g$} (6.8,.7) node{\small $c_g$};
\end{tikzpicture}
\caption{Rectangular diagram of $(\CD,\matchingo,\pointw)$}
\label{figplansplit}
\end{center}
\end{figure}
The rectangle has the standard parallelization of the plane. Then there is a map ``unit tangent vector'' from each partial projection of a beta curve $\beta_j$ in the plane to $S^1$. The total degree of this map for the curve $\beta_j$ is denoted by $\tilth(\beta_{j})$.
For a crossing $c \in \beta_j$, $\tilth(\halfbra c_{\rho(j)},c \halfbra_{\beta}) \in \frac12 \ZZ$\index{Nt}{de@$\tilth$} denotes the degree of the restriction of this map to the arc $\halfbra c_{\rho(j)},c \halfbra_{\beta}$. This degree is the average of the degrees of this map at the upward vertical vector and at the downward one.
For every $c \in \CaC$, define
$$\tilth(c)=\tilth(\halfbra c_{\rho(j(c))},c \halfbra_{\beta})
-\sum_{(r,s) \in \{1,\dots,g\}^2}\CJ_{sr}\langle \alpha_r,\halfbra c_{\rho(j(c))},c \halfbra_{\beta}\rangle\tilth(\beta_{s}),$$
where $\halfbra c,c \halfbra_{\beta}=\emptyset$. Then set $$e(\pointw,\matchingo)=e(\CD,\pointw,\matchingo)=\sum_{c \in \CaC} \CJ_{j(c)i(c)}\sigma(c)\tilth(c).$$\index{Nt}{ew@$e(\pointw,\matchingo)=e(\CD,\pointw,\matchingo)$}
In Section~\ref{subYepseta}, $e(\pointw,\matchingo)$ will be identified with an Euler class. See Proposition~\ref{propcontrob}.

\begin{example}
\label{exaheegfour}
For the Heegaard diagram $\CD_1$ equipped with the matching $\matchingo=\{c\}$, there are two choices for an exterior point $\pointw$ up to isotopy, the choice $\pointw$ of Figure~\ref{figheegrpt}, and the choice of a point $\pointw^{\prime}$ in the other connected component of $\partial \handlebodA \setminus \left( \alpha_1 \cup \beta_1\right)$.
These choices give rise to the two rectangular diagrams of $(\CD_1,\matchingo,\pointw)$ and $(\CD_1,\matchingo,\pointw^{\prime})$ shown in Figure~\ref{figrecgenone}.

\begin{figure}[h]
\begin{center}
\begin{tikzpicture}
\useasboundingbox (0,-.2) rectangle (8,1.9);
\draw [red,->] (1,.5) arc (-90:270:.4);
\draw [red] (1,.7) node{\scriptsize $\alpha^{\prime}_1$};
\draw  (.7,.9) node{\scriptsize $d$} (1.3,.9)node{\scriptsize $c$}  (3.1,.9) node{\scriptsize $d$} (2.5,.9)node{\scriptsize $c$};
\draw [blue] (2.7,.23) node{\scriptsize $\beta_1$}  (.9,1.63) node{\scriptsize $\beta_1$};
\draw [red,->] (2.8,.5) arc (270:-90:.4);
\draw [red] (2.85,.7) node{\scriptsize $\alpha^{\prime\prime}_1$};
\draw (0,-.1) -- (0,1.8) -- (3.8,1.8) -- (3.8,-.1) -- (0,-.1);
\draw [blue,thick,->]  (2.4,.9) .. controls (1.9,.9) and (1.4,1.5) .. (1,1.5);
\draw [blue,thick,->]  (1,1.5) .. controls (.6,1.5) and (.3,1.25) .. (.3,1.05) .. controls (.3,.9) .. (.6,.9);
\draw [blue,thick,->]  (3.2,.9) .. controls (3.5,.9) and (3.5,.7) .. (3.5,.5) .. controls (3.5,.3) and (3.3,.1) .. (2.6,.1);
\draw [blue,thick,->]  (2.6,.1) .. controls (2.3,.1) and (1.8,.9) .. (1.4,.9);
\draw (0,.1) node[right]{\scriptsize $(\CD_1,\matchingo,\pointw)$};
\fill (2.4,.5) circle (0.05);
\draw (2.25,.6) node{\scriptsize $\pointw^{\prime}$};
\begin{scope}[xshift=4.2cm]
\draw (3.8,.1) node[left]{\scriptsize $(\CD_1,\matchingo,\pointw^{\prime})$};
\draw [red,->] (1,.5) arc (-90:270:.4);
\draw [red] (1,.7) node{\scriptsize $\alpha^{\prime}_1$};
\fill (2.2,1.2) circle (0.05);
\draw (2.35,1.3) node{\scriptsize $\pointw$};
\draw  (.7,.9) node{\scriptsize $d$} (1.3,.9)node{\scriptsize $c$}  (3.1,.9) node{\scriptsize $d$} (2.5,.9)node{\scriptsize $c$};
\draw [blue]  (.9,1.68) node{\scriptsize $\beta_1$};
\draw [red,->] (2.8,.5) arc (270:-90:.4);
\draw [red] (2.85,.7) node{\scriptsize $\alpha^{\prime\prime}_1$};
\draw (0,-.1) -- (0,1.8) -- (3.8,1.8) -- (3.8,-.1) -- (0,-.1);
\draw [blue,thick,->]  (2.4,.9) .. controls (1.9,.9) and (1.4,1.4) .. (1,1.4);
\draw [blue,thick,->]  (1,1.4) .. controls (.6,1.4) and (.3,1.25) .. (.3,1.05) .. controls (.3,.9) .. (.6,.9);
\draw [blue,thick,->]  (3.2,.9) .. controls (3.5,.9) and (3.5,1.1) .. (3.5,1.3) .. controls (3.5,1.45) and (3,1.55) .. (1,1.55);
\draw [blue,thick,->]  (1,1.55) .. controls (.6,1.55) and (.1,1.4) .. (.1,.5) .. controls (.1,.3) and (.5,.2) .. (1,.2)  .. controls (1.5,.2) and (1.9,.9) .. (1.4,.9);
\end{scope}
\end{tikzpicture}
\caption{Rectangular diagrams of $(\CD_1,\{c\},\pointw)$ and $(\CD_1,\{c\},\pointw^{\prime})$}
\label{figrecgenone}
\end{center}
\end{figure}
For both rectangular diagrams, we have $\tilth(\halfbra c,c \halfbra_{\beta})=0$, $\tilth(c)=0$  and $\tilth(\halfbra c,d \halfbra_{\beta})=\frac12$ while $\tilth(\beta_1)=0$ for $(\CD_1,\{c\},\pointw)$ and $\tilth(\beta_1)=2$ for $(\CD_1,\{c\},\pointw^{\prime})$ so that $\tilth(d)=\frac12$ for $(\CD_1,\{c\},\pointw)$
and $\tilth(d)=-\frac12$ for $(\CD_1,\{c\},\pointw^{\prime})$.
Thus $e(\pointw^{\prime},\{c\})=-\frac14$ and $e(\pointw,\{c\})=\frac14$.
\end{example}

\subsection{Statement of the main theorem}

The main result of this article is the following theorem.

\begin{theorem}
\label{thmmain}
For any Heegaard diagram $\CD$ of a rational homology sphere $M$, for any exterior point $\pointw$ of $\CD$, and for any matching $\matchingo$ of $\CD$,
$$\Theta(M,\ComX(\CD,\pointw,\matchingo))=\lambdh + lk(\Link(\CD,\matchingo),\Link(\CD,\matchingo)_{\parallel}) -e(\CD,\pointw,\matchingo).$$
\end{theorem}

\begin{example}
\label{exaheegfive}
According to the computations of Examples~\ref{exaheegthree} and \ref{exaheegfour},
$$\Theta(\rpt,\ComX(\pointw,\{c\}))=-\frac14$$ and $\Theta(\rpt,\ComX(\pointw^{\prime},\{c\}))=\frac14$.
Since $\lambda(\rpt)=0$, this implies that $p_1(\ComX(\pointw^{\prime},\{c\}))=1$ and $p_1(\ComX(\pointw,\{c\}))=-1$.

Let us now evaluate the ingredients of our formula for the rectangular genus two diagram $(\CD_2,\{c,e\},\pointw)$ of Figure~\ref{figheegdtwo}. 
Recall from Example~\ref{exaheegthree} that $lk(\Link(\{c,e\}),\Link(\{c,e\})_{\parallel})=0$, and
$\ell_2(\CD_2)=0$ and observe $e(\CD_2,\pointw,\{c,e\})=\frac14$ so that $\Theta(\rpt,\ComX(\CD_2,\pointw,\{c,e\}))=-\frac14$.
\begin{figure}[h]
\begin{center}
\begin{tikzpicture}
\useasboundingbox (0,-.2) rectangle (7.7,1.9);
\draw [red,->] (1,.5) arc (-90:270:.4);
\draw  (.7,.9) node{\scriptsize $d$} (1.3,.9)node{\scriptsize $c$}  (3.1,.9) node{\scriptsize $d$} (2.5,.9)node{\scriptsize $c$};
\draw [blue] (.9,1.63) node{\scriptsize $\beta_1$};
\draw [red,->] (2.8,.5) arc (270:-90:.4);
\draw  [red]  (1,.7) node{\scriptsize $\alpha^{\prime}_1$} (2.85,.7) node{\scriptsize $\alpha^{\prime\prime}_1$};
\draw (0,-.1) -- (0,1.8) -- (7.6,1.8) -- (7.6,-.1) -- (0,-.1);
\draw [blue,thick,->]  (2.4,.9) .. controls (1.9,.9) and (1.4,1.5) .. (1,1.5);
\draw [blue,thick,->]  (1,1.5) .. controls (.6,1.5) and (.3,1.25) .. (.3,1.05) .. controls (.3,.9) .. (.6,.9);
\draw [blue,very thick,->]  (3.7,1.1) node{\scriptsize $\beta_1$}   (4.2,.9) -- (3.7,.9) (3.2,.9) -- (3.7,.9);
\begin{scope}[xshift=3.6cm]
\draw [red,->] (1,.5) arc (-90:270:.4);
\draw  (.7,.9) node{\scriptsize $f$} (1.3,.9)node{\scriptsize $e$}  (3.1,.9) node{\scriptsize $f$} (2.5,.9)node{\scriptsize $e$};
\draw [blue]  (1.9,1.1) node{\scriptsize $\beta_2$};
\draw [red,->] (2.8,.5) arc (270:-90:.4);
\draw [red] (1,.7) node{\scriptsize $\alpha^{\prime}_2$}(2.85,.7) node{\scriptsize $\alpha^{\prime\prime}_2$};
\draw [blue,very thick,dotted,->]  (1.4,.9) -- (1.9,.9) (2.4,.9) -- (1.9,.9);
\end{scope}
\draw [blue,thick,->] (3.7,.35) node{\scriptsize $\beta_1$}  (6.8,.9) .. controls (7.2,.9) and (8.2,.15) .. (3.7,.15);
\draw [blue,thick,->] (3.7,.15) .. controls (2.6,.15) and (1.9,.9) .. (1.4,.9);
\end{tikzpicture}
\caption{$(\CD_2,\{c,e\},\pointw)$}
\label{figheegdtwo}
\end{center}
\end{figure}

Consider the diagram $(\CD_3,\{c,e\},\pointw)$ of Figure~\ref{figheegdthree} obtained from $(\CD_2,\{c,e\},\pointw)$ 
by an isotopy of $\beta_2$ on $\partial \handlebodA$. 
\begin{figure}[h]
\begin{center}
\begin{tikzpicture}
\useasboundingbox (0,-.2) rectangle (7.7,1.9);
\draw [red,->] (1,.5) arc (-90:270:.4);
\draw (.8,1.1) node{\scriptsize $d$}   (.7,.9) node{\scriptsize $g$}   (.8,.73) node{\scriptsize $h$}  (1.3,.9)node{\scriptsize $c$}  (3.05,1.1) node{\scriptsize $d$}    (3.1,.9) node{\scriptsize $g$}   (3.05,.73) node{\scriptsize $h$}  (2.5,.9)node{\scriptsize $c$};
\draw [blue]  (.9,1.63) node{\scriptsize $\beta_1$};
\draw [red,->] (2.8,.5) arc (270:-90:.4);
\draw [red] (1,.33) node{\scriptsize $\alpha^{\prime}_1$}  (2.93,.33) node{\scriptsize $\alpha^{\prime\prime}_1$};
\draw (0,-.1) -- (0,1.8) -- (7.6,1.8) -- (7.6,-.1) -- (0,-.1);
\draw [blue,thick,->]  (2.4,.9) .. controls (1.9,.9) and (1.4,1.5) .. (1,1.5);
\draw [blue,thick,->]  (1,1.5) .. controls (.3,1.5) and (.2,1.06) .. (.62,1.06);
\draw [blue,thick,->]  (3.7,1.18) node{\scriptsize $\beta_1$}   (4.2,.9) .. controls (3.9,.9) and (3.9,.98) .. (3.69,.98) (3.18,1.06)  .. controls (3.4,1.06) and (3.4,.98) .. (3.69,.98);
\draw  [blue,thick,dotted]  (.4,.82) .. controls (.4,.9) .. (.6,.9) (.62,.74) .. controls (.4,.74) .. (.4,.82);
\draw [blue,thick] (.2,.78) node{\scriptsize $\beta_2$}  (.36,.78) -- (.4,.82) -- (.44,.78) ;
\begin{scope}[xshift=3.6cm]
\draw [red,->] (1,.5) arc (-90:270:.4);
\draw  (.7,.9) node{\scriptsize $f$} (1.3,.9)node{\scriptsize $e$}  (3.1,.9) node{\scriptsize $f$} (2.5,.9)node{\scriptsize $e$};
\draw [blue]  (1.75,.9) node{\scriptsize $\beta_2$};
\draw [red,->] (2.8,.5) arc (270:-90:.4);
\draw [red]  (1,.7) node{\scriptsize $\alpha^{\prime}_2$}  (2.85,.7) node{\scriptsize $\alpha^{\prime\prime}_2$};
\draw  [blue,thick,dotted,->]  (2.4,.9) .. controls (1.9,.9) and (1.6,.2) .. (1,.2);
\end{scope}
\draw  [blue,thick,dotted,->]  (4.6,.2) .. controls (4.2,.2) and (3.5,.74) .. (3.18,.74);
\draw  [blue,thick,dotted,->]  (3.2,.9) .. controls (3.6,.9) and (4.2,.35) .. (4.6,.35);
\draw  [blue,thick,dotted,->]  (4.6,.35) .. controls (5,.35) and (5.5,.8) .. (5,.9);
\draw [blue,thick,->] (3.7,.2) node{\scriptsize $\beta_1$}  (6.8,.9) .. controls (7.2,.9) and (8.2,0) .. (3.7,0);
\draw [blue,thick,->] (3.7,0) .. controls (2.6,0) and (1.9,.9) .. (1.4,.9);
\end{tikzpicture}
\caption{$(\CD_3,\{c,e\},\pointw)$}
\label{figheegdthree}
\end{center}
\end{figure}
The $\CJ_{ji}$ are the same as for $\CD_2$, and $\Link(\CD_3,\{c,e\})=\frac12(\gamma(c)-\gamma(d)) + \frac12(\gamma(g)-\gamma(h))$.
Again, choosing $\pointp(\alpha_1)=\pointp(\beta_1)=c$ 
and $\pointp(\alpha_2)=\pointp(\beta_2)=e$, $\ell(c,c)=\ell(c,d)=\ell(d,c)=\ell(d,d)=\frac18$ and $\ell(e,e)=0$.
For any crossing $x \in \{c,d,e,f\}$, $\ell(g,x)=\ell(h,x)$ and $\ell(x,g)=\ell(x,h)$.
Furthermore, $\ell(g,h)=\ell(h,g)$, $\ell(g,g)=\ell(h,g) + \frac14$ and $\ell(h,h)=\ell(h,g) - \frac14$
so that $lk(\Link(\CD_3,\{c,e\}),\Link(\CD_3,\{c,e\})_{\parallel})=0$ and
$\ell_2(\CD_3)=\CJ_{21}(\ell(h,h)-\ell(g,g))=\frac14$. Thus $\lambdh$ is not an invariant of Heegaard splittings.
Since $\tilth(g)=-\frac12$,
$e(\CD_3,\pointw,\{c,e\})=\frac14 + \frac14=\frac12$. Again $\Theta(\rpt,\ComX(\CD_3,\pointw,\{c,e\}))=-\frac14$.

A systematic study of the variations of the three ingredients of the formula under the moves that relate two Heegaard diagrams of a rational homology $3$-sphere is performed in \cite{lesHCcomb}.
\end{example}

\section{Propagators associated with Morse functions}
\setcounter{equation}{0}
\label{secpropout}

In this section, we introduce a propagator $\prop(\funcf,\metrigo)$ associated with a Morse function $\fMorse$ without minima and maxima of $\check{M}$, and with a metric $\metrigo$ that is standard outside $B_M$. This Morse propagator has been constructed in a joint work with Greg Kuperberg.
The pair $(\funcf,\metrigo)$ is supposed to give rise to the Heegaard diagram $\CD$ of Section~\ref{secformula} 
 as in Subsection~\ref{subonhd}.

We use the propagator $\prop(\funcf,\metrigo)$ (whose boundary is not associated with a combing) to prove Proposition~\ref{propeval}.
Similar propagators associated with more general Morse functions have been constructed by Watanabe in \cite{watanabeMorse}, independently.

\subsection{The Morse function \texorpdfstring{$\fMorse$}{f}}
Start with $\RR^3$ equipped with its standard height function $\fMorse_0$ and replace the parallelepiped $[0,2g] \times [0,4] \times [0,6]$ with a rational homology cube $C_M$\index{Nt}{CM@$C_M$} (which has the rational homology of a point) equipped with a Morse function $\fMorse$ that coincides with $\fMorse_0$ on $\partial \left([0,2g] \times [0,4] \times [0,6]\right)$, and that has
$2g$ critical points, $g$ points $a_1$, \dots, $a_g$ of index $1$, which are mapped to $1$ by
$\fMorse$, and $g$ points $b_1$, \dots, $b_g$ of index $2$, which are mapped to $5$ by
$\fMorse$. 
Let $\check{M}$ be the associated open manifold, and let $M$ be its one-point compactification. Equip $\check{M}$ with a Riemannian metric $\metrigo$ that coincides with the standard one outside $[0,2g] \times [0,4] \times [0,6]$.

The preimage $\handleboda$\index{Nt}{Ha@$\handleboda$} of $]-\infty, 2]$ under $\fMorse$ in $C_M$ has the standard representation of the bottom part of Figure~\ref{figHa}. Our standard representation of the preimage $\handlebodb$ of $[4,+\infty[$ under $\fMorse$ in $C_M$ is shown in the upper part of Figure~\ref{figHa}. It can be thought of as the complement of the bottom part in $[0,2g] \times [0,4] \times [0,6]$.\\

\begin{figure}[h]
\begin{center}
\begin{tikzpicture}
\useasboundingbox (-.2,0) rectangle (12.8,-8.2);
\draw (12.8,-3.5) -- (2.8,-3.5) -- (-.2,-1) -- (-.2,0) -- (9.8,0) -- (9.8,-1) -- (12.8,-3.5) -- (12.8,-2.5) -- (9.8,0) (9.8,-1) -- (-.2,-1);
\draw (12.3,-2.7) node{$\handlebodb$};
\begin{scope}[xshift=-5cm]
\draw (8,-1.95) ellipse (.15 and .25);
\draw [blue,thick] (8,-3.2) arc (-90:0:.5);
\draw [blue,thick,->] (8,-2.2) arc (90:0:.5);
\draw [blue,thick,dashed] (8,-2.2) arc (90:270:.5);
\draw [blue] (8.5,-2.7) node[right=-.1]{\scriptsize $\beta_1$};
\draw (7,-1.1) .. controls (7.2,-1.35) and (7.2,-3.2) .. (8,-3.2) .. controls (8.5,-3.2) and (9.4,-3.2) .. (9.5,-3.45);
\end{scope}
\draw (5.5,-2.8) node{$\dots$};
\draw (8,-1.95) ellipse (.15 and .25);
\draw [blue,thick] (8,-3.2) arc (-90:0:.5);
\draw [blue,thick,->] (8,-2.2) arc (90:0:.5);
\draw [blue] (8.5,-2.7) node[right=-.1]{\scriptsize$\beta_g$};
\draw [blue,thick,dashed] (8,-2.2) arc (90:270:.5);
\draw (7,-1.1) .. controls (7.2,-1.35) and (7.2,-3.2) .. (8,-3.2) .. controls (8.5,-3.2) and (9.4,-3.2) .. (9.5,-3.45);
\begin{scope}[yshift=-8.2cm]
\draw (1.3,1.5) .. controls (2,1.5) and (2,3.2) .. (3,3.2);
\draw [->] (4.7,1.5) .. controls (4,1.5) and (4,3.2) .. (3,3.2) ;
\draw (12.8,3.5) -- (2.8,3.5) -- (-.2,1) -- (-.2,0) -- (9.8,0) -- (9.8,1) -- (12.8,3.5) -- (12.8,2.5) -- (9.8,0) (9.8,1) -- (-.2,1);
\draw (12.3,2.6) node{$\handleboda$};
\draw (3,2) circle (.5);
\draw (5.5,1.5) node{$\dots$};
\draw (6.3,1.5) .. controls (7,1.5) and (7,3.2) .. (8,3.2) .. controls (9,3.2) and (9,1.5) .. (9.7,1.5);
\draw (8,2) circle (.5);
\draw [red,thick,->] (3,3.2) .. controls (2.9,3.2) and (2.7,3) .. (2.7,2.85);
\draw [red,thick] (3,2.5) .. controls (2.9,2.5) and (2.7,2.7) .. (2.7,2.85) node[right=-.1]{\scriptsize $\alpha_1$};
\draw [red,thick,dashed] (3,3.2) .. controls (3.1,3.2) and (3.3,3) .. (3.3,2.85) .. controls (3.3,2.7) and (3.1,2.5) .. (3,2.5);
\draw [red,thick,->] (8,3.2) .. controls (7.9,3.2) and (7.7,3) .. (7.7,2.85);
\draw [red,thick] (8,2.5) .. controls (7.9,2.5) and (7.7,2.7) .. (7.7,2.85) node[right=-.1]{\scriptsize $\alpha_g$};
\draw [red,thick,dashed] (8,3.2) .. controls (8.1,3.2) and (8.3,3) .. (8.3,2.85) .. controls (8.3,2.7) and (8.1,2.5) .. (8,2.5);
\end{scope}
\end{tikzpicture}
\caption{$\handleboda$ and $\handlebodb$}
\label{figHa}
\end{center}
\end{figure}

The two-dimensional ascending manifold of $a_i$ is oriented arbitrarily, its closure is denoted by $\Asc_i$\index{Nt}{Ai@$\Asc_i$}. Its intersection with $\handleboda$ is denoted by $D(\alpha_i)$. The boundary of $D(\alpha_i)$ is denoted by $\alpha_i$.
The descending manifold of $a_i$ is made of two half-lines $\linhp(a_i)$\index{Nt}{Laiplus@$\linhp(a_i)$} and $\linhn(a_i)$ starting as vertical lines and ending at $a_i$. The one with the orientation of the positive normal to $\Asc_i$ is called $\linhp(a_i)$\index{Nt}{Lai@$\linh(a_i)$}. Thus $\linh(a_i)=\linhp(a_i)\cup (-\linhn(a_i))$ is the descending manifold of $a_i$.

\begin{figure}[h]
\begin{center}
\begin{tikzpicture}
\useasboundingbox (.8,.8) rectangle (10.2,4.2);
\draw [thick,->] (1,1) -- (1,2) node[left]{\scriptsize $\linhp(a_i)$};
\draw [thick] (1,2) .. controls (1,2.7) and (1.3,3) .. (2,3);
\draw [fill=gray!20,draw=white] (2,3) ellipse (.6 and 1);
\draw [red]  (2,3) ellipse (.6 and 1);
\draw [thick,->] (3,1) -- (3,2) node[right]{\scriptsize $\linhn(a_i)$};
\draw [draw=white,double=black,very thick] (3,2) .. controls (3,2.7) and (2.7,3) .. (2,3);
\draw [->] (2,3) node[left] {\small $a_i$} -- (2,4.3);
\fill (2,3) circle (0.1);
\draw [->] (2,3) -- (2,1.7);
\draw [->] (2,3) -- (2.3,1.85);
\draw (2,2.5) node[rectangle,rotate=20,fill=white] {\tiny $D(\alpha_i)$};
\draw [->] (2,3) -- (2.3,4.15);
\draw [->] (2,3) -- (2.65,3.7);
\draw [->] (2,3) -- (2.8,3.2);
\draw [->] (2,3) -- (1.7,4.15);
\draw [red,->] (2,4) arc (90:180:.6 and 1);
\draw [red] (1.4,3) node[left]{\scriptsize $\alpha_i$};
\draw [thick,-<] (9,4.4) -- (9,3) node[right]{\scriptsize $\linhp(b_j)$};
\draw [thick] (8,2) .. controls (9,2.4) .. (9,3);
\draw [fill=gray!20,draw=white] (8,2) circle (1); \draw [blue]  (8,2) circle (1);
\draw [->] (6.7,2) -- (7.5,2);
\draw (7.5,2) -- (8,2);
\draw [->] (9.3,2) -- (8.5,2);
\draw (8.5,2) -- node[above]{\small $b_j$} (8,2);
\draw [->] (8,3.3) -- (8,2.5);
\draw (8,2.5) -- (8,2);
\draw [->] (8,.7) -- (8,1.2);
\draw (8,1.2) -- (8,2);
\draw [thick,-<] (7,3.6) -- (7,3) node[left]{\scriptsize $\linhn(b_j)$};
\draw [draw=white,double=black,very thick] (8,2) .. controls (7,1.6) .. (7,3);
\fill (8,2) circle (0.1);
\draw (8,1.5) node[rectangle,fill=white] {\tiny $D(\beta_j)$};
\draw [blue] (7.4,3.1) node{\scriptsize $\beta_j$};
\draw [blue,-<] (9,2) arc (0:130:1);
\end{tikzpicture}
\caption{$\linhp(a_i)$, $\linhn(a_i)$, $\linhp(b_j)$, $\linhn(b_j)$}
\label{figda}
\end{center}
\end{figure}

Symmetrically, the two-dimensional descending manifold of $b_j$ is oriented arbitrarily, its closure is denoted by $\Bdesc_j$. The $\Bdesc_j$ are assumed to be transverse to the $\Asc_i$ outside the critical points.
The ascending manifold of $b_j$ is made of two half-lines $\linhp(b_j)$ and $\linhn(b_j)$ starting at $b_j$ and ending as vertical lines. The one with the orientation of the positive normal to $\Bdesc_j$ is called $\linhp(b_j)$. Thus $\linh(b_j)=\linhp(b_j)-\linhn(b_j)$ is the ascending manifold of $b_j$. See Figure~\ref{figda}.

Let $$H_{a,2}=C_M \cap \fMorse^{-1}(2)$$\index{Nt}{Hatwo@$H_{a,2}$} and similarly define
$H_{b,4}=C_M \cap \fMorse^{-1}(4)$.
The preimage of $[2,4]$ in $C_M$ is the product $H_{a,2} \times [2,4]$. Its intersection with
$\Asc_i$ is $-\alpha_i \times [2,4]$ and its intersection with $\Bdesc_j$ is $\beta_j \times [2,4]$.
Each crossing $c$ of $\alpha_i \cap \beta_j$ has a sign $\sigma(c)$ and an associated flow line $\gamma(c)$ from $a_i$ to $b_j$ oriented as such.

Note the following lemma.
\begin{lemma}
\label{lemnorpos}
Let $c\in \alpha_i \cap \beta_j$.
Along $\gamma(c)$, $\Asc_i$ is cooriented by $\sigma(c)\beta_j$ and $\Bdesc_j$ is cooriented by $\sigma(c)\alpha_i$.
$$\Bdesc_j\cap \Asc_i=\sum_{c \in \alpha_i \cap \beta_j}\sigma(c) \gamma(c).$$
\end{lemma}
\eop

\subsection{The propagator \texorpdfstring{$\prop(\funcf,\metrigo)$}{P(f,g)}}

Let $s_{\phi}(\check{M})$\index{Nt}{sphi@$s_{\phi}(\check{M})$} be the closure in $U\check{M}$ of the (graph of the) section of $U\check{M}_{|\check{M} \setminus \{a_i,b_i;i \in \{1,\dots,g\}\}}$ directed by the gradient of $\fMorse$. This closure contains the restriction of the unit tangent bundle to the critical points, up to orientation. Let $\phi$ be the flow associated with the gradient of $\fMorse$.
Let $\preprop_{\phi}$\index{Nt}{Pphi@$\preprop_{\phi}$} be the closure in $C_2(M)$ of the image of
$$\begin{array}{lll}
\left(\check{M} \setminus \{a_i,b_i;i \in \{1,\dots,g\}\} \right)\times ]0,+\infty[ & \rightarrow & C_2(M)\\
(x,t) & \mapsto & (x,\phi_t(x)),
\end{array}$$
let $\left((\Bdesc_j \times \Asc_i)\cap C_2(M)\right)$ denote the closure of $\left((\Bdesc_j \times \Asc_i)\cap (\check{M}^2 \setminus \mbox{diagonal}) \right)$ in $C_2(M)$,
set 
$$\preprop_{\CI} = \sum_{(i,j) \in \{1,\dots,g\}^2} \CJ_{ji} \left((\Bdesc_j \times \Asc_i)\cap C_2(M)\right) \;\;\;\;\;\;\;\;\mbox{and}\;\;\;\;\;\;\;\;\prop(\funcf,\metrigo)=\preprop_{\phi}+ \preprop_{\CI}$$\index{Nt}{PI@$\preprop_{\CI}$}\index{Nt}{Pfg@$\prop(\funcf,\metrigo)$}

Let $\vec{v}$ be the upward vector in $S^2$, and let \index{Nt}{delod@$\partial_{od}$}$$\partial_{od}=p_{\infty}^{-1}(\vec{v}) \cap \left(\partial C_2(M) \setminus U\check{M} \right)$$ be a boundary part {\em outside the diagonal of $\check{M}^2$.\/}
(If $\vec{v}_{\infty}$ denotes the upward vertical vector in the boundary of the compactification $C_1(M)$ of $\check{M}$, then
$\partial_{od}$ contains
$\left(- \check{M} \times \vec{v}_{\infty} - \left((-\vec{v}_{\infty}) \times \check{M}\right)\right)$.)

\begin{theorem}[Kuperberg--Lescop]
\label{thmbordFfii}
The $4$--chain $\prop(\funcf,\metrigo)$ is a propagator and its boundary, which lies in $\partial C_2(M)$, is
$$ \partial \prop(\funcf,\metrigo) = \partial_{od} + \sum_{c \in \CaC} \CJ_{j(c)i(c)} \sigma(c) U\check{M}_{|\gamma(c)} +\overline{s_{\phi}(\check{M})}$$
where $\overline{s_{\phi}(\check{M})}$ is the closure of $s_{\phi}(\check{M})$ in $\partial C_2(M)$.
\end{theorem}
\bp
The expression of $\partial \prop(\funcf,\metrigo)$ is the immediate consequence of the following two lemmas. Then it is easy to see that, for a tiny sphere $\partial B(x)$ around a point $x$ outside the $\gamma(c)$, $\langle (x \times \partial B(x)),\prop(\funcf,\metrigo) \rangle_{C_2(M)}$ is the algebraic intersection in $U\check{M}$ of a fiber and the section $s_{\phi}(\check{M})$, which is one.
\eop

Note that $U\check{M}_{|\gamma(c)}$ is diffeomorphic to $S^2 \times \gamma(c)$. For simplicity, $U\check{M}_{|\gamma(c)}$ will sometimes be simply denoted by $S^2 \times \gamma(c)$, or by $S^2 \times_{\tau} \gamma(c)$ when the parallelization $\tau$ that induces such a diffeomorphism matters.

\begin{lemma}
$$ \partial \preprop_{\phi}=\partial_{od} +\overline{s_{\phi}(\check{M})} -
\sum_{i=1}^g \linh(a_i)\times \Asc_i
- \sum_{j=1}^g \Bdesc_j \times \linh(b_j)$$
\end{lemma}
\bp
The boundary of $\preprop_{\phi}$ is made of $\left(\partial_{od} +\overline{s_{\phi}(\check{M})}\right)$
and some other parts coming from the critical points.
Let us look at the part coming from $a_i$, where the closures $\linhp(a_i)$ and $\linhn(a_i)$ of flow lines stop and closures of flow lines of $\Asc_i$ start.
Consider a tubular neighborhood $$D^2 \times \linhp(a_i)= \{(u\exp(i\theta),y);u\in [0,1],\theta \in [0,2\pi[,y \in \linhp(a_i)\}$$ around $\linhp(a_i)$, where $\phi_t((u\exp(i\theta),y))$ reads $(u^{\prime}\exp(i\theta),y^{\prime})$ for some $u^{\prime}\geq u$, for $t\geq 0$ and for $u$ small enough, so that $\theta$ is preserved by the flow.
When $u$ approaches $0$, the flow line through $(u\exp(i\theta),y)$ approaches 
$\linhp(a_i) \cup \lintheta(\Asc_i)$ where $\lintheta(\Asc_i)$ is the closure of a flow line in $\Asc_i$ determined by $\theta$, for {\em generic\/} $\theta$ (which are $\theta$ such that this closure does not end at a $b_j$).
In particular, $\preprop_{\phi}$ contains $\pm (\linhp(a_i) \times \Asc_i)$, and we examine more closely what $\preprop_{\phi}$ looks like near $\left( \linhp(a_i) \times \fMorse^{-1}([1,+\infty[)\right)$.

Blow up $0$ in $D^2$ to obtain an annulus $B\!\ell(D^2,0)$.
Blow up $\linhp(a_i)$ in $D^2 \times \linhp(a_i)$ to
replace $\linhp(a_i)$ by its unit normal bundle $S^1 \times \linhp(a_i)=\{(\exp(i\theta),y)\}$.
Let $B\!\ell(D^2,0) \times \linhp(a_i)$ denote the blown-up tubular neighborhood.
Fix a fiber $B\!\ell(D^2,0)_0=\{(u,\exp(i\theta));u\in [0,1],\exp(i\theta) \in S^1\}$ of $B\!\ell(D^2,0) \times \linhp(a_i)$, and its natural projection onto the disk $D^2_0=\{u \exp(i\theta)\}$. Then there are topological embeddings

 $$\begin{array}{llll} E_1 \colon &D^2_0 \times ]-\infty,1[&\rightarrow & \fMorse^{-1}(]-\infty,1[)\\
    &(u\exp(i\theta),x)& \mapsto & m=E_1(u\exp(i\theta),x)
   \end{array}$$
such that $m$ is on the flow line through the point $u\exp(i\theta)$ of $D^2_0$ and $\fMorse(m)=x$, and
 $$\begin{array}{llll} E_2 \colon &B\!\ell(D^2,0)_0 \times ]1,5[&\rightarrow & \fMorse^{-1}(]1,5[)\\
    &(u,\exp(i\theta),x)& \mapsto & n=E_2(u,\exp(i\theta),x)
   \end{array}$$
such that $\fMorse(n)=x$, $n$ is on the flow line through the point $u\exp(i\theta)$ of $B\!\ell(D^2,0)_0$ if $u\neq 0$, and $E_2(0,\exp(i\theta),x) \in \lintheta(\Asc_i)$.
Then $\preprop_{\phi}$ intersects $\fMorse^{-1}(]-\infty,1[) \times \fMorse^{-1}(]1,5[)$ near $\linhp(a_i) \times \fMorse^{-1}(]1,5[)$ as the image of the continuous embedding
$$\begin{array}{llll} E \colon &B\!\ell(D^2,0)_0 \times ]-\infty,1[ \times ]1,5[ &\rightarrow & \check{M}^2\\
    &(u,\exp(i\theta),x_1,x_2)& \mapsto & \left(E_1(u\exp(i\theta),x_1),E_2(u,\exp(i\theta),x_2)\right)
   \end{array}$$
and the boundary of $\preprop_{\phi}$ contains $E(\partial_b B\!\ell(D^2,0)_0 \times ]-\infty,1[ \times ]1,5[)$ where $$\partial_b B\!\ell(D^2,0)_0=-S^1$$ is the preimage of $(0 \in D^2_0)$. The closure of $]-\infty,1[$ is naturally identified with $\linhp(a_i)$ via $E_1$,
so that the boundary of $\preprop_{\phi}$ contains $\linhp(a_i) \times E_2(S^1 \times ]1,5[)$ and it is easy to conclude that the boundary part coming from $a_i$ near $\linhp(a_i) \times \fMorse^{-1}([1,+\infty[)$ is $(-\linhp(a_i))\times \Asc_i$ (with a minor $2$--dimensional abuse of notation around $a_i$). We similarly find $\linhn(a_i)\times \Asc_i$ in $\partial \preprop_{\phi}$, and the part of $\partial \preprop_{\phi}$ coming from $a_i$ is $(\linhn(a_i)-\linhp(a_i))\times \Asc_i$.

For $\linhp(b_j)$, we similarly get a part of $\partial \preprop_{\phi}$ $$-\bigcup_{\exp(i\theta) \in S^1} \mbox{flow line}\; \lintheta(\Bdesc_j)   \times \linhp(b_j),$$ locally oriented as
$(\mbox{flow line}\; \lintheta(\Bdesc_j)) \times (S^1 \times \linhp(b_j))$
where $\Bdesc_j$ locally reads $(-\lintheta(\Bdesc_j) \times S^1)$, and the boundary part coming from $b_j$ is $\Bdesc_j \times (\linhn(b_j)-\linhp(b_j))$.
The two boundary parts $(-\linh(a_i))\times \Asc_i$ and $\Bdesc_j \times (-\linh(b_j))$ intersect along a two-dimensional locus, and the $3$-cycle $\partial \preprop_{\phi}$ is completely described in the statement.
\eop

\begin{figure}[h]
\begin{center}
\begin{tikzpicture}[scale=.6]
\useasboundingbox (-4.2,-4.1) rectangle (4.2,4.2);
\draw [fill=gray!20, draw=white] (0,0) circle (4);
\draw (4,0) arc (0:220:4);
\draw (4,0) arc (0:-40:4);
\draw [dotted] (0,-4) arc (-90:-40:4);
\draw [dotted] (0,-4) arc (-90:-140:4);
\draw [fill=white, draw=white] (4,0) circle (1.2);
\draw [fill=white, draw=white] (-1.2,4) arc (-180:0:1.2) -- (0,4.2) -- (-1.2,4);
\draw [fill=white, draw=white] (-4,0) circle (1.2);
\draw (-2.8,0) arc (0:83:1.2);
\draw (-2.8,0) arc (0:-83:1.2);
\draw [->] (-2.8,0) arc (0:60:1.2);
\draw [->] (-2.8,0) arc (0:-60:1.2);
\draw (2.8,0) arc (180:97:1.2);
\draw (2.8,0) arc (-180:-97:1.2);
\draw [->] (2.8,0) arc (180:120:1.2);
\draw [->] (2.8,0) arc (-180:-120:1.2);
\draw (3.9,-.7) node{\scriptsize $\linhn(b_j)$} (3.9,.67) node{\scriptsize $\linhp(b_j)$};
\draw (0,2.8) arc (-90:-173:1.2);
\draw (0,2.8) arc (-90:-7:1.2);
\draw [->] (0,2.8) arc (-90:-150:1.2);
\draw [->] (0,2.8) arc (-90:-30:1.2);
\draw (-1.15,3.56) node[right]{\scriptsize $\linhn(b_k)$} (.8,3.1) node[right]{\scriptsize $\linhp(b_k)$};
\draw (0,0) circle (1);
\draw [->] (-1,0) arc (180:270:1);
\draw (0,-1.5) node{\scriptsize $\alpha_i$};
\draw [->] (0,0) -- (-1.9,0);
\draw [->] (0,0) -- (1.9,0);
\draw [->] (0,0) -- (0,1.9);
\draw (0,0) -- (45:4);
\draw [->] (0,0) -- (45:1.9);
\draw (0,0) -- (35:4);
\draw [->] (0,0) -- (35:1.9);
\draw [->] (0,0) -- (15:1.9);
\draw (15:1.9) .. controls (15:2.5) and (25:3) .. (25:4);
\draw (-1.9,-.5) node{\scriptsize $\gamma(c_3)$} (1.9,-.5) node{\scriptsize $\gamma(c_1)$} (0,1.9) node[left]{\scriptsize $\gamma(c_2)$} (7,1.3);
\draw (-1.9,0) -- (-2.8,0);
\draw (1.9,0) -- (2.8,0);
\draw (0,1.9) -- (0,2.8);
\fill (0,0) node[below]{\scriptsize $a_i$} circle (0.1);
\end{tikzpicture}
\caption{The interior of $\Asc_i$ (In the figure $\sigma(c_1)=1=-\sigma(c_2)$.)}
\label{figinA}
\end{center}
\end{figure}
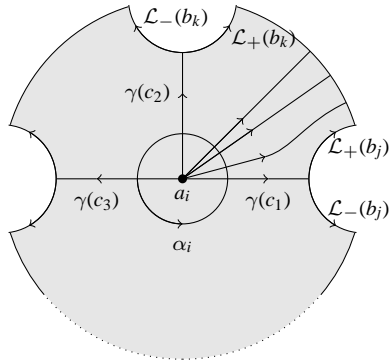

\begin{lemma}
$$ \partial \preprop_{\CI}=
\sum_{i=1}^g \linh(a_i)\times \Asc_i
+ \sum_{j=1}^g \Bdesc_j \times \linh(b_j)
+\sum_{c \in \CaC} \CJ_{j(c)i(c)} \sigma(c) (S^2 \times \gamma(c))$$
\end{lemma}
\bp
The interior of a figure similar to Figure~\ref{figinA} embeds in the closure $\Asc_i$ of the ascending manifold of $a_i$ in $\check{M}$. The whole closure is obtained by attaching such an open disk to the ascending manifolds $(\linh(b_j)=\linhp(b_j)-\linhn(b_j))$ of the $b_j$.

Recall that when the sign $\sigma(c)$ of a crossing $c \in \alpha_i \cap \beta_j$ is $1$, $\beta_j$ is positively normal to
$\Asc_i$ and $\alpha_i$ is positively normal to
$\Bdesc_j$ along the interior of $\gamma(c)$. See Lemma~\ref{lemnorpos}.

When $\Asc_i$ arrives at $b_j$ by a line $\gamma(c)$, it opens to $\linh(b_j)$ and we find
$$\partial \Asc_i= \sum_{j=1}^g \sum_{c \in \alpha_i \cap \beta_j} \sigma(c) \linh(b_j)=\sum_{j=1}^g \langle \alpha_i,\beta_j \rangle_{H_{a,2}}\linh(b_j)$$
$$\partial \Bdesc_j= \sum_{i=1}^g \langle \alpha_i,\beta_j \rangle_{H_{a,2}}\linh(a_i).$$
Near a connecting flow line  $\gamma(c)$, $\Bdesc_j$ is parametrized by $\beta_j\times\gamma(c)(]1,5[) $ and $\Asc_i$ is parametrized by $\gamma(c)(]1,5[) \times \alpha_i$.
Near the diagonal of such a line, $\Bdesc_j \times \Asc_i$ is parametrized by the height of the first point in $[1,5]$ followed by the tiny difference (second point minus first point), which is parametrized by $(\mbox{height difference}, \alpha_i, -(-\beta_j))$,
where one minus sign in front of $\beta_j$ comes from the permutation of the parameters, and the other one comes from the fact that $\beta_j$ is now used to parametrize the difference,
so that we get
$ \sum_{c \in \CaC} \CJ_{j(c)i(c)} \sigma(c) (S^2 \times \gamma(c)) $
in the boundary. \eop

\subsection{Using the propagator to prove Proposition~\ref{propeval}}
\label{subpropeval}

Let $\iota$\index{Nt}{iota@$\iota$} denote the continuous involution of $C_2(M)$ that exchanges two points in a pair of $(\check{M}^2\setminus \mbox{diag})$. Note that $\iota$ reverses the orientation of $C_2(M)$.

\begin{lemma}
\label{lemsymlink}
For any $2$--cycle $\twocycG=\sum_{(c,d)\in \CaC^2}g_{c d} (\gamma(c) \times \gamma(d)_{\parallel})$ of $C_2(M)$,
$$[\twocycG]=\left[\sum_{(c,d)\in \CaC^2}g_{c d} (\gamma(d) \times \gamma(c)_{\parallel})\right].$$
\end{lemma}
\bp
With the notation of Subsection~\ref{subparflow}, for $\varepsilon=\pm $ and $\eta=\pm $,
let $$\gamma(c)_{\nu^{\varepsilon}(\Asc)\nu^{\eta}(\Bdesc)}=\gamma_{\Asc}^{\varepsilon}(c) + [c^{\varepsilon}_{\Asc},c_{\Bdesc}^{\eta}] + \gamma_{\Bdesc}^{\eta}(c)$$
so that
$$\gamma(c)_{\parallel}=\frac14\left(\gamma(c)_{\nu^+(\Asc)\nu^+(\Bdesc)} + \gamma(c)_{\nu^+(\Asc)\nu^-(\Bdesc)} +\gamma(c)_{\nu^-(\Asc)\nu^+(\Bdesc)} +\gamma(c)_{\nu^-(\Asc)\nu^-(\Bdesc)} \right).$$
Then for any $\varepsilon$ and for any $\eta$,
$$\twocycG^{\varepsilon,\eta}=\sum_{(c,d)\in \CaC^2}g_{c d} \gamma(c) \times \gamma(d)_{\nu^{\varepsilon}(\Asc)\nu^{\eta}(\Bdesc)}$$
is a $2$--cycle homotopic to
$$\twocycG_s^{\varepsilon,\eta}=\sum_{(c,d)\in \CaC^2}g_{c d} \gamma(c)_{\nu^{-\varepsilon}(\Asc)\nu^{-\eta}(\Bdesc)} \times \gamma(d).$$
Now,
$$\iota(\twocycG_s^{\varepsilon,\eta})=-\sum_{(c,d)\in \CaC^2}g_{c d}\gamma(d) \times \gamma(c)_{\nu^{-\varepsilon}(\Asc)\nu^{-\eta}(\Bdesc)},$$ and,
since $[\iota_{\ast}(S)]=-[S]$, $\iota_{\ast}$ is the multiplication by $(-1)$ in $H_2(C_2(M);\QQ)$, and $(-\iota(\twocycG_s^{\varepsilon,\eta}))$ is homologous to $\twocycG^{\varepsilon,\eta}$.
Since $\twocycG$ is the average of the $\twocycG^{\varepsilon,\eta}$, and since $\left(\sum_{(c,d)\in \CaC^2}g_{c d} \gamma(d) \times \gamma(c)_{\parallel}\right)$ is the average of the $(-\iota(\twocycG_s^{\varepsilon,\eta}))$, the lemma is proved.
\eop

In order to prove Proposition~\ref{propeval}, we are now left with the proof that $$[\twocycG]=\sum_{(c,d)\in \CaC^2}g_{c d}
\ell(c,d)[S].$$ We prove this by transforming
the $\gamma(c)$
into $$\gamma(c)_{\nu(\Bdesc)}=\frac12 \left(\gamma(c)_{\nu^+(\Bdesc)}+\gamma(c)_{\nu^-(\Bdesc)}\right)$$ where $\gamma(c)_{\nu^+(\Bdesc)}$ (resp. $\gamma(c)_{\nu^-(\Bdesc)}$) is obtained from $\gamma(c)$ by pushing it much less than the distance between $\gamma(c)_{\nu^{\varepsilon}(\Asc)\nu^{\eta}(\Bdesc)}$ and $\gamma(c)$, in the direction of the positive (resp. negative) normal to $\Bdesc_{j(c)}$, except in the neighborhood of $a_{i(c)}$, where
\begin{itemize}
\item $\gamma(c)_{\nu(\Bdesc)}$ is in $\Asc_{i(c)}$ and it is transverse to the $\Bdesc_j$,
\item the starting points near $a_i$ of all the $\gamma(c)_{\nu^+(\Bdesc)}$ and the $\gamma(c)_{\nu^-(\Bdesc)}$ for which $i(c)=i$ coincide, they are denoted by $a_{i,\nu(\Bdesc)}$,
\item this starting point $a_{i,\nu(\Bdesc)}$ does not belong to the sheets of the $\Bdesc_j$ corresponding to crossings of $\alpha_i$ and the $\beta_j$, (these sheets meet along $\linh(a_i)$),
\item the first encountered sheet from $a_{i,\nu(\Bdesc)}$ when turning around $\linh(a_i)$ like $\alpha_i$ is the sheet of $\pointp(\alpha_i)$.
\end{itemize}
See the local infinitesimal picture of Figure~\ref{figNB}.
Recall from Lemma~\ref{lemnorpos} that $\alpha_i$ is the positive normal to $\Bdesc_j$ along flow lines through positive crossings.

\begin{figure}[h]
\begin{center}
\begin{tikzpicture}[scale=.77]
\useasboundingbox (-5,-4.1) rectangle (5,4.3);
\begin{scope}[xshift=-4.4cm]
\draw [fill=gray!20, draw=white] (0,0) circle (4);
\draw [dashed, -latex] (-4,0) arc (180:360:4);
\draw [-latex] (108:1) -- (108:4) (0,0) --(108:1);
\draw [-latex] (0,0) -- (96:4);
\draw [-latex] (0,0) -- (78:4);
\draw [dashed, -latex] (.4,0) -- (113:4) ;
\draw [dashed, -latex] (.4,0) -- (91:4);
\draw [dashed, -latex] (.4,1) -- (.4,3.98) (.4,0) -- (.4,1);
\draw (-.2,1) node[left]{\scriptsize $\gamma(f)$} (-.2,-.2) node{\scriptsize  $a_i$} (.2,-.2) node[right]{\scriptsize  $a_{i,\nu(\Bdesc)}$};
\draw (.4,1) node[right]{\scriptsize  $\gamma(\pointp(\alpha_i))_{\nu^+(\Bdesc)}$} (110:4) node[left]{\scriptsize  $\gamma(f)_{\nu^+(\Bdesc)}$}
(79:3.6) node[right]{\scriptsize $\gamma(\pointp(\alpha_i))$} (93:4.2) node[left]{\scriptsize  $\gamma(e)$} (95:4.2) node[right]{\scriptsize  $\gamma(e)_{\nu^+(\Bdesc)}$} (4,0) node[left]{\scriptsize  $\alpha_i$};
\end{scope}
\begin{scope}[xshift=4.4cm]
\draw [fill=gray!20, draw=white] (0,0) circle (4);
\draw [dashed, -latex] (-4,0) arc (180:360:4);
\draw [-latex] (108:1) -- (108:4) (0,0) --(108:1);
\draw [-latex] (0,0) -- (96:4);
\draw [-latex] (0,0) -- (78:4);
\draw [dashed, -latex] (.4,0) -- (105:4) ;
\draw [dashed, -latex] (.4,0) -- (99:4);
\draw [dashed, -latex] (.4,0) -- (72:4);
\draw (-.2,1) node[left]{\scriptsize  $\gamma(f)$} (-.2,-.2) node{\scriptsize  $a_i$}  (.2,-.2) node[right]{\scriptsize  $a_{i,\nu(\Bdesc)}$};
\draw (73:4) node[right]{\scriptsize  $\gamma(\pointp(\alpha_i))_{\nu^-(\Bdesc)}$} (102:4.1) node[left]{\scriptsize  $\gamma(f)_{\nu^-(\Bdesc)}$}
(104:4.3) node[right]{\scriptsize  $\gamma(e)_{\nu^-(\Bdesc)}$} (4,0) node[left]{\scriptsize  $\alpha_i$};
\end{scope}
\end{tikzpicture}
\caption{The $\gamma(c)_{\nu^+(\Bdesc)}$ and the $\gamma(c)_{\nu^-(\Bdesc)}$ near $a_i$ (where $\sigma(\pointp(\alpha_i))=\sigma(f)=1=-\sigma(e)$)}
\label{figNB}
\end{center}
\end{figure}

\begin{figure}[h]
\begin{center}
\begin{tikzpicture}[scale=.77]
\useasboundingbox (-5,-4.1) rectangle (5,4.3);
\begin{scope}[xshift=-4.4cm]
\draw [fill=gray!20, draw=white] (0,0) circle (4);
\draw [dashed, -latex] (-4,0) arc (180:360:4);
\draw [-<] (108:1) -- (108:4) (0,0) --(108:1);
\draw [-<] (0,0) -- (96:4);
\draw [-<] (0,0) -- (78:4);
\draw [dashed, -<] (.4,0) -- (113:4) ;
\draw [dashed, -<] (.4,0) -- (91:4);
\draw [dashed, -<] (.4,1) -- (.4,3.98) (.4,0) -- (.4,1);
\draw (-.2,1) node[left]{\scriptsize $\gamma(f)$} (-.2,-.2) node{\scriptsize $b_j$} (.8,-.2) node{\scriptsize $b_{j}^{\pm}$};
\draw (.4,1) node[right]{\scriptsize $\gamma(\pointp(\beta_j))_{\nu^+(\Asc)}$} (110:4) node[left]{\scriptsize $\gamma(f)_{\nu^+(\Asc)}$}
(79:3.6) node[right]{\scriptsize $\gamma(\pointp(\beta_j))$} (93:4.2) node[left]{\scriptsize $\gamma(e)$} (94:4.3) node[right]{\scriptsize $\gamma(e)_{\nu^+(\Asc)}$} (4,0) node[left]{\scriptsize $\beta_j$};
\end{scope}
\begin{scope}[xshift=4.4cm]
\draw [fill=gray!20, draw=white] (0,0) circle (4);
\draw [dashed, -latex] (-4,0) arc (180:360:4);
\draw [-<] (108:1) -- (108:4) (0,0) --(108:1);
\draw [-<] (0,0) -- (96:4);
\draw [-<] (0,0) -- (78:4);
\draw [dashed, -<] (.4,0) -- (105:4) ;
\draw [dashed, -<] (.4,0) -- (99:4);
\draw [dashed, -<] (.4,0) -- (72:4);
\draw (-.2,1) node[left]{\scriptsize $\gamma(f)$} (-.2,-.2) node{\scriptsize $b_j$} (.8,-.2) node{\scriptsize $b_{j}^{\pm}$};
\draw (73:4) node[right]{\scriptsize $\gamma(\pointp(\beta_j))_{\nu^-(\Asc)}$} (102:4.1) node[left]{\scriptsize $\gamma(f)_{\nu^-(\Asc)}$}
(102:4.3) node[right]{\scriptsize $\gamma(e)_{\nu^-(\Asc)}$} (4,0) node[left]{\scriptsize $\beta_j$};
\end{scope}
\end{tikzpicture}
\caption{The orthogonal projections of the $\gamma(d)_{\parallel}$ on $\Bdesc_j$ near $b_j$ (where $\sigma(\pointp(\beta_j))=\sigma(f)=1=-\sigma(e)$)}
\label{figNA}
\end{center}
\end{figure}
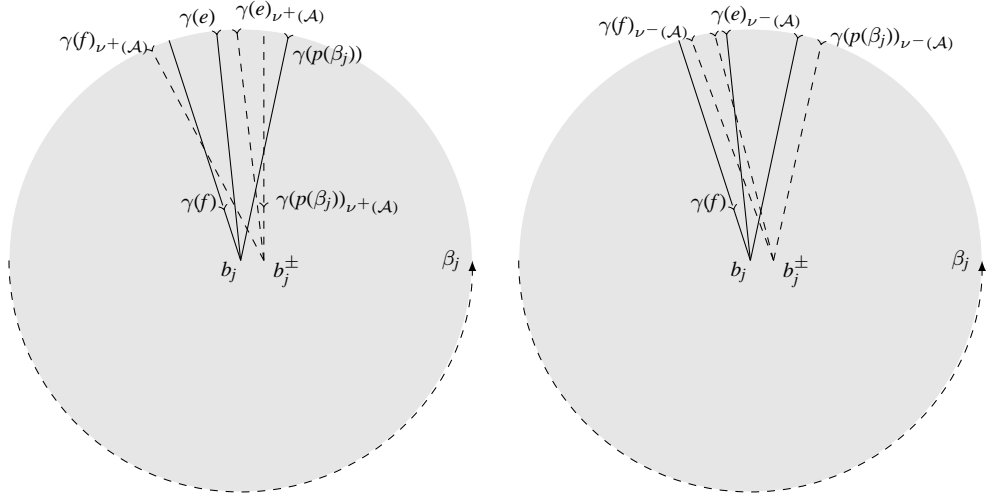

We shall similarly fix the positions of the
$$\gamma(d)_{\parallel}=\frac14\left(\gamma(d)_{\nu^+(\Asc)\nu^+(\Bdesc)} + \gamma(d)_{\nu^+(\Asc)\nu^-(\Bdesc)} +\gamma(d)_{\nu^-(\Asc)\nu^+(\Bdesc)} +\gamma(d)_{\nu^-(\Asc)\nu^-(\Bdesc)} \right)$$
by homotopies of the $\gamma(d)_{\nu^{\varepsilon}(\Asc)\nu^{\eta}(\Bdesc)}=\gamma_{\Asc}^{\varepsilon}(d) + [d^{\varepsilon}_{\Asc},d_{\Bdesc}^{\eta}] + \gamma_{\Bdesc}^{\eta}(d)$, with the notation of Subsection~\ref{subparflow}, so that:
\begin{itemize}
\item for any $d$, $\partial \gamma(d)_{\nu^{\varepsilon}(\Asc)\nu^{\eta}(\Bdesc)} = b_{j(d)}^{\eta} - a_{i(d)}^{\varepsilon}$ is fixed, 
\item $\gamma(d)_{\nu^{\varepsilon}(\Asc)\nu^{\eta}(\Bdesc)}$ is on the $\varepsilon$ side of $\Asc_{i(d)}$ except near $b_{j(d)}$ where its orthogonal projection $\gamma(d)_{\nu^{\varepsilon}(\Asc)}$ on $\Bdesc_{j(d)}$ is shown in Figure~\ref{figNA},
\item $\gamma(d)_{\nu^{\varepsilon}(\Asc)\nu^{\eta}(\Bdesc)}$ is on the $\eta$ side of $\Bdesc_{j(d)}$ except near $a_{i(d)}$ where its orthogonal projection on $\Asc_{i(d)}$ behaves like the projection of $\gamma(d)_{\nu^{\eta}(\Bdesc)}$ in Figure~\ref{figNB} at a larger scale.
\end{itemize}
In particular, the orthogonal projections on $\Bdesc_{j(d)}$ of $b_{j(d)}^+$ and $b_{j(d)}^-$ both coincide with the intersection point of the dashed segments in Figure~\ref{figNA}, and the orthogonal projections on $\Asc_{i(d)}$ of $a_{i(d)}^+$ and $a_{i(d)}^-$ both coincide with the intersection point of the dashed segments in Figure~\ref{figNB} at a larger scale.

These positions being fixed, we have the following proposition that implies Proposition~\ref{propeval}.

\begin{proposition}
\label{proplkplusmieux}

$$\langle \gamma(c)_{\nu(\Bdesc)} \times \gamma(d)_{\parallel},\prop(\funcf,\metrigo) \rangle=\ell(c,d).$$
\end{proposition}

Recall $\prop(\funcf,\metrigo)=\preprop_{\phi}+ \preprop_{\CI}$. We prove the proposition by computing the intersections with
$\preprop_{\CI}$ and $\preprop_{\phi}$ in Lemmas~\ref{lemintPI} and \ref{lemintPphi} below.
\begin{lemma}
\label{lemintPI}
$$\langle \gamma(c)_{\nu(\Bdesc)} \times \gamma(d)_{\parallel},\Bdesc_j \times \Asc_i \rangle=-\langle [\pointp(\alpha(c)),c\halfbra_{\alpha},\beta_j \rangle \langle \alpha_i ,[\pointp(\beta(d)),d\halfbra_{\beta} \rangle$$
\end{lemma}
\bp
In any case, $\langle \gamma(c)_{\nu(\Bdesc)} \times \gamma(d)_{\parallel},\Bdesc_j \times \Asc_i\rangle_{C_2(M)}=\langle \gamma(c)_{\nu(\Bdesc)},\Bdesc_j\rangle_{M} \langle \gamma(d)_{\parallel}, \Asc_i\rangle_{M}$.

The only intersection points of $\gamma(c)_{\nu(\Bdesc)}$ with $\Bdesc_j$
are shown in Figure~\ref{figNB}. Then since the $\gamma(c)_{\nu(\Bdesc)}$ cross the $\Bdesc_j$ like the $\alpha_i$, which are positive normals for $\Bdesc_j$ along flow lines associated to positive crossings $$\langle \gamma(c)_{\nu(\Bdesc)},\Bdesc_j\rangle_{M}=\langle [\pointp(\alpha(c)),c\halfbra_{\alpha(c)},\beta_j \rangle.$$

The computation of $\langle \gamma(d)_{\parallel}, \Asc_i\rangle_{M}$ is similar since the position of the $\gamma(d)_{\parallel}$ with respect to $\Bdesc_j$ does not matter. The only difference comes from the fact that the flow lines are oriented towards $b_{j(d)}$ so that they cross the $\Asc_i$ like $(-\beta_j)$, which is the positive normal along flow lines associated to negative crossings. See Figure~\ref{figNA}.
$$\langle \gamma(d)_{\parallel},\Asc_i\rangle_{M}=-\langle \alpha_i, [\pointp(\beta(d)),d\halfbra_{\beta(d)} \rangle.$$
\eop

\begin{lemma}
\label{lemintPphi}
$$\langle \gamma(c)_{\nu(\Bdesc)} \times \gamma(d)_{\parallel },\preprop_{\phi}\rangle=\langle [\pointp(\alpha(c)),c\halfbra_{\alpha}, [\pointp(\beta(d)),d\halfbra_{\beta}\rangle.$$
\end{lemma}
\bp Assume $c \in \alpha_{i} \cap \beta_{j(c)}$ and $d \in \alpha_{i(d)} \cap \beta_{j}$.
When the first $\check{M}$-coordinate of a point of $\preprop_{\phi}$ is in $\gamma(c) \setminus a_i$, its second $\check{M}$-coordinate is in $\left(\gamma(c) \cup \linh(b_{j(c)})\right)$, and therefore it is not in $\gamma(d)_{\parallel}$. 
Since the first $\check{M}$-coordinate of a point in $\gamma(c)_{\nu(\Bdesc)} \times \gamma(d)_{\parallel }$ is very close to $\gamma(c)$, $\gamma(c)_{\nu(\Bdesc)} \times \gamma(d)_{\parallel }$ intersects $\preprop_{\phi}$ in a small neighborhood of $a_i \times \Asc_i$.

Thus, the intersection points will be very close to pairs of points on flow rays from $a_i$ on $\Asc_i$, the closest point to $a_i$ being on
$\gamma(c)_{\nu(\Bdesc)}$ and the second one on $\gamma(d)_{\parallel}$. Then, for a given
$\gamma(c)$, the second point must be on the subsurface $D(\gamma(c))$ of $\Asc_i$ made of the points $x$ such that the flow ray from $a_i$ to $x$ intersects $\gamma(c)_{\nu^+(\Bdesc)}$ or $\gamma(c)_{\nu^-(\Bdesc)}$. This interaction locus of $\gamma(c)_{\nu^+(\Bdesc)}$, $D(\gamma(c))$, is shown in
Figure~\ref{figinterac1}.
The interaction locus of $\gamma(c)_{\nu^-(\Bdesc)}$ is similar.

\begin{figure}[h]
\begin{center}
\begin{tikzpicture}[scale=.6]
\useasboundingbox (-9,-2) rectangle (7,4.4);
\begin{scope}[xshift=2cm]
\draw [fill=gray!40, draw=white] (0,0) circle (4);
\draw [fill=white, draw=white] (-4.2,-4.2) rectangle (4.2,0);
\draw [fill=white, draw=white] (-4.2,0) -- (.4,0) -- (113:4) -- (-4.2,4.2) -- cycle;
\draw [dashed, -latex] (2,0) arc (0:360:2);
\draw [-latex] (108:1) -- (108:4) (0,0) --(108:1);
\draw [-latex] (0,0) -- (96:4);
\draw [-latex] (0,0) -- (78:4);
\draw [dashed, -latex] (.4,0) -- (113:4) ;
\draw [dashed, -latex] (.4,0) -- (91:4);
\draw [dashed, -latex] (.4,1) -- (.4,3.98) (.4,0) -- (.4,1);
\draw (.4,1) node[right]{\scriptsize $\gamma(\pointp(\alpha_i))_{\nu^+(\Bdesc)}$} (110:4) node[left]{\scriptsize $\gamma(f)_{\nu^+(\Bdesc)}$} (79:3.6) node[right]{\scriptsize $\gamma(\pointp(\alpha_i))$} (93:4.2) node[left]{\scriptsize $\gamma(e)$} (94:4.3) node[right]{\scriptsize $\gamma(e)_{\nu^+(\Bdesc)}$} (2,0) node[right]{\scriptsize $\alpha_i$};
\end{scope}
\begin{scope}[xshift=-6cm]
\draw [fill=gray!40, draw=white] (.4,0) -- (91:4) arc (91:0:4) -- cycle;
\draw [dashed, -latex] (2,0) arc (0:360:2);
\draw [-latex] (108:1) -- (108:4) (0,0) --(108:1);
\draw [-latex] (0,0) -- (96:4);
\draw [-latex] (0,0) -- (78:4);
\draw [dashed, -latex] (.4,0) -- (113:4) ;
\draw [dashed, -latex] (.4,0) -- (91:4);
\draw [dashed, -latex] (.4,1) -- (.4,3.98) (.4,0) -- (.4,1);
\draw (-.2,1) node[left]{\scriptsize $\gamma(f)$};
\draw (93:4.2) node[left]{\scriptsize $\gamma(e)$} (94:4.3) node[right]{\scriptsize $\gamma(e)_{\nu^+(\Bdesc)}$} (2,0) node[right]{\scriptsize $\alpha_i$};\end{scope}
\end{tikzpicture}
\caption{Interaction loci of $\gamma(e)_{\nu^+(\Bdesc)}$ and $\gamma(f)_{\nu^+(\Bdesc)}$ on $\Asc_i$ (where $\sigma(f)=1=-\sigma(e)$)}
\label{figinterac1}
\end{center}
\end{figure}

The only intersection points of $\gamma(d)_{\parallel }$ with the domain $D(\gamma(c))$ of $\Asc_i$ are near the $b_j$ and they are shown in Figure~\ref{figNA}.

The curve $\gamma(d)_{\parallel}$ meets $\Asc_i$ near a crossing line $\gamma(e)$, where {\em near\/} means
in the sheet of $\gamma(e)$ around $\linh(b_j)$,
\begin{itemize}
 \item with probability $1$ if $i(e)=i$ and if $e \in [\pointp(\beta_j),d[_{\beta_j}$,
\item with probability $1/2$ (depending on the side of $\Asc_i$ for $\gamma(d)_{\parallel }$ near $b_j$) if $i(e)=i$ and if $e=d$, (this is also valid when $e=\pointp(\beta_j)=d$),
\item  with probability $0$ in the other cases.
\end{itemize}

The corresponding intersection point is in $D(\gamma(c))$
if $e \in [\pointp(\alpha_i),c[_{\alpha_i}$, or if $e=c$ and $\gamma(d)_{\parallel }$ is on the correct side of $\Bdesc_j$ (the $(-\alpha_i)$ side), that is with a probability $1/2$ independent of the previous one.

Then $M$ is oriented as $(\mbox{flow line} \times \gamma(c)_{\nu(\Bdesc)} \times \nu^+(\Asc_i))$ near $a_i$ and $\preprop_{\phi}$ is oriented as $$(\mbox{beginning of flow line} \times \mbox{diag}(\gamma(c)_{\nu(\Bdesc)} \times \nu^+(\Asc_i)) \times \mbox{end of flow line}),$$
which is intersected negatively by $\gamma(c)_{\nu(\Bdesc)} \times \nu^+(\Asc_i)$, where $\nu^+(\Asc_i)$ is oriented like $\sigma(e)\beta_j$ and like $(-\sigma(e))\gamma(d)_{\parallel}$ near a point in $\left(\gamma(c)_{\nu(\Bdesc)} \times \gamma(d)_{\parallel}\right) \cap \preprop_{\phi}$ corresponding to a crossing $e$ of $[\pointp(\alpha(c)),c\halfbra_{\alpha} \cap [\pointp(\beta(d)),d\halfbra_{\beta}$.
\eop

\section{The combing associated with \texorpdfstring{$\matchingo$}{m} and its associated propagator}
\label{seccombCP}
\setcounter{equation}{0}

In this section, we first define the combing $\ComX(\pointw,\matchingo)$ of $\check{M}$. Next we introduce correction $4$-chains $\preprop_{\homotop}$ and $\preprop_{\Sigma}$ in $U\check{M} \subset \partial C_2(M)$ such that the sum $\prop=\prop(\funcf,\metrigo)+\preprop_{\homotop}+\preprop_{\Sigma}$ is a propagator associated with $\ComX(\pointw,\matchingo)$.

\subsection{The combing \texorpdfstring{$\ComX(\pointw,\matchingo)$}{X(w,m)}}
\label{subframedsel}
Consider the matching $\matchingo$ introduced in Subsection~\ref{subframedHd}.
Up to renumbering and reorienting the $\Bdesc_j$, assume that $c_i \in \alpha_i \cap \beta_i$ and that $\sigma(c_i)=1$. Set $\gamma_i=\gamma(c_i)$.

There is a combing $\ComX=\ComX(\pointw,\matchingo)$\index{Nt}{Xwm@$\ComX(\pointw,\matchingo)$} (section of the unit tangent bundle) of $\check{M}$ that coincides with the direction $s_{\phi}$ of the flow (and the gradient of $\fMorse$) outside the union of regular neighborhoods $N(\gamma_i)$ of
the $\gamma_i$, that is opposite to $s_{\phi}$ along the interiors of the $\gamma_i$ and
that is obtained as follows on $N(\gamma_i)$.
Choose a natural trivialization $(X_1,X_2,X_3)$ of $T\check{M}$ on a regular neighborhood $N(\gamma_i)$ of $\gamma_i$, such that:
\begin{itemize}
\item $\gamma_i$ is directed by $X_1$,
\item the other flow lines never have $X_1$ as an oriented tangent vector,
\item $(X_1,X_2)$ is tangent to $\Asc_i$ (except on the parts of $\Asc_i$ near $b_i$ that
come from other crossings of $\alpha_i \cap \beta_i$), and $(X_1,X_3)$ is tangent to $\Bdesc_i$ (except on the parts of $\Bdesc_i$ near $a_i$ that
come from other crossings of $\alpha_i \cap \beta_i$).
\end{itemize}
This parallelization identifies the unit tangent bundle $UN(\gamma_i)$ of $N(\gamma_i)$ with $S^2 \times N(\gamma_i)$.

\medskip
\noindent There is a homotopy  $\homotop\colon [0,1]\times (N(\gamma_i) \setminus \gamma_i) \rightarrow S^2$,  such that
\begin{itemize}
\item
$\homotop_0$ is the unit tangent vector to the flow lines of $\phi$,
\item
$\homotop_1$ is the constant map to $(-X_1)$ and \item
$\homotop_t(y)$ goes from $\homotop_0(y)=s_{\phi}(y)$ to $(-X_1)$ along the shortest geodesic arc of $S^2$ from $s_{\phi}(y)$ to $(-X_1)$, which is denoted by $[s_{\phi}(y),-X_1]$.
\end{itemize}
Let $2\eta$ be the distance between $\gamma_i$ and $\partial N(\gamma_i)$ and
let $\ComX(y) = \homotop(\mbox{max}(0,1-d(y,\gamma_i)/\eta),y)$ on $N(\gamma_i)\setminus \gamma_i$, and
$\ComX=-\ComX_1$ along $\gamma_i$.

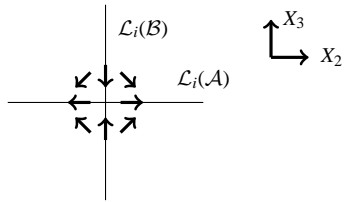
\begin{figure}[h]
\begin{center}
\begin{tikzpicture}
\useasboundingbox (-1.4,-1.4) rectangle (1.4,1.4);
\draw (-1.3,0) -- (1.3,0) (0,-1.3) -- node[very near end, right]{\scriptsize $\linh_i(\Bdesc)$} (0,1.3);
\draw [very thick,->] (0,-.5) -- (0,-.2);
\draw [very thick,->] (0,.5) -- (0,.2);
\draw [very thick,->] (-.2,0) -- (-.5,0);
\draw [very thick,->] (.2,0) -- (.5,0);
\draw [very thick,->] (.2,.4) -- (.4,.2);
\draw [very thick,->] (.2,-.4) -- (.4,-.2);
\draw [very thick,->] (-.2,.4) -- (-.4,.2);
\draw [very thick,->] (-.2,-.4) -- (-.4,-.2);
\draw [very thick,->] (2.2,.6) -- (2.2,1.1);
\draw [very thick,->] (2.2,.6) -- (2.7,.6);
\draw (1.3,.3) node{\scriptsize $\linh_i(\Asc)$} (2.7,.6) node[right]{\scriptsize $X_2$} (2.2,1.1) node[right]{\scriptsize $X_3$};
\end{tikzpicture}
\caption{Projection of $\ComX$}
\label{figprojX}
\end{center}
\end{figure}

Note that $\ComX$ is tangent to $\Asc_i$ on $N(\gamma_i)$ (except on the parts of $\Asc_i$ near $b_i$ that
come from other crossings of $\alpha_i \cap \beta_i$), and that $\ComX$ is tangent to $\Bdesc_i$ on $N(\gamma_i)$ (except on the parts of $\Bdesc_i$ near $a_i$ that
come from other crossings of $\alpha_i \cap \beta_i$).
More generally, project the normal bundle to $\gamma_i$ to $\RR^2$ in the $\ComX_1$--direction by sending $\gamma_i$ to $0$, $\Asc_i$ to an axis $\linh_i(\Asc)$ and $\Bdesc_i$ to an axis $\linh_i(\Bdesc)$. Then the projection of $\ComX$ goes towards $0$ along $\linh_i(\Bdesc)$ and starts from $0$ along $\linh_i(\Asc)$, it has the direction of $\sigma_a(y)$ at a point $y$ of $\RR^2$ near $0$, where
$\sigma_a$\index{Nt}{sa@$\sigma_a$} is the planar reflexion that fixes $\linh_i(\Asc)$ and reverses $\linh_i(\Bdesc)$. See Figure~\ref{figprojX}.

Then $\ComX(y)$ is on the half great circle that contains $\sigma_a(y)$, $X_1$ and $(-X_1)$.
In Figure~\ref{figgammai} (and in Figure~\ref{figHa}), $\gamma_i$ is a vertical segment, all the other flow lines corresponding to crossings involving $\alpha_i$ go upwards from $a_i$, and $\ComX$ is simply the upward vertical field. See also Figure~\ref{figinthandles}.

\begin{figure}[h]
\begin{center}
\begin{tikzpicture}
\useasboundingbox (1,1.4) rectangle (5,3.3);
\draw (1.3,1.5) .. controls (2,1.5) and (2,3.2) .. (3,3.2);
\draw (4.7,1.5) .. controls (4,1.5) and (4,3.2) .. (3,3.2) ;
\draw [thick,->,fill=gray!20,draw=blue] (2.5,2) arc (180:0:.5);
\draw (2.5,2) arc (-180:0:.5);
\draw [blue] (3.5,2) node[right]{\scriptsize $\beta_i$};
\draw [red,thick,->] (3,3.2) .. controls (2.9,3.2) and (2.8,3) .. (2.8,2.85);
\draw [red,thick] (3,2.5) .. controls (2.9,2.5) and (2.8,2.7) .. (2.8,2.85) (2.6,2.85) node{\scriptsize $\alpha_i$};
\draw [red,thick,dashed] (3,3.2) .. controls (3.1,3.2) and (3.2,3) .. (3.2,2.85) .. controls (3.2,2.7) and (3.1,2.5) .. (3,2.5);
\draw [->] (3,2.85) -- node[very near end,right=.1]{\scriptsize $\gamma_i$} (3,2);
\end{tikzpicture}
\caption{$\gamma_i$}
\label{figgammai}
\end{center}
\end{figure}

\subsection{The propagator associated with a combed Heegaard splitting}
\label{subFcomb}

Recall that $UN(\gamma_i)$ is identified with $S^2 \times N(\gamma_i)$.
Let
$\preprop_{\homotop}=\preprop_{\homotop}(\matchingo)$\index{Nt}{Ph@$\preprop_{\homotop}$} be the closure in $\partial C_2(M)$ of the image of 
$\{(t,y); t \in [0,\mbox{max}(0,1-d(y,\gamma_i)/\eta)], y \in N(\gamma_i)\setminus \gamma_i\}$
in $S^2 \times (N(\gamma_i))$ under $\left((t,y) \mapsto (\homotop(t,y),y)\right)$.

\begin{lemma}
$\partial \preprop_{\homotop}=\overline{\ComX(\check{M})}-\overline{s_{\phi}(\check{M})} -\sum_{i=1}^g U\check{M}_{|\gamma_i} $
\end{lemma}
\bp We explain the $(U\check{M}_{|\gamma_i}= S^2\times\gamma_i )$ part of $\partial \preprop_{\homotop}$, with its sign.
The homotopy $\homotop$ naturally extends to $[0,1] \times B\!\ell(N(\gamma_i),\gamma_i)$, where $B\!\ell(N(\gamma_i),\gamma_i)$ is obtained from $N(\gamma_i)$ by blowing up $\gamma_i$, so that $(-B\!\ell(N(\gamma_i),\gamma_i))$
contains the unit normal bundle $S^1 \times \gamma_i$ of $\gamma_i$ in $C_M$, in its boundary. Then $\partial \preprop_{\homotop}$ contains
$\{(\homotop(t,y),p_{\gamma_i}(y)) \in S^2 \times \gamma_i; t \in [0,1], y \in S^1 \times \gamma_i\}$, where $S^1$, which is the blown-up center of the fiber $D^2$ of $N(\gamma_i)$, is mapped by $\sigma_a$ to the equator of $S^2$ so that the image of $([0,1] \times S^1)$ covers a fiber $S^2$ of $U\check{M}_{|\gamma_i}$ with degree $(-1)$.
\eop

Recall the $1$--cycle $\Link(\matchingo) = \sum_{i=1}^g \gamma_i- \sum_{c \in \CaC} \CJ_{j(c)i(c)} \sigma(c) \gamma(c)$.
Let $\Sigma(\matchingo)$ be a two--chain bounded by $\Link(\matchingo)$ in $\check{M}$ and
 let $$\preprop_{\Sigma}=U\check{M}_{|\Sigma(\matchingo)}.$$\index{Nt}{PSigma@$\preprop_{\Sigma}$}
Note that $\preprop_{\Sigma}$ is homeomorphic to $S^2 \times \Sigma(\matchingo)$.

\begin{proposition}
 $$\prop=\prop(\funcf,\metrigo)+\preprop_{\homotop}+\preprop_{\Sigma}$$
is a propagator associated with the combing $\ComX(\pointw,\matchingo)$.
\end{proposition}
\bp 
The boundary of $\prop$
is $(\overline{\ComX(\pointw,\matchingo)(\check{M})}+\partial_{od})$. \eop

Recall that $\iota$ denotes the involution of $C_2(M)$ that exchanges two points in a pair. Then
$\iota(\prop)$ is also a propagator associated with the combing $(-\ComX(\pointw,\matchingo))$. Theorem~\ref{thmdefinvcomb} defines
$\Theta(M,\ComX(\pointw,\matchingo))$ from the algebraic intersection of $\prop$ and $\iota(\prop)$, which we compute from now on in order to prove Theorem~\ref{thmmain}.

\section{Computation of \texorpdfstring{$[\prop_{\ComX(\pointw,\matchingo)} \cap \prop_{-\ComX(\pointw,\matchingo)}]$}{the propagators' intersection}}
\label{seccompint}
\setcounter{equation}{0}

\subsection{A description of \texorpdfstring{$[\prop_{\ComX(\pointw,\matchingo)} \cap \prop_{-\ComX(\pointw,\matchingo)}]$}{the propagators' intersection}}
\label{subdescY}

Fix $\pointw$, $\matchingo$, $\ComX=\ComX(\pointw,\matchingo)$, $\Link=\Link(\matchingo)$ and $\Sigma=\Sigma(\matchingo)$ such that $\partial \Sigma= \Link$.

Consider a vector field $\ComY$ of $\ComX^{\perp}$ on $\check{M}$ such that
\begin{itemize}
\item $\ComY$ vanishes outside $C_M$,
\item the norm of $\ComY$ is one on the $\gamma(c)$,
\item for every $i$, $\ComY(a_i)$ is tangent to the line $\linh(a_i)$, which is the descending manifold of $a_i$, (but $\ComY(a_i)$ does not necessarily direct the line),
\item for every $j$, $\ComY(b_j)$ is tangent to the line $\linh(b_j)$, (again, $\ComY(b_j)$ does not necessarily direct the line),
\end{itemize}

Then $\Link_{\parallel \ComY}$ denotes the link parallel to $\Link$ obtained by pushing $\Link$ in the $\ComY$ direction.
Along $\gamma(c)$, $\sigma_a$ is the symmetry of $\ComX^{\perp}$ with respect to $\Asc_{i(c)}$ that preserves the vectors
tangent to $\Asc_{i(c)}$ and reverses the vectors tangent to $\Bdesc_{j(c)}$.
Define $\gamma(c) \times \gamma(d)_{\parallel \sigma_a(-\ComY)}$ as the product of $\gamma(c)$ and a parallel of $\gamma(d)$ ``infinitely'' close to $\gamma(d)$ in the direction of $\sigma_a(-\ComY)$.
This can be formalised as follows.
When $c\neq d$, $\gamma(c) \times \gamma(d)_{\parallel \sigma_a(-\ComY)}=\gamma(c) \times \gamma(d)$ (away from the possibly coinciding ends). For $x \in \gamma(c)$, let $\gamma^{\prime}_x(c)$ denote the unit tangent vector of $\gamma(c)$ at $x$ that orients $\gamma(c)$, and let $[-\gamma^{\prime}_x(c),\gamma^{\prime}_x(c)]_{\sigma_a(-\ComY)}$ denote the half great circle in the fiber $U\check{M}_{|x}$ through $\sigma_a(-\ComY(x))$ towards $\gamma^{\prime}_x(c)$. Let $s_{[-\gamma^{\prime}(c),\gamma^{\prime}(c)]_{\sigma_a(-\ComY)}}(\gamma(c))$ be the total space of the bundle over $\gamma(c)$ of these half-circles. 
Then $$\gamma(c) \times \gamma(c)_{\parallel \sigma_a(-\ComY)}=\overline{\gamma(c)^2 \setminus \mbox{diag}(\gamma(c)^2)} - s_{[-\gamma^{\prime}(c),\gamma^{\prime}(c)]_{\sigma_a(-\ComY)}}(\gamma(c)).$$
Similarly, $s_{[-\ComX,\ComX]_{\sigma_a(-\ComY)}}(\partial \Sigma)$ is the total space of the bundle over $\partial \Sigma$ of the half-circles $[-\ComX,\ComX]_{\sigma_a(-\ComY)}$.
In this section, we prove the following proposition.

\begin{proposition}
\label{propcycintdouble}
Let $\ComY$ be a vector field of $\ComX^{\perp}$ as above.
There exists a two-chain $O(\sigma_a(-\ComY))$ in the hemispheres centered at $\sigma_a(-\ComY)$ in $U\check{M}_{|\cup_i a_i \cup (\cup_j b_j)}$ such that
\index{Nt}{Gdui@$\Gdu^i(\ComY)$}$$\begin{array}{ll}\Gdu^i(\ComY)=&\sum_{(i,j,k,\ell) \in \{1,\dots,g\}^4}\CJ_{ji}\CJ_{\ell k}\left((\Bdesc_j\cap \Asc_{k}) \times (\Bdesc_{\ell}\cap \Asc_i)_{\parallel \sigma_a(-\ComY)}\right) \\
&-\sum_{c \in \CaC} \CJ_{j(c)i(c)}\sigma(c)\left(\gamma(c) \times \gamma(c)_{\parallel \sigma_a(-\ComY)}\right)\\&+O(\sigma_a(-\ComY))\end{array}
$$ is a $2$--cycle of $C_2(M)$ whose homology class is unambiguously defined.
Let $S$ be a fiber of $U\check{M}$ and let $\ComX(\Sigma)$ denote the graph of $\ComX_{|\Sigma}$ in $U\check{M}$. Set
\index{Nt}{GdubXY@$\Gdu^b(\ComX,\ComY)$}$$\Gdu^b(\ComX,\ComY)=lk (\Link,\Link_{\parallel \ComY})S-\left(\ComX(\Sigma) - (-\ComX)(\Sigma) -s_{[-\ComX,\ComX]_{\sigma_a(-\ComY)}}(\partial \Sigma)\right).$$
Then the cycle $$\Gdu=\Gdu^i(\ComY) + \Gdu^b(\ComX,\ComY)$$\index{Nt}{Gdu@$\Gdu$}
represents the homology class of $\prop_{\ComX(\pointw,\matchingo)} \cap \prop_{-\ComX(\pointw,\matchingo)}$.
\end{proposition}

\subsection{Introduction to specific chains \texorpdfstring{$\prop_{\ComX}$}{P(X)} and \texorpdfstring{$\prop_{-\ComX}$}{P(-X)}}
\label{subdeformprop}

In this subsection, we deform the propagators $\prop$ and $\iota(\prop)$ constructed in Section~\ref{subFcomb} to propagators $\prop_{\ComX}$ and $\prop_{-\ComX}$ that are transverse to each other, in order to determine their algebraic intersection.

Let $[-1,0] \times \partial C_2(M)$ be a (topological) collar of $\partial C_2(M)$ in $C_2(M)$.
Then $C_2(M)$ is homeomorphic to $\tilde{C}_2(M) = C_2(M) \setminus (]-1/2,0] \times \partial C_2(M))$ by the {\em shrinking homeomorphism\/} 
$$\begin{array}{llll} h_s\colon &C_2(M) &\rightarrow &\tilde{C}_2(M)\\
   &(t,x) \in [-1,0] \times \partial C_2(M) & \mapsto & ((t-1)/2,x) \in [-1,-1/2] \times \partial C_2(M)
  \end{array}
$$
that is the identity map outside the collar. 
Identifying $[-1/2,0]$ with $[0,6]$ by the appropriate affine monotonous transformation identifies $C_2(M)$ with $$\tilde{C}_2(M) \cup_{\partial \tilde{C}_2(M)} \left([0,6]\times \partial C_2(M)\right),$$ which is our space $C_2(M)$ from now on.

Use $h_s$ to shrink $\prop(\funcf,\metrigo)$ and $\iota(\prop(\funcf,\metrigo))$ into $\tilde{C}_2(M)$, and construct transverse $\prop_{\ComX}$ and $\prop_{-\ComX}$ with respective boundaries $\{6\} \times \partial \prop_{\ComX}$ and $\{6\} \times \partial \prop_{-\ComX}$ as follows:

$$\begin{array}{ll}\prop_{-\ComX}=&h_s(\iota(\prop(\funcf,\metrigo))) +[0,1] \times \partial \iota(\prop(\funcf,\metrigo)) \\
&+\{1\} \times \iota(\preprop_{\homotop}) +[1,3] \times (\iota(-S^2 \times \Link + \partial_{od}) + \overline{(-\ComX)(\check{M})})\\
&+\{3\} \times \iota(S^2 \times \Sigma) +[3,6] \times (\overline{(-\ComX)(\check{M})}+\iota(\partial_{od}))\end{array}$$
while the following expression of $\prop_{\ComX}$, which is partially schematically drawn in Figure~\ref{figschematic}, will require a perturbating diffeomorphism $\Psi$ of $C_2(M)$ isotopic and very close to the identity map in order to get transversality near the diagonal.
$$\begin{array}{ll}\prop_{\ComX}=&h_s(\Psi(\prop(\funcf,\metrigo))) +[0,2] \times \partial \Psi(\prop(\funcf,\metrigo))\\& +\{2\} \times \Psi(\preprop_{\homotop})
+[2,4] \times \Psi(-S^2 \times \Link + \overline{\ComX(\check{M})} + \partial_{od})
\\& +\{4\} \times \Psi(S^2 \times \Sigma )+[4,5] \times \Psi( \overline{\ComX(\check{M})} + \partial_{od}) \\&+ \{5\} \times \Psi_{[\varepsilon,0]}( \partial \prop_{\ComX}) + [5,6] \times \partial \prop_{\ComX} \end{array}$$
where $\Psi_{[\varepsilon,0]}( \partial \prop_{\ComX})$ is the small cobordism between $\Psi( \overline{\ComX(\check{M})} + \partial_{od})$
and $\partial \prop_{\ComX}$ induced by the isotopy between $\Psi$ and the identity map.
We describe $\Psi$ in the next subsection.

\begin{figure}[h]
\begin{center}
\begin{tikzpicture}
\useasboundingbox (1,0) rectangle (14.1,11.8);
\draw (1.4,1.6) -- (1,1.2) (1,.2) -- (6,.2) -- (7.6,1.8) -- (7.4,1.8);
\fill[gray!20] (1.4,1.6) -- (2,1.6) .. controls (2.3,1.6) .. (2.4,.9) .. controls (3,.4) .. (3.4,.6) ..  controls (3.58,.75) .. (3.6,.9) .. controls  (3.62,1.05) .. (3.7,1.2) -- (4.7,1.2) .. controls (5.2,1.2) and (6,1.6) .. (6.6,1.6) -- (7.4,1.6) -- (7.4,4.8) -- (6.6,4.8) .. controls (6,4.8) and (5.2,4.4) .. (4.7,4.4) -- (3.7,4.4) ..  controls (3.62,4.25) .. (3.6,4.1)  .. controls (3.58,3.95) .. (3.4,3.8) ..  controls (3,3.6) ..  (2.4,4.1) .. controls (2.3,4.8) .. (2,4.8) -- (1.4,4.8);
\fill[pattern=vertical lines, pattern color=gray!50] (1.4,1.6) -- (2,1.6) .. controls (2.3,1.6) .. (2.4,.9) .. controls (3,.4) .. (3.4,.6) ..  controls (3.58,.75) .. (3.6,.9) .. controls  (3.62,1.05) .. (3.7,1.2) -- (4.7,1.2) .. controls (5.2,1.2) and (6,1.6) .. (6.6,1.6) -- (7.4,1.6) -- (7.4,4.8) -- (6.6,4.8) .. controls (6,4.8) and (5.2,4.4) .. (4.7,4.4) -- (3.7,4.4) ..  controls (3.62,4.25) .. (3.6,4.1)  .. controls (3.58,3.95) .. (3.4,3.8) ..  controls (3,3.6) ..  (2.4,4.1) .. controls (2.3,4.8) .. (2,4.8) -- (1.4,4.8);
\draw[thick,->] (1.4,1.6) -- (1.7,1.6) (2,1.6) .. controls (2.3,1.6) .. (2.4,.9) (2,1.6) -- (1.7,1.6);
\draw (1.7,1.4) node{\scriptsize $\partial_{od}$};
\draw[thick,<-] (2.4,.9) .. controls (3,.4) .. (3.4,.6) ..  controls (3.58,.75) .. (3.6,.9);
\draw (2.5,.8) node[left]{\scriptsize $\Psi(s_{\phi}(\check{M}))$};
\draw[thick,>-] (3.6,.9) .. controls  (3.62,1.05) .. (3.7,1.2) -- (4.7,1.2);
\draw (3.55,.8) node[right]{\scriptsize $\Psi(U\check{M}_{|L \setminus (\cup_i \gamma_i)})$};
\draw[thick,<-] (4.7,1.2) .. controls (5.2,1.2) and (6,1.6) .. (6.6,1.6);
\draw (4.7,1.4) node{\scriptsize $\Psi(s_{\phi}(\check{M}))$};
\draw[thick,<-] (6.6,1.6) -- (7,1.6) (7,1.6) -- (7.4,1.6);
\draw (6.85,1.4) node{\scriptsize $\partial_{od}$};
\fill (6.6,1.6) circle (0.065) (2,1.6) circle (0.065) (3.4,.6) circle (0.065) (3.7,1.2) circle (0.065);
\begin{scope}[yshift=1.6 cm]
\draw[gray!60,thick,->] (1.4,1.6) -- (1.7,1.6) (2,1.6) .. controls (2.3,1.6) .. (2.4,.9) (2,1.6) -- (1.7,1.6);
\draw[gray!60,thick,<-] (2.4,.9) .. controls (3,.4) .. (3.4,.6) ..  controls (3.58,.75) .. (3.6,.9);
\draw[gray!90] (2.4,.9) node[left]{\scriptsize $[0,2] \times $};
\draw[gray!60,thick,<-] (3.6,.9) .. controls  (3.62,1.05) .. (3.7,1.2) -- (4.7,1.2);
\draw[gray!60,thick,<-] (4.7,1.2) .. controls (5.2,1.2) and (6,1.6) .. (6.6,1.6);
\draw[gray!60,thick,<-] (6.6,1.6) -- (7,1.6) (7,1.6) -- (7.4,1.6);
\fill[gray!60] (6.6,1.6) circle (0.065) (2,1.6) circle (0.065) (3.4,.6) circle (0.065) (3.7,1.2) circle (0.065);
\end{scope}
\begin{scope}[yshift=3.2 cm,xshift=6.5cm]
\draw (1,.2) -- (6,.2) -- (7.6,1.8) -- (1.6,1.8) -- (1,1.4);
\fill[gray!30] (2,1.6) .. controls (2.3,1.6) and (3.6,1.4) .. (4,1.4) .. controls (4.3,1.4) and (4.6,1.2) .. (4.7,1.2) -- (3.7,1.2) .. controls (3.5,1.2) and (3.2,1) .. (3.2,.9) .. controls (3.2,.6) .. (3.4,.6) .. controls (3,.4) .. (2.4,.9) .. controls (2.3,1.6) ..(2,1.6);
\draw (2.8,1.2) node{\scriptsize $\Psi(\preprop_{\homotop})$};
\draw (2.4,1.2) node[left]{\scriptsize $\{2\} \times $};
\draw[thick] (2,1.6) .. controls (2.3,1.6) .. (2.4,.9);
\draw[thick,-] (2.4,.9) .. controls (3,.4) .. (3.4,.6);
\draw[thick,>-] (3.2,.9) .. controls (3.2,1) and (3.5,1.2) .. (3.7,1.2)  (3.2,.9) .. controls (3.2,.6) .. (3.4,.6);
\draw (3.2,.8) node[right]{\scriptsize $\Psi(U\check{M}_{|\cup_i \gamma_i})$};
\draw[thick] (3.7,1.2) -- (4.7,1.2);
\draw[thick,-] (2,1.6) .. controls (2.3,1.6) and (3.6,1.4) .. (4,1.4) (4.7,1.2) .. controls (4.6,1.2) and (4.3,1.4) .. (4,1.4) ;
\fill (3.4,.6) circle (0.065) (3.7,1.2) circle (0.065);
\end{scope}
\begin{scope}[yshift=3.2 cm]
\fill[gray!20] (1.4,1.6) -- (2,1.6) .. controls (2.3,1.6) and (3.6,1.4) .. (4,1.4) .. controls (4.3,1.4) and (4.6,1.2) .. (4.7,1.2) .. controls (5.2,1.2) and (6,1.6) .. (6.6,1.6) -- (7.4,1.6) -- (7.4,6.4) -- (6.6,6.4) .. controls (6,6.4) and (5.2,6) .. (4.7,6) .. controls (4.6,6) and (4.3,6.2) .. (4,6.2) .. controls (3.6,6.2) and (2.3,6.4)   .. (2,6.4) -- (1.4,6.4);
\fill[pattern=vertical lines, pattern color=gray!50] (1.4,1.6) -- (2,1.6) .. controls (2.3,1.6) and (3.6,1.4) .. (4,1.4) .. controls (4.3,1.4) and (4.6,1.2) .. (4.7,1.2) .. controls (5.2,1.2) and (6,1.6) .. (6.6,1.6) -- (7.4,1.6) -- (7.4,6.4) -- (6.6,6.4) .. controls (6,6.4) and (5.2,6) .. (4.7,6) .. controls (4.6,6) and (4.3,6.2) .. (4,6.2) .. controls (3.6,6.2) and (2.3,6.4)   .. (2,6.4) -- (1.4,6.4);
\fill[gray!40] (2,1.6) .. controls (2.3,1.6) and (3.6,1.4) .. (4,1.4) .. controls (4.3,1.4) and (4.6,1.2) .. (4.7,1.2) -- (3.7,1.2) .. controls (3.5,1.2) and (3.2,1) .. (3.2,.9) .. controls (3.2,.6) .. (3.4,.6) .. controls (3,.4) .. (2.4,.9) .. controls (2.3,1.6) ..(2,1.6);
\fill[pattern=north east lines, pattern color=gray!70] (2,1.6) .. controls (2.3,1.6) and (3.6,1.4) .. (4,1.4) .. controls (4.3,1.4) and (4.6,1.2) .. (4.7,1.2) -- (3.7,1.2) .. controls (3.5,1.2) and (3.2,1) .. (3.2,.9) .. controls (3.2,.6) .. (3.4,.6) .. controls (3,.4) .. (2.4,.9) .. controls (2.3,1.6) ..(2,1.6);
\fill[gray!20] (7.4,6.4) -- (1.4,6.4) -- (1.4,8) -- (7.4,8) -- (7.4,6.4);
\fill[pattern=vertical lines, pattern color=gray!50] (7.4,6.4) -- (1.4,6.4) -- (1.4,8) -- (7.4,8) -- (7.4,6.4);
\draw[thick] (2,6.4) -- (6.6,6.4);
\draw[gray!70,thick] (1.4,6.4) -- (2,6.4) (6.6,6.4) -- (7.4,6.4);
\draw[thick,->] (4.4,8) -- (1.4,8) (7.4,8) -- (4.4,8);
\draw (4.4,7.9) node[above]{\scriptsize $\partial \prop_{\ComX}$ };
\draw[gray!90] (4,7.3) node[below]{\scriptsize $[5,6] \times $};
\draw[gray!70,thick,->] (4.4,7.2) -- (1.4,7.2) (7.4,7.2) -- (4.4,7.2);
\fill[gray!40] (6.6,6.4) .. controls (6,6.4) and (5.2,6) .. (4.7,6) .. controls (4.6,6) and (4.3,6.2) .. (4,6.2) .. controls (3.6,6.2) and (2.3,6.4)   .. (2,6.4) -- (6.6,6.4);
\fill[pattern=north east lines, pattern color=gray!70] (6.6,6.4) .. controls (6,6.4) and (5.2,6) .. (4.7,6) .. controls (4.6,6) and (4.3,6.2) .. (4,6.2) .. controls (3.6,6.2) and (2.3,6.4)   .. (2,6.4) -- (6.6,6.4);
\draw[thick](6.6,6.4) .. controls (6,6.4) and (5.2,6) .. (4.7,6) .. controls (4.6,6) and (4.3,6.2) .. (4,6.2) .. controls (3.6,6.2) and (2.3,6.4)   .. (2,6.4) -- (6.6,6.4);
\draw[thick,-](2,1.6) .. controls (2.3,1.6) .. (2.4,.9);
\draw[gray!70, thick] (1.4,1.6) -- (2,1.6);
\fill[gray!50,semitransparent] (3.2,.9) .. controls (3.2,.6) .. (3.4,.6) .. controls (3.58,.75) .. (3.6,.9) .. controls  (3.62,1.05) .. (3.7,1.2) -- (3.7,4.4) .. controls  (3.62,4.25) .. (3.6,4.1) ..  controls (3.58,3.95) .. (3.4,3.8) .. controls (3.2,3.8) .. (3.2,4.1) -- (3.2,.9);
\fill[pattern=vertical lines, pattern color=gray!60,semitransparent] (3.2,.9) .. controls (3.2,.6) .. (3.4,.6) .. controls (3.58,.75) .. (3.6,.9) .. controls  (3.62,1.05) .. (3.7,1.2) -- (3.7,4.4) .. controls  (3.62,4.25) .. (3.6,4.1) ..  controls (3.58,3.95) .. (3.4,3.8) .. controls (3.2,3.8) .. (3.2,4.1) -- (3.2,.9);
\draw[thick,-] (2.4,.9) .. controls (3,.4) .. (3.4,.6);
\draw[gray!70, thick] (3.4,.6) ..  controls (3.58,.75) .. (3.6,.9);
\draw[thick,dashed] (3.2,.9) .. controls (3.2,1) and (3.5,1.2) .. (3.7,1.2) ;
\draw[thick,<-](3.2,.9) .. controls (3.2,.6) .. (3.4,.6);
\draw[gray!70,thick,-] (3.6,.9) .. controls  (3.62,1.05) .. (3.7,1.2);
\draw[thick] (3.7,1.2) -- (4.7,1.2);
\draw[thick,-] (2,1.6) .. controls (2.3,1.6) and (3.6,1.4) .. (4,1.4) (4.7,1.2) .. controls (4.6,1.2) and (4.3,1.4) .. (4,1.4) ;
\draw[gray!70,thick] (4.7,1.2) .. controls (5.2,1.2) and (6,1.6) .. (6.6,1.6);
\draw[gray!70,thick] (6.6,1.6) -- (7.4,1.6);
\end{scope}
\begin{scope}[yshift=7.2 cm]
\draw[gray!90] (4.15,1.2) node[left]{\scriptsize $[2,5] \times $};
\draw[gray!70,thick] (1.4,1.6) -- (2,1.6);
\draw[gray!70,thick,->] (2,1.6) .. controls (2.3,1.6) and (3.6,1.4) .. (4,1.4) (4.7,1.2) .. controls (4.6,1.2) and (4.3,1.4) .. (4,1.4);
\draw[gray!90] (4,1.6) node{\scriptsize $\Psi(\ComX(\check{M}))$};
\draw[gray!70,thick] (6.6,1.6) -- (7.4,1.6);
\draw[gray!70,thick] (4.7,1.2) .. controls (5.2,1.2) and (6,1.6) .. (6.6,1.6);
\fill[gray!70] (6.6,1.6) circle (0.065) (2,1.6) circle (0.065);
\end{scope}
\begin{scope}[yshift=4.8 cm]
\draw[gray!90] (3.2,.9) node[left]{\scriptsize $[2,4] \times $};
 \draw[gray!70,thick,<-] (3.4,.6) ..  controls (3.58,.75) .. (3.6,.9);
\draw[gray!70,thick,dashed] (3.2,.9) .. controls (3.2,1) and (3.5,1.2) .. (3.7,1.2);  \draw[gray!70,thick,<-] (3.2,.9) .. controls (3.2,.6) .. (3.4,.6);
\draw[gray!70,thick,<-] (3.6,.9) .. controls  (3.62,1.05) .. (3.7,1.2);
\draw (3.6,.9) node[right]{\scriptsize $\Psi(U\check{M}_{|-L})$};
\end{scope}
\begin{scope}[yshift=6.4 cm,xshift=6.5cm]
\draw (1,.2) -- (6,.2) -- (7.6,1.8) -- (1.6,1.8) -- (1,1.4);
\draw (3.3,.9) node[left]{\scriptsize $\{4\} \times \Psi(S^2 \times \Sigma )$};
\fill[gray!30] (3.7,1.2) .. controls (3.5,1.2) and (3.2,1) .. (3.2,.9) .. controls (3.2,.6) .. (3.4,.6) ..  controls (3.58,.75) .. (3.6,.9) .. controls  (3.62,1.05) .. (3.7,1.2);
\draw[thick] (3.4,.6) ..  controls (3.58,.75) .. (3.6,.9);
\draw[thick] (3.2,.9) .. controls (3.2,1) and (3.5,1.2) .. (3.7,1.2)  (3.2,.9) .. controls (3.2,.6) .. (3.4,.6);
\draw[thick] (3.6,.9) .. controls  (3.62,1.05) .. (3.7,1.2);
\end{scope}
\begin{scope}[yshift=6.4 cm]
\draw (3.25,.9) node[left]{\scriptsize $\{4\} \times \Psi(S^2 \times \Sigma )$};
\fill[gray!40] 
(3.7,1.2) .. controls (3.5,1.2) and (3.2,1) .. (3.2,.9) .. controls (3.2,.6) .. (3.4,.6) ..  controls (3.58,.75) .. (3.6,.9) .. controls  (3.62,1.05) .. (3.7,1.2);
\fill[pattern=north east lines, pattern color=gray!70] 
(3.7,1.2) .. controls (3.5,1.2) and (3.2,1) .. (3.2,.9) .. controls (3.2,.6) .. (3.4,.6) ..  controls (3.58,.75) .. (3.6,.9) .. controls  (3.62,1.05) .. (3.7,1.2);
\draw[thick,>-] (3.4,.6) ..  controls (3.58,.75) .. (3.6,.9);
\draw[thick,>-] (3.2,.9) .. controls (3.2,1) and (3.5,1.2) .. (3.7,1.2)  (3.2,.9) .. controls (3.2,.6) .. (3.4,.6);
\draw[thick,>-] (3.6,.9) .. controls  (3.62,1.05) .. (3.7,1.2);
\end{scope}
\begin{scope}[yshift=8 cm,xshift=6.5cm]
\draw (1,.2) -- (6,.2) -- (7.6,1.8) -- (1.6,1.8) -- (1,1.4);
\end{scope}
\begin{scope}[yshift=3.2 cm,xshift=6.5cm]
\fill[gray!40] (6.6,6.4) .. controls (6,6.4) and (5.2,6) .. (4.7,6) .. controls (4.6,6) and (4.3,6.2) .. (4,6.2) .. controls (3.6,6.2) and (2.3,6.4)   .. (2,6.4) -- (6.6,6.4);
\draw (6.6,6.4) .. controls (6,6.4) and (5.2,6) .. (4.7,6) .. controls (4.6,6) and (4.3,6.2) .. (4,6.2) .. controls (3.6,6.2) and (2.3,6.4)   .. (2,6.4) -- (4.3,6.4);
\draw[-] (6.6,6.4) -- (4.3,6.4);
\draw (4.2,5.8) node[left]{\scriptsize $\{5\} \times \Psi_{[\varepsilon,0]}( \partial \prop_{\ComX})$};
\end{scope}
\end{tikzpicture}
\caption{$\prop_{\ComX} \cap \left([0,6]\times \partial C_2(M)\right)$ and its horizontal pieces}
\label{figschematic}
\end{center}
\end{figure}

\subsection{The perturbating diffeomorphism \texorpdfstring{$\Psi_{\ComY,\varepsilon}$ of $C_2(M)$}{ }}
\label{subpertudif}
Recall that $\ComY$ is a field like in Section~\ref{subdescY}.
For $\eta$ small enough,
we have an isotopy $\psi_{\ComY} \colon [0,\eta] \times \check{M} \rightarrow \check{M}$ such that $\frac{d}{dt} \psi_{\ComY}(t,y)=\ComY(y)$ and $\psi_0$ is the identity.

Let $$\begin{array}{llll}\chi_{\varepsilon}\colon &[0,\varepsilon] &\rightarrow& [0,\varepsilon]\\
&0& \mapsto & \varepsilon\\
&\varepsilon& \mapsto & 0
\end{array}$$ be a smooth family of decreasing functions with horizontal tangents at $0$ and $\varepsilon$ for $\varepsilon \in [0,\eta]$.

Fix $\varepsilon$.
Consider the diffeomorphism $\Psi=\Psi_{\ComY,\varepsilon}$ of $C_2(\check{M})$ that
is the identity outside a neighborhood $U\check{M}\times [0,\varepsilon]$ of the blown-up diagonal, where the second coordinate stands for the distance between two points in a pair, and that reads
$$(v \in U\check{M}_{|m},u) \mapsto (\tang \psi_{\ComY}(\chi_{\varepsilon}(u),m)(0,v),u)$$
on $U\check{M}\times [0,\varepsilon]$, so that it coincides with $ \tang \psi$
on $(U\check{M}=U\check{M} \times \{0\})$, where $\psi=\psi_{\ComY}(\varepsilon,.)$.

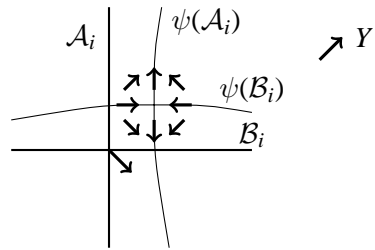
\begin{figure}[h]
\begin{center}
\begin{tikzpicture}
\useasboundingbox (-2,-2) rectangle (1.4,1.4);
\draw [thick] (-1.9,-.6) -- (1.3,-.6) (-.6,-1.9) -- node[very near end, left]{$\Asc_i$} (-.6,1.3);
\draw (.1,1.3) .. controls (0,.7) .. node[very near start, right]{$\psi(\Asc_i)$} (0,0) .. controls (0,-.7) .. (.2,-1.9)
(-1.9,-.2) .. controls (-.7,0) .. (0,0) .. controls (.7,0) .. (1.3,-.1);
\draw [very thick,<-] (0,-.5) -- (0,-.2);
\draw [very thick,<-] (0,.5) -- (0,.2);
\draw [very thick,<-] (-.2,0) -- (-.5,0);
\draw [very thick,<-] (.2,0) -- (.5,0);
\draw [very thick,<-] (.2,.4) -- (.4,.2);
\draw [very thick,<-] (.2,-.4) -- (.4,-.2);
\draw [very thick,<-] (-.2,.4) -- (-.4,.2);
\draw [very thick,<-] (-.2,-.4) -- (-.4,-.2);
\draw [very thick,->] (2.2,.6) -- (2.5,.9);
\draw [very thick,->] (-.6,-.6) -- (-.3,-.9);
\draw (1.3,-.4) node{$\Bdesc_i$} (2.5,.9) node[right]{$\ComY$} (1.3,.2) node{$\psi(\Bdesc_i)$};
\end{tikzpicture}
\caption{Horizontal directions of the flow lines of $\psi_{\ast}(\phi)$}
\label{figprojcrit}
\end{center}
\end{figure}

Define the flow $\psi\phi \psi^{-1}$ ($(t,m) \mapsto \psi\phi_t\psi^{-1}(m)$) on $\check{M}$.
Observe $$\Psi(\overline{s_{\phi}(\check{M})})= s_{\psi\phi \psi^{-1}}(\check{M}).$$
The projections of the directions of the flow lines of $\psi_{\ast}(\phi)=\psi\phi \psi^{-1}$ onto a fiber of the tubular neighborhood of a line $\gamma(c)$ are shown in Figure~\ref{figprojcrit}.
We shall refer to the directions of these projections as {\em horizontal\/} directions.

Without loss, assume that the isotopy $\psi_{\ComY}$ moves the critical points $a_i$ along the lines $\linh(a_i)$ 
and the $b_j$ along the $\linh(b_j)$ (recall that $\ComY$ is tangent to these lines).
Let $\overline{\phi}$ denote the flow $\phi$ reversed so that $\iota(\preprop_{\phi})=\preprop_{\overline{\phi}}$.

\begin{lemma}
\label{lemdirsa}
For $\varepsilon$ small enough, the direction of $\psi_{\ast}(\phi)$ (which is the direction of $s_{\psi_{\ast}(\phi)}$) along $\gamma(c)$ is very close to a geodesic arc between the direction of $\phi$ and $\sigma_a(-\ComY)$, so that its distance in $S^2$ from $\sigma_a(\ComY)$ is at least $\pi/4$.

The direction of $\overline{\phi}$ along $\psi(\gamma(c))$ is very close to a geodesic arc between the direction of $(-T(\psi(\gamma(c))))$
and $\sigma_a(-\ComY)$, so that its distance in $S^2$ from $\sigma_a(\ComY)$ is at least $\pi/4$.

Furthermore, the direction of $\psi_{\ast}(\phi)$ at the critical points and the direction of $\overline{\phi}$  at their images under $\psi$ coincide with $\sigma_a(-\ComY)$.
\end{lemma}
\bp Away from the ends of $\gamma(c)$, the direction of $\psi_{\ast}(\phi)$ along $\gamma(c)$ is very close to the tangent direction of $\gamma(c)$,  and it is slightly deviated in the orthogonal direction of $\sigma_a(-\ComY)$ since $\gamma(c)$ is obtained from $\psi(\gamma(c))$ by a translation of 
$-\ComY$. See Figure~\ref{figprojcrit} and Subsection~\ref{subframedsel}. Near the critical points, the direction of $\psi_{\ast}(\phi)$ approaches the direction
of $\sigma_a(-\ComY)$, and it reaches it at the critical points.
Similarly, the direction of $\overline{\phi}$ along $\psi(\gamma(c))$  is very close to the direction of $(-T(\gamma(c)))$ away from the ends and it is slightly deviated in the orthogonal direction of $(-\sigma_a(\ComY))$. Near the critical points, the direction of $\overline{\phi}$ approaches the direction
of $\sigma_a(-\ComY)$, and it reaches it at the critical points.
\eop

\begin{lemma}
\label{lemfifi}
$\lim_{\varepsilon \to 0} \Psi(\preprop_{\phi}) \cap \iota (\preprop_{\phi})$ is discrete located at the points $s_{\sigma_a(-\ComY)}(a_i)$ and $s_{\sigma_a(-\ComY)}(b_j)$ of $U\check{M}$, which are the unit tangent vectors directed by $\sigma_a(-\ComY)$ at the critical points.
\end{lemma}
\bp
Observe that $\preprop_{\phi} \cap \iota (\preprop_{\phi})$ is supported on the restrictions of $U\check{M}$ to the critical points. Therefore, for $\varepsilon$ small enough, $\Psi(\preprop_{\phi}) \cap \iota (\preprop_{\phi})$ will be near the restrictions of $U\check{M}$ to the critical points. There are 4g points of type $s_{\overline{\phi}}(\psi(a_i))$, $s_{\psi_{\ast}(\phi)}(a_i)$, $s_{\overline{\phi}}(\psi(b_j))$ and $s_{\psi_{\ast}(\phi)}(b_j)$ in the intersection. They have the wanted direction thanks to Lemma~\ref{lemdirsa}. Except for those points we have to look for flow lines for $\phi$ and flow lines for $\psi_{\ast}(\phi)$ that intersect twice and that connect the intersection points with opposite directions.
Under our assumptions, this can only happen on the lines $\linh(c)$ between $c$ and $\psi(c)$ for a critical point $c$. Indeed, outside $\linh(c)$, $\phi$ and $\psi_{\ast}(\phi)$ both escape from the neighborhoods of $\linh(c)$ if $c=a_i$, or both get closer if $c=b_i$. On these lines, the only parts where $\phi$ and $\psi_{\ast}(\phi)$ have opposite directions is between $c$ and $\psi(c)$, and the tangent direction to $\overline{\phi}$ is the direction of $\sigma_a(-\ComY)$.
\eop

\subsection{Reduction of the proof of Proposition~\ref{propcycintdouble}}

Consider a regular neighborhood $N$ of the union of the $\gamma(c)$ that contains the $\psi(\gamma(c))$, and consider the fiber
bundle over $N$ whose fibers are the complement of an open disk of radius $\pi/4$ around $\sigma_a(\ComY)$ in the fibers of $UN$. Let $E$ be the total space of this bundle and let $\CN = [-1,0] \times E \subset [-1,0] \times \partial C_2(M) \subset C_2(M)$. Then $H_2(\CN;\ZZ)=0$.

Without loss, the chains $\prop_{\ComX}$ and $\prop_{-\ComX}$ are now assumed to be transverse so that their intersection $I$ is a $2$--cycle of $C_2(M)$, which we are going to compute piecewise. We shall neglect the pieces in $\cal N$ and write them as $O(\CN)$ in the statements. Sometimes, we shall also add arbitrary pieces in $\cal N$ in order to close some  $2$--chains and find some $2$--cycle $I^{\prime}$ such that $$I^{\prime} = I + O(\CN)$$
so that $I^{\prime}$ will be homologous to $I$.

We shall also consider continuous limits when possible to simplify the expressions as in Lemma~\ref{lemfifi}, which now reads:
$$\lim_{\varepsilon \to 0} \Psi(\preprop_{\phi}) \cap \iota (\preprop_{\phi}) = O(\CN)$$ or, \\
for $\varepsilon > 0$ small enough,
$\Psi(\preprop_{\phi}) \cap \iota (\preprop_{\phi}) = O(\CN)$.

For example,
$$\begin{array}{lll} \prop_{\ComX} \cap \prop_{-\ComX} \cap ([5/2,6]\times \partial C_2(M)) &= &[5/2,3] \times (\psi_{\ast}(\ComX)(\Link)-(-\ComX)(\psi(\Link))) \\
  & & +\{3\} \times (-\psi_{\ast}(\ComX)(\Sigma) + S^2 \times(\psi(\Link) \cap \Sigma) )\\
&&-[3,4] \times (-\ComX)(\psi(\Link))\\
&& + \{4\} \times (-\ComX)(\psi(\Sigma))\\
&=& \{3\} \times (-\psi_{\ast}(\ComX))(\Sigma) + \{4\} \times (-\ComX)(\psi(\Sigma))\\
&&+  \{3\} \times S^2 \times (\psi(\Link) \cap \Sigma) + O(\CN).
\end{array}$$
Then $S^2 \times (\psi(\Link) \cap \Sigma)$ is a disjoint union of spheres homologous to $lk(\Link,\Link_{\parallel \ComY}) [S]$.
Let $$\ell=\lim_{\varepsilon \to 0} \left(-\{3\} \times (\psi_{\ast}(\ComX))(\Sigma) + \{4\} \times (-\ComX)(\psi(\Sigma))\right).$$
$$\begin{array}{lll}\ell&=&-\{3\} \times \ComX(\Sigma) + \{4\} \times (-\ComX)(\Sigma)\\
&=&-\{3\} \times \ComX(\Sigma) + \{4\} \times (-\ComX)(\Sigma) \\&& -[3,4] \times (-\ComX)(\Link)
+\{3\} \times s_{[-\ComX,\ComX]_{\sigma_a(-\ComY)}}(\Link) + O(\CN)
  \end{array}
$$
where the last equality comes from the fact that both $[3,4] \times (-\ComX)(\Link)$ and 
$\{3\} \times s_{[-\ComX,\ComX]_{\sigma_a(-\ComY)}}(\Link)$ are in $\CN$.
Then $\prop_{\ComX} \cap \prop_{-\ComX} \cap ([5/2,6]\times \partial C_2(M))$ is homologous to 
$\Gdu^b(\ComX,\ComY)$ mod $\CN$
and the proof of Proposition~\ref{propcycintdouble} is reduced to the proof of the two following propositions.

\begin{proposition}
\label{propcycintdoubleint}
$$\prop_{\ComX} \cap \prop_{-\ComX} \cap \tilde{C}_2(M)=\Gdu^i(\ComY) + O(\CN).$$
\end{proposition}

\begin{proposition}
\label{propcycintdoublebord}
$$\prop_{\ComX} \cap \prop_{-\ComX} \cap \left([0,5/2]\times \partial C_2(M)\right)= O(\CN).$$
\end{proposition}

In particular, $\Gdu^i(\ComY)$ may be thought of as the intersection of $\prop(\funcf,\metrigo) \cap \prop(-\funcf,\metrigo)$ in the interior of $C_2(M)$, while $\Gdu^b(\ComY)$ collects the intersection coming from the boundary corrections.

\subsection{Proof of Proposition~\ref{propcycintdoubleint}}

\begin{lemma}
\label{lemii}
$$\lim_{\varepsilon \to 0}\Psi(\preprop_{\CI}) \cap \iota (\preprop_{\CI})=\sum_{(i,j,k,\ell) \in \{1,\dots,g\}^4}\CJ_{ji}\CJ_{\ell k}(\Bdesc_j\cap \Asc_{k}) \times (\Bdesc_{\ell}\cap \Asc_i)_{\parallel \sigma_a(-\ComY)}+O(\CN).
$$
\end{lemma}
\bp
The intersection
$\Bdesc_j \times \Asc_i \cap (\Asc_{k}\times \Bdesc_{\ell})$ is cooriented by the positive normals of $\Bdesc_j$, $\Asc_i$, $\Asc_k$ and $\Bdesc_{\ell}$ in this order. Therefore the intersection reads as in the statement of the lemma away from the diagonal.
Near the diagonal and away from the critical points, $\Asc_i$ and $\Bdesc_j$ are moved in the direction of $\ComY$. If $\ComY =\vec{a} + \vec{b}$ where $\vec{a}$ is tangent to $\Asc_{k}$ and $\vec{b}$ is tangent to $\Bdesc_{\ell}$, then abusively write $\Asc_i= \Asc_{k} + \vec{b}$ and $\Bdesc_j= \Bdesc_{\ell} + \vec{a}$ and see that the difference of the two points is moved in the direction $(\vec{b}-\vec{a})$ of $\sigma_a(-\ComY)$, so that the corresponding intersection sits inside the neglected part $\CN$. 
(When two points vary along the same $\gamma(c)$, the second one will be deviated in the direction of $\sigma_a(-\ComY)$ so that the limit pairs of points describe an arc in $U\check{M}_{|\gamma(c)}$ from $-\gamma^{\prime}(c)$ to $\gamma^{\prime}(c)$ through $\sigma_a(-\ComY)$, that is along the half great circle $[-\gamma^{\prime}(c),\gamma^{\prime}(c)]_{\sigma_a(-\ComY)}$.)

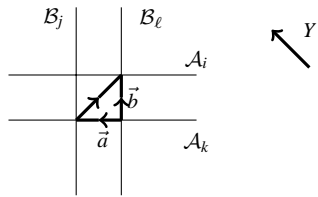
\begin{figure}[h]
\begin{center}
\begin{tikzpicture}
\useasboundingbox (-1.6,-1.1) rectangle (1.1,1.6);
\draw (-1.5,0) -- (1,0) (0,-1) -- (0,1.5);
\draw (-1.5,.6) -- (1,.6) (-.6,-1) -- (-.6,1.5);
\draw [very thick,->] (0,.3) -- (0,.6) (0,0) -- (0,.3);
\draw [very thick,->] (-.3,0) -- (-.6,0) (0,0) -- (-.3,0);
\draw [very thick,->] (2.5,.7) -- (2,1.2);
\draw [very thick,->] (-.3,.3) -- (0,.6) (-.6,0) -- (-.3,.3);
\draw (-.25,-.25) node{\scriptsize $\vec{a}$} (2.5,1.2) node{\scriptsize $\ComY$} (.15,.3) node{\scriptsize $\vec{b}$} (.4,1.35) node{\scriptsize $\Bdesc_{\ell}$} (-.9,1.35) node{\scriptsize $\Bdesc_j$} (1,.8) node{\scriptsize $\Asc_i$} (1,-.3) node{\scriptsize $\Asc_{k}$};
\end{tikzpicture}
\caption{Deviation near the diagonal}
\label{figdevdiag}
\end{center}
\end{figure}
Near a critical point, two points can come from different crossings.
Then the direction between them in $\overline{(\Bdesc_j\cap \Asc_{k}) \times (\Bdesc_{\ell}\cap \Asc_i) \setminus \mbox{diag}}$ is orthogonal to $\ComY = \pm \sigma_a(\ComY)$.
The field $\ComY$ can be assumed to preserve the $\Bdesc$-sheets near the $a_i$ and the $\Asc$-sheets near the $b_j$.
Then the difference of the two points is moved in the direction of $\sigma_a(-\ComY)$ so that it belongs to the hemisphere centered at $\sigma_a(-\ComY)$.
\eop

\begin{lemma}
\label{lemfii}
$$ \lim_{\varepsilon \to 0}\Psi(\preprop_{\phi}) \cap \iota (\preprop_{\CI})=\sum_{c \in \CaC} \CJ_{j(c)i(c)}(-\sigma(c))\overline{ \{(\gamma(c)(t_1),\gamma(c)(t_2));t_1 < t_2\}} + O(\CN).$$
\end{lemma}
\bp
The intersection $\preprop_{\phi} \cap \left(\iota (\preprop_{\CI})=\sum_{(i,j) \in \{1,\dots,g\}^2} \CJ_{ji} \Asc_i \times \Bdesc_j\right)$ is supported on the $$\overline{ \{(\gamma(c)(t_1),\gamma(c)(t_2));t_1 < t_2\}}$$ away from the unit bundles of the critical points.
It is transverse except near these unit bundles.

Let $c\in \alpha_i \cap \beta_j$.
Along $\gamma(c)$, $\Asc_i \times \Bdesc_j$ is cooriented by $\beta_j \times \alpha_i$. Then
$\preprop_{\phi} \cap (\Asc_i \times \Bdesc_j)$ will be oriented as $(-\sigma(c)) \{(\gamma(c)(t_1),\gamma(c)(t_2));t_1 < t_2\}$.
Since $\psi_{\ast}(\phi)$ is almost vertical away from the critical points, we are left with the behaviour near the critical points. Near $a_i$ on $\Asc_i$, (or near $b_j$ on $\Bdesc_j$) the direction of $\psi_{\ast}(\phi)$ is in the hemisphere centered at $\sigma_a(-\ComY)$,
according to Lemma~\ref{lemdirsa}, so that the pairs of points of $\Asc_i \times \Bdesc_j$ connected by flow lines of $\psi_{\ast}(\phi)$ near a critical point are in $\CN$.
\eop

Similarly, we have
\begin{lemma}
\label{lemifi}
$$ \lim_{\varepsilon \to 0}\Psi(\preprop_{\CI}) \cap \iota (\preprop_{\phi})=\sum_{c \in \CaC} \CJ_{j(c)i(c)}(-\sigma(c))\overline{ \{(\gamma(c)(t_1),\gamma(c)(t_2));t_1 > t_2\}} + O(\CN).$$
\end{lemma}
\bp
Away from the unit bundles of the critical points, it is clear.
According to Lemma~\ref{lemdirsa}, the direction of $\overline{\phi}$ on $\psi(\Asc_i)$ near $\psi(a_i)$ (or on $\psi(\Bdesc_j)$ near $\psi(b_j)$) is in the hemisphere centered at $\sigma_a(-\ComY)$, so that the pairs of points of $(\psi(\Bdesc_j) \times \psi(\Asc_i)) \cap \iota (\preprop_{\phi})$ near  the critical points are again in $\CN$.
\eop

Proposition~\ref{propcycintdoubleint} is a direct corollary of Lemmas~\ref{lemfifi}, \ref{lemii}, \ref{lemfii}, \ref{lemifi}. \eop

\subsection{Proof of Proposition~\ref{propcycintdoublebord}}

We prove that
$\left(\prop_{\ComX} \cap \prop_{-\ComX} \cap ([0,5/2]\times \partial C_2(M))\right)$ is in $\CN$.

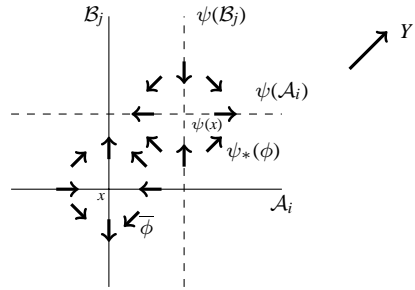
\begin{figure}
\begin{center}
\begin{tikzpicture}
\useasboundingbox (-2.4,-2.4) rectangle (1.4,1.4);
\draw [dashed] (-2.3,0) -- (1.3,0) (0,-2.3) --  (0,1.3);
\draw (0,1.3) node[right]{\scriptsize $\psi(\Bdesc_j)$} (-.9,1.3) node[left]{\scriptsize $\Bdesc_j$};
\draw [very thick,->] (0,-.7) -- (0,-.4);
\draw [very thick,->] (0,.7) -- (0,.4);
\draw [very thick,->] (-.4,0) -- (-.7,0);
\draw [very thick,->] (.4,0) -- (.7,0);
\draw [very thick,->] (.3,.5) -- (.5,.3);
\draw [very thick,->] (.3,-.5) -- (.5,-.3);
\draw [very thick,->] (-.3,.5) -- (-.5,.3);
\draw [very thick,->] (-.3,-.5) -- (-.5,-.3);
\draw [very thick,->] (2.2,.6) -- (2.7,1.1);
\draw (1.3,.3) node{\scriptsize $\psi(\Asc_i)$} (1.3,-1.2) node{\scriptsize $\Asc_i$} (2.7,1.1) node[right]{\scriptsize $\ComY$} (.3,-.15) node{\tiny $\psi(x)$} (-.5,-1.5) node{\scriptsize $\overline{\phi}$} (-1.1,-1.1) node{\tiny $x$} (.4,-.5) node[right]{\scriptsize $\psi_{\ast}(\phi)$};
\draw (-2.3,-1) -- (1.3,-1) (-1,-2.3) --  (-1,1.3);
\draw [very thick,<-] (-1,-1.7) -- (-1,-1.4);
\draw [very thick,<-] (-1,-.3) -- (-1,-.6);
\draw [very thick,<-] (-1.4,-1) -- (-1.7,-1);
\draw [very thick,<-] (-.6,-1) -- (-.3,-1);
\draw [very thick,<-] (-.7,-.5) -- (-.5,-.7);
\draw [very thick,<-] (-.8,-1.5) -- (-.6,-1.3);
\draw [very thick,<-] (-1.3,-.5) -- (-1.5,-.7);
\draw [very thick,<-] (-1.3,-1.4) -- (-1.5,-1.2);
\end{tikzpicture}
\caption{Tangencies of the flow lines of $\overline{\phi}$ and $\psi_{\ast}(\phi)$ near some $\gamma(c)$}
\label{figpsiphi}
\end{center}
\end{figure}

According to Theorem~\ref{thmbordFfii}, $$ \partial \prop(\funcf,\metrigo) = \partial_{od} + \sum_{c \in \CaC} \CJ_{j(c)i(c)} \sigma(c) (S^2 \times \gamma(c)) +\overline{s_{\phi}(\check{M})}.$$
Therefore, according to Lemmas~\ref{lemfifi} and \ref{lemdirsa}, 
$$\Psi(\partial \prop(\funcf,\metrigo)) \cap \partial \iota(\prop(\funcf,\metrigo))=O(\CN).$$
Let us now show that
$$\Psi(\partial \prop(\funcf,\metrigo)) \cap \iota(\preprop_{\homotop})=O(\CN).$$
According to the construction of $\preprop_{\homotop}$ in Subsections~\ref{subframedsel} and \ref{subFcomb}, $\iota(\preprop_{\homotop})$ intersects $\Psi(S^2 \times \gamma(c))=S^2 \times \psi(\gamma(c))$ 
on $s_{[\overline{\phi},-\ComX]}(\psi(\gamma(c))$ where $[\overline{\phi},-\ComX]$ is the shortest geodesic arc between the tangent to $\overline{\phi}$ and $-\ComX$, which is in the hemisphere centered at $\sigma_a(-\ComY)$, according to Lemma~\ref{lemdirsa}. Now, look at the intersection of
$\iota(\preprop_{\homotop})$ and $s_{\psi_{\ast}(\phi)}(\check{M})$, where the direction of $\psi_{\ast}(\phi)$ must belong to $[\overline{\phi},-\ComX]$.
This can only happen in a tubular neighborhood of $\gamma_i$ at a place where
the flow lines of $\psi_{\ast}(\phi)$ and $\overline{\phi}$ have the same horizontal direction. This only happens between $\gamma_i$ and $ \psi(\gamma_i)$, more precisely in the preimage of the rectangle shown in Figure~\ref{figpsiphi} under the orthogonal projection directed by $\ComX_1$.
There the horizontal direction is close to the direction of $\sigma_a(-\ComY)$.

Similarly,
$$\Psi(\preprop_{\homotop}) \cap \left(S^2 \times \Link + \overline{(-\ComX)(\check{M})}\right) = O(\CN).$$
Indeed, since the horizontal component of the direction of $\psi_{\ast}(\phi)$ along $\gamma(c)$ is in the direction of $\sigma_a(-\ComY)$,
$\Psi(\preprop_{\homotop}) \cap (S^2 \times \Link) = O(\CN).$
Now, $(-\ComX)$ can belong to $[\psi_{\ast}(\phi),\psi_{\ast}(\ComX)]$ only if the direction of the horizontal component of $(-\ComX)$, which is the direction of the horizontal component of $s_{\overline{\phi}}$, is the same as the direction of the horizontal component of $s_{\psi_{\ast}(\phi)}$.
This can only happen in the same rectangles as before where $(-\ComX)$ is in the hemisphere centered at $\sigma_a(-\ComY)$.
Finally, $$\Psi(-S^2 \times \Link + \overline{\ComX(\check{M})}) \cap \left(S^2 \times \Link + \overline{(-\ComX)(\check{M})}\right) = O(\CN)$$ since it only consists of unit tangent vectors to $\check{M}$ over $\Link \cup \psi(\Link)$ in the direction of $\pm \ComX$.
\eop

\section{Concluding the proof of Theorem~\ref{thmmain}}
\setcounter{equation}{0}
\label{secproof}

Recall that $\pointw$, $\matchingo$, $\ComX=\ComX(\pointw,\matchingo)$, $\Link=\Link(\matchingo)$ and $\Sigma$ such that $\partial \Sigma= \Link$ are fixed. 
Note that $\ComX$ depends neither on the orientations of the $\alpha_i$ and the $\beta_j$, nor on their order.
Furthermore $e(\pointw,\matchingo)$ is independent of the order of the $\beta_j$.
Thus, the permutation $\rho$ of $\{1,2,\dots,g\}$ associated with $\matchingo$ is assumed to be the identity, without loss.

\subsection{Reducing the proof of Theorem~\ref{thmmain} to an Euler class computation}
\label{subYepseta}

Define four nowhere zero fields $\ComY^{++}$, $\ComY^{+-}$, $(\ComY^{-+}=-\ComY^{+-})$ and $(\ComY^{--}=-\ComY^{++})$ of $\ComX^{\perp}$, up to homotopy among nowhere zero fields, over a neighborhood of the $\gamma(c)$, so that
\begin{itemize}
\item $\ComY^{++}$ and $\ComY^{+-}$ are \emph{positive normals} for $\Asc_i$ on 
$H_{a,\leq 3}=C_M \cap \fMorse^{-1}(]-\infty,3])$ --meaning that $\ComY^{++}$ (or $\ComY^{+-}$) followed by an oriented basis of the tangent space to $\Asc_i$ gives rise to an oriented basis of the tangent space to $M$--, and
\item $\ComY^{++}$ and $\ComY^{-+}$ are positive normals for $\Bdesc_j$ on $H_{b,\geq 3}=C_M \cap \fMorse^{-1}([3,+\infty[)$,
\end{itemize}
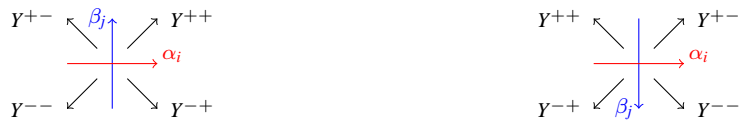
\begin{figure}[h]
\begin{center}
\begin{tikzpicture}
\useasboundingbox (-4.2,-.8) rectangle (4.2,.8);
\begin{scope}[xshift=-3.5cm]
\draw [blue,->] (-.2,.6) node{\scriptsize $\beta_j$} (0,-.6) -- (0,.6);
\draw [red,->] (.5,.1) node[right]{\scriptsize $\alpha_i$} (-.6,0) -- (.6,0);
\draw [->] (.2,.2) -- (.6,.6) node[right]{\scriptsize $\ComY^{++}$};
\draw [->] (-.2,.2) -- (-.6,.6) node[left]{\scriptsize $\ComY^{+-}$};
\draw [->] (-.2,-.2) -- (-.6,-.6) node[left]{\scriptsize $\ComY^{--}$};
\draw [->] (.2,-.2) -- (.6,-.6) node[right]{\scriptsize $\ComY^{-+}$};
\end{scope}
\begin{scope}[xshift=3.5cm]
\draw [blue,->] (-.2,-.6) node{\scriptsize $\beta_j$} (0,.6) -- (0,-.6);
\draw [red,->] (.5,.1) node[right]{\scriptsize $\alpha_i$} (-.6,0) -- (.6,0);
\draw [->] (.2,.2) -- (.6,.6) node[right]{\scriptsize $\ComY^{+-}$};
\draw [->] (-.2,.2) -- (-.6,.6) node[left]{\scriptsize $\ComY^{++}$};
\draw [->] (-.2,-.2) -- (-.6,-.6) node[left]{\scriptsize $\ComY^{-+}$};
\draw [->] (.2,-.2) -- (.6,-.6) node[right]{\scriptsize $\ComY^{--}$};
\end{scope}
\end{tikzpicture}
\caption{The fields $\ComY^{\varepsilon,\eta}$ on $\fMorse^{-1}(\{3\})$}
\label{figypm}
\end{center}
\end{figure}
More explicitly, such fields $\ComY^{\varepsilon,\eta}$ can be pictured on $\fMorse^{-1}(\{3\})$ as in Figure~\ref{figypm} and the field $\ComY^{++}$ becomes closer to the ``actual positive orthogonal normal'' to $\Asc_i$ as we approach $a_i$, and closer to the ``actual positive orthogonal normal'' to $\Bdesc_j$, as we approach $b_j$, where Figure~\ref{figypm} cannot be drawn anymore. (In order to determine these fields up to homotopy among nowhere zero fields, it is enough to determine an open half-space where they lie, continuously.)

Then with the notation of Subsections~\ref{subparflow} and \ref{subframedHd},
$$lk(\Link(\matchingo),\Link(\matchingo)_{\parallel})=\frac14 \sum_{(\varepsilon,\eta) \in \{+,-\}^2}lk (\Link,\Link_{\parallel \ComY^{\varepsilon,\eta}})$$
and, with the notation of Proposition~\ref{propcycintdouble},  $$[\Gdu]=\frac14 \sum_{(\varepsilon,\eta) \in \{+,-\}^2}\left[\Gdu^i(\ComY^{\varepsilon,\eta}) + \Gdu^b(\ComX,\ComY^{\varepsilon,\eta})\right]$$
where $\sigma_a(-\ComY^{\varepsilon,\eta})=\ComY^{\varepsilon, (-\eta)}$,
so that the collection of the $\sigma_a(-\ComY^{\varepsilon,\eta})$ is the same as the collection of the $\ComY^{\varepsilon,\eta}$ and, thanks to Lemma~\ref{lemnorpos}, 
$$[\twocycG(\CD)]=\frac14 \sum_{(\varepsilon,\eta) \in \{+,-\}^2}[\Gdu^i(\ComY^{\varepsilon,\eta})]$$
with the notation of Proposition~\ref{propChMorse}.

Therefore, thanks to Proposition~\ref{propcycintdouble}, the proof of Theorem~\ref{thmmain} is reduced to the proof of the following equality in $H_2(U\check{M};\QQ)$.
$$\left[\ComX(\Sigma) - (-\ComX)(\Sigma) -\frac14 \sum_{(\varepsilon,\eta) \in \{+,-\}^2}s_{[-\ComX,\ComX]_{\ComY^{\varepsilon,\eta}}}(\partial \Sigma)\right]=e(\pointw,\matchingo)[S].$$

Consider the rank $2$ sub-vector bundle $\ComX^{\perp}$ of $T\check{M}$ of the planes orthogonal to $\ComX$. Let $\ComX^{\perp}(\Sigma)$ be the total space of the restriction of $\ComX^{\perp}$ to our surface $\Sigma$. 
Let $\ComY$ be a  nowhere zero section of $\ComX^{\perp}$ on $\partial \Sigma$. The {\em relative Euler class $e(\ComX^{\perp}(\Sigma),\ComY)$
of $\ComY$ in $\ComX^{\perp}(\Sigma)$} is the obstruction to extending $\ComY$ as a nonzero section of $\ComX^{\perp}(\Sigma)$ over $\Sigma$. If $\tilde{\ComY}$ is an extension of $\ComY$ as a section of $\ComX^{\perp}(\Sigma)$ transverse to the zero section $s_0(\ComX^{\perp}(\Sigma))$,
then $$e(\ComX^{\perp}(\Sigma),\ComY)=\langle \tilde{\ComY}(\Sigma), s_0(\ComX^{\perp}(\Sigma))\rangle_{\ComX^{\perp}(\Sigma)}.$$
\begin{lemma}
Under the assumptions above,
$$\left[\ComX(\Sigma) - (-\ComX)(\Sigma) - s_{[-\ComX,\ComX]_{\ComY}}(\partial \Sigma)\right]=e(\ComX^{\perp}(\Sigma),\ComY)[S]$$ in $H_2(C_2(M))$.
\end{lemma}
\bp If $\ComY$ extends as a nonzero section of $\ComX^{\perp}(\Sigma)$ still denoted by $\ComY$, then the cycle of the left-hand side bounds $s_{[-\ComX,\ComX]_{\ComY}}(\Sigma)$. This allows us to reduce the proof to the case when $\Sigma$ is a neighborhood of a zero of the extension $\tilde{\ComY}$ above, that is when $\Sigma$ is a disk $\Delta$ equipped with a trivial $D^2$-bundle, and when $\ComY \colon \partial \Delta \rightarrow \partial D^2$ has degree $d=\pm 1$. Then
$d=e(\ComX^{\perp}(\Delta),\ComY)$, and $\left[\ComX(\Delta) - (-\ComX)(\Delta) - s_{[-\ComX,\ComX]_{\ComY}}(\partial \Delta)\right]=d[S]$.
\eop

Thus,
$$\left[\ComX(\Sigma) - (-\ComX)(\Sigma) -\frac14 \sum_{(\varepsilon,\eta) \in \{+,-\}^2}s_{[-\ComX,\ComX]_{\ComY^{\varepsilon,\eta}}}(\partial \Sigma)\right]$$ $$=\frac14 \sum_{(\varepsilon,\eta) \in \{+,-\}^2}e(\ComX^{\perp}(\Sigma),\ComY^{\varepsilon,\eta})[S].$$

The proof of Theorem~\ref{thmmain} is now reduced to the proof of the following proposition, which occupies the end of this section.

\begin{proposition}
\label{propcontrob}
$$e(\pointw,\matchingo)=\frac14 \sum_{(\varepsilon,\eta) \in \{+,-\}^2}e(\ComX(\pointw,\matchingo)^{\perp}(\Sigma),\ComY^{\varepsilon,\eta}).$$
\end{proposition}

\begin{remark}
\label{rkproofind}
Note that this proposition provides a combinatorial formula for the average of the Euler classes in the right-hand side. In this formula, the $\tilth(\beta_j)$ and $\tilth(\halfbra c_{j(c)},c \halfbra_{\beta})$ depend on our rectangular diagram of $(\CD,\matchingo,\pointw)$ in Figure~\ref{figplansplit}. Thus, the proposition implies that the sum $e(\pointw,\matchingo)$ is independent of our special picture of the Heegaard diagram.
\end{remark}

\subsection{A surface \texorpdfstring{$\Sigma(\Link(\matchingo))$}{ associated with L}}
\label{subsurL}

Let $H_{b,\geq 2}=C_M \cap \fMorse^{-1}([2,+\infty[)$. For any crossing $c$ of $\CaC$, define the triangle $T_{\beta}(c)$ in the disk $(D_{\geq 2}(\beta_{j(c)})=\Bdesc_{j(c)}\cap H_{b,\geq 2})$ such that
$$ \partial T_{\beta}(c) =[ c_{j(c)},c ]_{\beta} + (\gamma(c) \cap H_{b,\geq 2})- (\gamma_{j(c)} \cap H_{b,\geq 2}).$$
Similarly, define the triangle $T_{\alpha}(c)$ in the disk $(D_{\leq 2}(\alpha_{i(c)})=\Asc_{i(c)}\cap \handleboda)$ such that
$$ \partial T_{\alpha}(c) =-[ c_{i(c)},c ]_{\alpha} + (\gamma(c) \cap \handleboda)- (\gamma_{i(c)} \cap \handleboda).$$
\begin{proposition}
\label{propSigL} Recall $H_{a,2}=C_M \cap \fMorse^{-1}(2)$.
There exists a $2$-chain $F(\matchingo)$ in $H_{a,2}$ such that the boundary of
$$\Sigma(\Link(\matchingo))=F(\matchingo)-\sum_{c \in \CaC} \CJ_{j(c)i(c)}\sigma(c)(T_{\beta}(c) + T_{\alpha}(c)) $$
$$+\sum_{(j,i) \in \{1,\dots,g\}^2}\sum_{c \in \CaC}\CJ_{j(c)i(c)}\sigma(c)\CJ_{ji} \left(\langle \alpha_i, \halfbra c_{j(c)},c \halfbra_{\beta} \rangle D_{\geq 2}(\beta_j)
-\langle \halfbra c_{i(c)},c \halfbra_{\alpha},\beta_j\rangle D_{\leq 2}(\alpha_i)\right)$$
is $\Link(\matchingo)$.
\end{proposition}
\bp
The boundary of the defined pieces reads $(\Link(\matchingo) + u)$
where the cycle $u$ is
$$\begin{array}{ll}u=&\sum_{c \in \CaC} \CJ_{j(c)i(c)}\sigma(c)\left([ c_{i(c)},c ]_{\alpha} - [ c_{j(c)},c ]_{\beta}\right)\\
  &+\sum_{(j,i) \in \{1,\dots,g\}^2}\sum_{c \in \CaC}\CJ_{j(c)i(c)}\sigma(c)\CJ_{ji} \left(\langle \alpha_i, \halfbra c_{j(c)},c \halfbra_{\beta} \rangle \beta_j
-\langle \halfbra c_{i(c)},c \halfbra_{\alpha},\beta_j\rangle \alpha_i\right).
\end{array}$$
Compute $\langle\alpha_k, u\rangle$, by pushing $u$ in the direction of the positive normal to $\alpha_k$ and in the direction of the negative normal, and by averaging.
Since $\alpha_k$ intersects neither the pushed $[ c_{i(c)},c ]_{\alpha}$ nor the pushed $\alpha_i$, and since its intersection with the above average of the pushed $[ c_{j(c)},c ]_{\beta}$ is $\langle \alpha_k,  \halfbra c_{j(c)},c \halfbra_{\beta}\rangle$,
$$\langle\alpha_k, u\rangle=-\sum_{c \in \CaC} \CJ_{j(c)i(c)}\sigma(c) \langle\alpha_k, \halfbra c_{j(c)},c \halfbra_{\beta}\rangle+\sum_{c \in \CaC}\CJ_{j(c)i(c)}\sigma(c) \langle \alpha_k, \halfbra c_{j(c)},c \halfbra_{\beta} \rangle = 0.
$$
Similarly, $\langle u, \beta_{\ell}\rangle=0$ for any $\ell$ so that $(-u)$ bounds a $2$-chain $F(\matchingo)$ in $H_{a,2}$. 
\eop

\subsection{Proof of the combinatorial formula for the Euler classes}
\label{subobs}

In this section, we prove Proposition~\ref{propcontrob}.

\begin{figure}[h]
\begin{center}
\begin{tikzpicture}
\useasboundingbox (1,1.4) rectangle (5,3.3);
\draw (1.3,1.5) .. controls (2,1.5) and (2,3.2) .. (3,3.2);
\draw [->] (4.7,1.5) .. controls (4,1.5) and (4,3.2) .. (3,3.2) ;
\draw (3,2) circle (.5);
\draw [blue,thick,fill=gray!20] (2.64,1.65) arc (225:-45:.5);
\draw [blue,thick,->] (2.5,2) arc (180:0:.5);
\draw [blue,thick] (1.8,1.5) .. controls (1.85,1.6) .. (2,2);
\draw [blue,thick] (4.2,1.5) .. controls (4.15,1.6) .. (4,2);
\draw [blue,thick,->] (2,2) .. controls (2.43,3.45) and (3.57,3.45) .. (4,2);
\draw [blue,thick] (1.9,1.45) .. controls (1.95,1.6) .. (2.1,2);
\draw [blue,thick] (4.1,1.45) .. controls (4.05,1.6) .. (3.9,2);
\draw [blue,thick,->] (2.1,2) .. controls (2.5,3.3) and (3.5,3.3) .. (3.9,2);
\draw [blue,thick] (2.05,1.4) .. controls (2.1,1.6) .. (2.25,2);
\draw [blue,thick] (3.95,1.4) .. controls (3.9,1.6) .. (3.75,2);
\draw [blue,thick,<-] (2.25,2) .. controls (2.6,3.1) and (3.4,3.1) .. (3.75,2);
\draw [blue] (3.5,2) node[left]{\scriptsize $\beta_1$};
\draw [red,thick,->] (3,3.2) .. controls (2.9,3.2) and (2.8,3) .. (2.8,2.85);
\draw [red,thick] (3,2.5) .. controls (2.9,2.5) and (2.8,2.7) .. (2.8,2.85) ;
\draw [red,thick,dashed] (3,3.2) .. controls (3.1,3.2) and (3.2,3) .. (3.2,2.85) .. controls (3.2,2.7) and (3.1,2.5) .. (3,2.5);
\draw [red,-<] (2.6,1.7) .. controls (2.4,1.6) and (2.3,1.55) ..
(2.15,1.55);
\draw [red] (2.15,1.55).. controls (2,1.55) and (1.9,1.6) .. (1.7,1.7);
\draw [red,dashed] (2.6,1.7) .. controls (2.4,1.8) and (2.3,1.85) ..
(2.15,1.85) .. controls (2,1.85) and (1.9,1.8) .. (1.7,1.7);
\draw [red,-<] (3.4,1.7) .. controls (3.6,1.6) and (3.7,1.55) ..
(3.85,1.55);
\draw [red] (3.8,1.4) node[left]{\scriptsize $\alpha^{\prime}_1$};
\draw [red] (2.2,1.4) node[right]{\scriptsize $\alpha^{\prime\prime}_1$};
\draw [red] (3.85,1.55).. controls (4,1.55) and (4.1,1.6) .. (4.3,1.7);
\draw [red,dashed] (3.4,1.7) .. controls (3.6,1.8) and (3.7,1.85) ..
(3.85,1.85) .. controls (4,1.85) and (4.1,1.8) .. (4.3,1.7);
\end{tikzpicture}
\caption{How the $\beta_j$ look like near the handles' cores}
\label{figinthandles}
\end{center}
\end{figure}
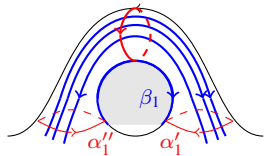

Represent $\handleboda$ like in Figure~\ref{figHa}, and assume that the curves $\beta_j$ intersect the handles as arcs parallel to Figure~\ref{figinthandles}, one below through the favourite crossing and the other ones above.

For each $\alpha_i$, remove the annular neighborhood of $\alpha_i$ bounded by $\alpha^{\prime \prime}_i \cup (-\alpha^{\prime}_i)$ in Figure~\ref{figinthandles}  from $H_{a,2}$ in order to get the rectangular diagram of $(\CD,\matchingo,\pointw)$ of Figure~\ref{figplansplit}, Subsection~\ref{subframedHd}.

\begin{figure}[h]
\begin{center}
\begin{tikzpicture}
\useasboundingbox (0,-.2) rectangle (11.4,1.9);
\draw [fill=gray!40, draw=white] (.9,.68) rectangle (2.9,1.12);
\draw [fill=gray!40, draw=white] (5.5,.68) rectangle (7.5,1.12);
\draw [->] (1.9,1.12) -- (1.9,.68);
\draw (2.05,.5) node{\scriptsize $\alpha_1$};
\draw [>-] (1,.5) arc (-90:90:.4);
\draw (1,.2) node{\scriptsize $\alpha^{\prime}_1$};
\draw [>-] (2.8,.5) arc (270:90:.4);
\draw (2.8,.2) node{\scriptsize $\alpha^{\prime\prime}_1$} (4.2,.9) node{\scriptsize $\dots$};
\draw [->] (6.5,1.12) -- (6.5,.68);
\draw (6.65,.45) node{\scriptsize $\alpha_g$};
\draw [>-] (5.6,.5) arc (-90:90:.4);
\draw (5.6,.2) node{\scriptsize $\alpha^{\prime}_g$};
\draw [>-] (7.4,.5) arc (270:90:.4);
\draw (7.4,.2) node{\scriptsize $\alpha^{\prime\prime}_g$}; 
\draw (1,1.3) .. controls (.75,1.3) and (.75,1.12) .. (.9,1.12) -- (2.9,1.12) .. controls (3.05,1.12) and (3.05,1.3) .. (2.8,1.3);
\draw (1,.5) .. controls (.75,.5) and (.75,.68) .. (.9,.68) -- (2.9,.68) .. controls (3.05,.68) and (3.05,.5) .. (2.8,.5);
\draw (5.6,1.3) .. controls (5.35,1.3) and (5.35,1.12) .. (5.5,1.12) -- (7.5,1.12) .. controls (7.65,1.12) and (7.65,1.3) .. (7.4,1.3);
\draw (5.6,.5) .. controls (5.35,.5) and (5.35,.68) .. (5.5,.68) -- (7.5,.68) .. controls (7.65,.68) and (7.65,.5) .. (7.4,.5);
\draw (0,-.1) -- (0,1.8) -- (8.4,1.8) -- (8.4,-.1) -- (0,-.1);
\draw [->] (7,.9) -- (6.85,.9);
\draw [->] (6.15,.9) -- (6,.9);
\draw [->] (2.4,.9) -- (2.25,.9);
\draw [->] (1.55,.9) -- (1.4,.9);
\draw [->] (3.3,.75) -- (.5,.75);
\draw [->] (.5,.95) -- (3.3,.95);
\draw [->] (.5,.85) -- (3.3,.85);
\draw [->] (5.1,.75) -- (7.9,.75);
\draw [->] (5.1,.85) -- (7.9,.85);
\draw [->] (7.9,.95) -- (5.1,.95); 
\draw [->] (5.1,1.05) -- (7.9,1.05);
\draw [thick,dotted] (2.25,.9) .. controls (1.9,.9) and (2.5,1.5) .. (2.8,1.5) .. controls (4,.9) and (4.9,1.05) .. (5.1,1.05) (7.9,1.05).. controls (8.4,1.05) and (8.2,1.7) .. (2.8,1.7) .. controls (2,1.7) and (1.9,.9) .. (1.55,.9);
\draw [very thick,->] (10.3,1.2) --  (10.7,1.2);
\draw [very thick,->] (10.7,.3) --  (10.3,.3);
\draw (10.4,1.5) node[right]{\scriptsize $\ComY^{++}=\ComY^{+-}$} (11.2,.6) node[left]{\scriptsize $\ComY^{--}=\ComY^{-+}$};
\end{tikzpicture}
\caption{The punctured surface $H_{a,2}^{\matchingo}$}
\label{figsurband}
\end{center}
\end{figure}
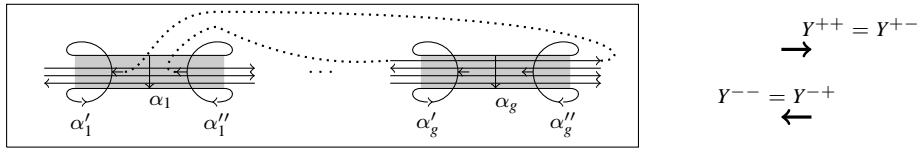

Let $H_{a,2}^{\matchingo}$ denote the complement of disk neighborhoods of the favourite crossings in the surface $H_{a,2}$.
See $H_{a,2}^{\matchingo}$ as the surface obtained from the rectangle of Figure~\ref{figplansplit} by adding a band of the handle's upper part for each $\alpha_i$, so that the band of $\alpha_i$ contains all the non-favourite crossings of $\alpha_i$. See Figure~\ref{figsurband} for an immersion of this surface in the plane. 

Extend every $\ComY=\ComY^{\varepsilon,\eta}$ on $\handleboda$
so that the fields $\ComY^{\varepsilon, \eta}$ are horizontal and their projections are the depicted constant fields in Figure~\ref{figsurband}.

\begin{figure}[h]
\begin{center}
\begin{tikzpicture}
\useasboundingbox (0,-.2) rectangle (11.4,1.9);
\draw [->] (1,.5) arc (-90:270:.4);
\draw (1,.2) node{\scriptsize $\alpha^{\prime}_1$};
\draw [->] (2.8,.5) arc (270:-90:.4);
\draw (2.8,.2) node{\scriptsize $\alpha^{\prime\prime}_1$} (4.2,.9) node{$\dots$};
\draw [->] (5.6,.5) arc (-90:270:.4);
\draw (5.6,.2) node{\scriptsize $\alpha^{\prime}_g$};
\draw [->] (7.4,.5) arc (270:-90:.4);
\draw (7.4,.2) node{\scriptsize $\alpha^{\prime\prime}_g$};
\draw (0,-.1) -- (0,1.8) -- (8.4,1.8) -- (8.4,-.1) -- (0,-.1);
\draw [->] (2.55,.9) -- (2.25,.9); \draw (2.55,.9) -- (2.8,.9);
\draw [->] (1.55,.9) -- (1.25,.9); \draw (1.25,.9) -- (1,.9);
\draw [->] (7.15,.9) -- (6.85,.9); \draw (7.15,.9) -- (7.4,.9);
\draw [->] (6.15,.9) -- (5.85,.9); \draw (5.85,.9) -- (5.6,.9);
\draw [->] (3.05,.85) -- (3.35,.85); \draw (3.05,.85) .. controls (2.95,.85) .. (2.8,.9);
\draw [->] (3,.75) -- (3.3,.75); \draw (3,.75) .. controls (2.9,.75)   .. (2.8,.9);
\draw [->] (3.25,.65) -- (2.95,.65); \draw (2.95,.65) .. controls (2.85,.65) .. (2.8,.9);
\draw [->] (.85,.65) -- (.55,.65); \draw (.85,.65) .. controls (.95,.65) .. (1,.9);
\draw [->] (.45,.85) -- (.75,.85); \draw (.75,.85) .. controls (.85,.85) .. (1,.9);
\draw [->] (.5,.75) -- (.8,.75); \draw (.8,.75) .. controls (.9,.75) .. (1,.9);
\draw [->] (5.1,.75) -- (5.4,.75); \draw (5.4,.75) .. controls (5.5,.75) .. (5.6,.9);
\draw [->] (5.05,.85) -- (5.35,.85); \draw (5.35,.85) .. controls (5.4,.85) .. (5.6,.9);
\draw [->] (5.35,.95) -- (5.05,.95); \draw (5.35,.95) .. controls (5.4,.95)  .. (5.6,.9);
\draw [->] (5.1,1.05) -- (5.4,1.05); \draw (5.4,1.05) .. controls (5.45,1.05) .. (5.6,.9);
\draw [->] (7.6,.75) -- (7.9,.75); \draw (7.6,.75) .. controls (7.5,.75)  .. (7.4,.9);
\draw [->] (7.65,.85) -- (7.95,.85); \draw (7.65,.85) .. controls (7.55,.85)  .. (7.4,.9);
\draw [->] (7.95,.95) -- (7.65,.95); \draw (7.65,.95) .. controls (7.55,.95) .. (7.4,.9);
\draw [->] (7.6,1.05) -- (7.9,1.05); \draw (7.6,1.05) .. controls (7.55,1.05)  .. (7.4,.9);
\draw [thick,dotted] (2.25,.9) .. controls (1.9,.9) and (2.5,1.5) .. (2.8,1.5) .. controls (4,.9) and (4.9,1.05) .. (5.1,1.05) (7.9,1.05).. controls (8.4,1.05) and (8.2,1.7) .. (2.8,1.7) .. controls (2,1.7) and (1.9,.9) .. (1.55,.9);
\draw [very thick,->] (10.3,1.2) --  (10.7,1.2);
\draw [very thick,->] (10.7,.3) --  (10.3,.3);
\draw (10.4,1.5) node[right]{\scriptsize $\ComY^{++}=\ComY^{+-}$} (11.2,.6) node[left]{\scriptsize $\ComY^{--}=\ComY^{-+}$};
\end{tikzpicture}
\caption{A typical slice of $[0,2g]\times [0,4] \times [-\infty,0]$}
\label{figplansplitinfty}
\end{center}
\end{figure}
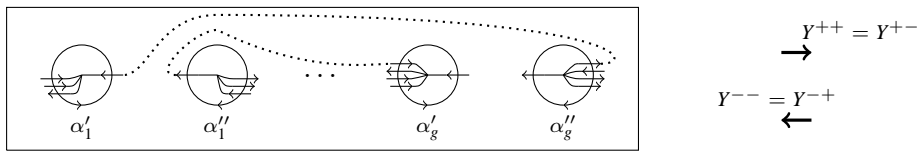

Note that $[0,2g]\times [0,4] \times [-\infty,0]$ is the product of Figure~\ref{figplansplitinfty} by $[-\infty,0]$ where all the flow lines are directed by $[-\infty,0]$.

\begin{figure}[h]
\begin{center}
\begin{tikzpicture}
\useasboundingbox (0,-.1) rectangle (10.4,4.5);
\draw (0,0) -- (0,4.4) -- (8.4,4.4) -- (8.4,0) -- (0,0);
\draw [->] (1.8,3.2) arc (0:360:.4);
\draw (1.8,3.2) node[right]{\scriptsize $\beta^{\prime}_1$};
\draw [->] (1.8,1.2) arc (360:0:.4);
\draw (1.8,1.2) node[right]{\scriptsize $\beta^{\prime\prime}_1$} (4.2,2.2) node{$\dots$};
\draw [->] (6.6,3.2) arc (-180:180:.4);
\draw (6.6,3.2) node[left]{\scriptsize $\beta^{\prime}_g$};
\draw [->] (7.4,1.2) arc (360:0:.4);
\draw (7.4,1.2) node[right]{\scriptsize $\beta^{\prime\prime}_g$};
\draw (1.4,3.2) -- (1.4,2.65) (1.4,1.75) -- (1.4,1.2);
\draw [thick,dashed,->] (1.4,2.2) -- (1.4,1.75) (1.4,2.65) -- (1.4,2.2);
\draw [->] (7,3.2) -- (7,2.65);
\draw [->] (7,1.45) -- (7,1.2) (7,1.75) -- (7,1.45);
\draw (1.18,3.8) -- (1.4,3.2) -- (1.62,3.8);
\draw (1.08,3.8) -- (1.4,3.2);
\draw (1.4,3.75) node{\tiny \dots};
\draw (1.18,.6) -- (1.4,1.2) -- (1.62,.6);
\draw (1.08,.6) -- (1.4,1.2);
\draw (1.4,.65) node{\tiny \dots};
\draw (6.78,3.8) -- (7,3.2) -- (7.22,3.8);
\draw (6.68,3.8) -- (7,3.2);
\draw (7,3.75) node{\tiny \dots};
\draw (6.78,.6) -- (7,1.2) -- (7.22,.6);
\draw (6.68,.6) -- (7,1.2);
\draw (7,.65) node{\tiny \dots};
\draw [dashed] (1.62,3.8) .. controls (1.8,4.15) and (4.7,2.2) .. (7,2.2) .. controls (8.2,2.2) and (7.4,4.15) .. (7.22,3.8);
\draw [dashed] (1.62,.6) .. controls (1.8,.25) and (7.4,.25) .. (7.22,.6);
\draw [thick,->] (10,2.8) -- node[very near end,right]{\scriptsize $\ComY^{++}=\ComY^{-+}$} (10,3.2);
\draw [thick,->] (10,2.5) -- node[very near end,right]{\scriptsize $\ComY^{--}=\ComY^{+-}$} (10,2.1);
\end{tikzpicture}
\caption{A typical slice of $[0,2g]\times [0,4] \times [6,\infty]$}
\label{figplansplitinftyb}
\end{center}
\end{figure}
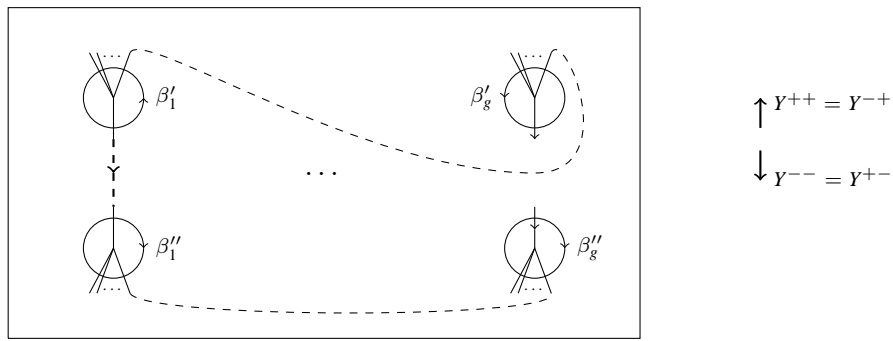

Similarly, assume that the $\alpha$-curves are orthogonal to the picture on the lower parts of the handles in the standard picture of $\handlebodb$ in Figure~\ref{figHa}, and draw a planar picture similar to Figure~\ref{figsurband} of $H_{a,4}^{\matchingo}$ (which is $\fMorse^{-1}(4) \cap C_M$ minus disk neighborhoods of the favourite crossings), by starting with Figure~\ref{figplansplitinftyb} and by adding a vertical band cut by a horizontal arc of $\beta_j$ oriented from right to left, for each $\beta_j$.

Again,
$[0,2g]\times [0,4] \times [6,\infty]$ is the product of Figure~\ref{figplansplitinftyb} by $[6,\infty]$ where all the flow lines are directed by $[6,\infty]$. Extend every $\ComY=\ComY^{\varepsilon,\eta}$ on $\handlebodb$ so that $\ComY$ looks constant and horizontal in our standard figure of $\handlebodb$ in Figure~\ref{figHa} and so that its projection on Figure~\ref{figplansplitinftyb} is the drawn constant field.

Also assume that every $\ComY=\ComY^{\varepsilon,\eta}$ varies in a quarter of horizontal plane in our tubular neighborhoods of the $\gamma_i$ in Figure~\ref{figgammai}. Similarly, extend every $\ComY=\ComY^{\varepsilon,\eta}$ in the product by $[2,4]$ of the bands of Figure~\ref{figsurband} so that $\ComY^{\varepsilon,\eta}$ is horizontal and is never a $(-\varepsilon)$-normal to the $\Asc_i$ there.

Let $H_{a,2}^{\CaC}$ denote the punctured rectangle of Figure~\ref{figplansplit}, which is a subsurface of $H_{a,2}$.
Now, $\ComY$ is defined everywhere except in $H_{a,2}^{\CaC}\times ]2,4[$
so that, for the surface $\Sigma=\Sigma(\Link(\matchingo))$ of Proposition~\ref{propSigL},\\
$$e(\ComX^{\perp}(\Sigma),\ComY)=e(\ComX^{\perp}(\Sigma \cap (H_{a,2}^{\CaC} \times [2,4])),\ComY)$$
$$\begin{array}{ll}=&-\sum_{c \in \CaC} \CJ_{j(c)i(c)}\sigma(c)e(\ComX^{\perp}([c_{j(c)},c]_{\beta} \times [2,4])),\ComY)\\
&+\sum_{(j,i) \in \{1,\dots,g\}^2}\sum_{c \in \CaC}\CJ_{j(c)i(c)}\sigma(c)\CJ_{ji} \langle \alpha_i, \halfbra c_{j(c)},c \halfbra_{\beta} \rangle e(\ComX^{\perp}(\beta_j \times [2,4]),\ComY).
\end{array}$$ 
Thus, Proposition~\ref{propcontrob} will be proved as soon as we have proved the following lemma.
\begin{lemma}
\label{lemthetaobs}
With the notation of Subsection~\ref{subframedHd},
$$\tilth(\beta_j)= -\frac14 \sum_{(\varepsilon,\eta) \in \{+,-\}^2}e(\ComX^{\perp}(\beta_j \times [2,4]),\ComY^{\varepsilon,\eta})$$
and
$$\tilth(\halfbra c_{j(c)},c \halfbra_{\beta})=- \frac14 \sum_{(\varepsilon,\eta) \in \{+,-\}^2}e(\ComX^{\perp}(\halfbra c_{j(c)},c \halfbra_{\beta}\times [2,4]),\ComY^{\varepsilon,\eta}).$$
\end{lemma}
\bp
Consider an arc $[c,d]_{\beta}$ between two consecutive crossings of $\beta$. Let $[c^{\prime},d^{\prime}]=[c,d]_{\beta} \cap H_{a,2}^{\CaC}$. On $[c^{\prime},d^{\prime}] \times [2,4]$, the field $\ComX$ is directed by $[2,4]$, the field $\ComY^{\varepsilon, \eta}$ is
defined on $\partial \left([c^{\prime},d^{\prime}] \times [2,4]\right)$, and it is in the hemisphere of the $\eta$-normal to $[c^{\prime},d^{\prime}] \times [2,4]$ along $\partial \left([c^{\prime},d^{\prime}] \times [2,4]\right) \setminus [c^{\prime},d^{\prime}] \times \{2\}$ (the $\eta$-normal is the positive normal when $\eta=+$ and the negative normal otherwise).
Then $e(\ComX^{\perp}([c^{\prime},d^{\prime}] \times [2,4]),\ComY^{\varepsilon, \eta})$ is the degree of $\ComY^{\varepsilon, \eta}$ at the $(-\eta)$-normal to $[c^{\prime},d^{\prime}]=[c^{\prime},d^{\prime}]\times \{2\}$, in the fiber of the unit tangent bundle $UH_{a,2}$ of $H_{a,2}$ trivialised by the normal to $[c^{\prime},d^{\prime}]$.
Thus, $e(\ComX^{\perp}([c^{\prime},d^{\prime}] \times [2,4]),\ComY^{\varepsilon, \eta})$ is the opposite of the degree of the $(-\eta)$-normal to $[c^{\prime},d^{\prime}]$ in the fiber of $UH_{a,2}$ at $\ComY^{\varepsilon, \eta}$ trivialised by $\ComY^{\varepsilon, \eta}$ (that is by Figure~\ref{figplansplit}) along $[c^{\prime},d^{\prime}]$.
This $(-\eta)$-normal starts and ends as vertical in this figure, and $\ComY^{\varepsilon, \eta}$ is horizontal with a direction that depends on the sign of $\varepsilon$.
The $(-\eta)$-normal to $[c^{\prime},d^{\prime}]$ makes $(\tilth(\halfbra c, d  \halfbra_{\beta}) \in \frac12 \ZZ)$ positive loops with respect to the parallelization induced by Figure~\ref{figplansplit}. Therefore the sum of the degrees of the $(-\eta)$ normal at the direction of $\ComY^{\varepsilon, \eta}$ and at the direction of $\ComY^{(-\varepsilon), \eta}$ is $2\tilth(\halfbra c, d  \halfbra_{\beta})$.

This shows that
$$\begin{array}{ll}\tilth(\halfbra c, d  \halfbra_{\beta})&=-\frac12 \left(e(\ComX^{\perp}([c^{\prime},d^{\prime}] \times [2,4]),\ComY^{\varepsilon, \eta})+e(\ComX^{\perp}([c^{\prime},d^{\prime}] \times [2,4]),\ComY^{(-\varepsilon), \eta})\right)
\\&=- \frac14 \sum_{(\varepsilon,\eta) \in \{+,-\}^2}e(\ComX^{\perp}([c^{\prime},d^{\prime}]\times [2,4]),\ComY^{\varepsilon,\eta}).\end{array}$$
The first equality of the statement follows since each side is the sum, over all the arcs of $\beta_j$ between consecutive crossings,
of the corresponding side of the equality above. The second equality follows similarly.
This concludes the proof of Lemma~\ref{lemthetaobs}, and therefore the proofs of Proposition~\ref{propcontrob} and Theorem~\ref{thmmain}. \eop


\def\cprime{$'$}


\@restonecoltrue\if@twocolumn\@restonecolfalse\fi
  \columnseprule \z@ \columnsep 35pt
 
\twocolumn
\begin{theindex}
 \addcontentsline{toc}{section}{Index of notations}

  \item $\Asc_i$, 1018
  \item $a_i$, 1009

  \indexspace

  \item $B_M$, 1006

  \indexspace

  \item $C_M$, 1017
  \item $C_2(M)$, 1003, 1005

  \indexspace

  \item $D(\alpha_i)$, 1008
  \item $\tilth$, 1015
  \item $\partial_{od}$, 1020

  \indexspace

  \item $e(\pointw,\matchingo)=e(\CD,\pointw,\matchingo)$, 1015

  \indexspace

  \item $\gamma(c)$, 1009
  \item $\Gdu$, 1031
  \item $\Gdu^b(\ComX,\ComY)$, 1030
  \item $\Gdu^i(\ComY)$, 1030

  \indexspace

  \item $\handleboda$, 1018
  \item $H_{a,2}$, 1019

  \indexspace

  \item $\iota$, 1023

  \indexspace

  \item $\CJ_{ji}$, 1010

  \indexspace

  \item $\linh(a_i)$, 1018
  \item $\linhp(a_i)$, 1018
  \item $\lambda $, 1001
  \item $\ell(.,.)$, 1014
  \item $\lambdh$, 1012, 1014
  \item $\Link(\matchingo)= \Link(\CD,\matchingo)$, 1010

  \indexspace

  \item $\matchingo$, 1009

  \indexspace

  \item $\prop(\funcf,\metrigo)$, 1020
  \item $\preprop_{\homotop}$, 1028
  \item $\preprop_{\CI}$, 1020
  \item $\preprop_{\phi}$, 1019
  \item $\preprop_{\Sigma}$, 1029
  \item $p_1$, 1001, 1007
  \item $p_{\infty}$, 1006

  \indexspace

  \item $[S]$, 1005
  \item $\sigma_a$, 1028
  \item $s_{\phi}(\check{M})$, 1019
  \item $\sigma(c)$, 1008

  \indexspace

  \item $U\check{M}$, 1005

  \indexspace

  \item $\pointw$, 1009

  \indexspace

  \item $\ComX(\pointw,\matchingo)$, 1027

\end{theindex}

\end{document}